\documentclass[12pt,a4paper,reqno]{amsart}
\usepackage{amsmath}
\usepackage{amssymb}
\usepackage{enumitem}
\usepackage{geometry}
\usepackage{mathtools,nccmath}
\usepackage{amsthm}
\usepackage{tikz}
\usetikzlibrary{shapes}
\usepackage{mathrsfs}
\usepackage{hyperref}
\usetikzlibrary{arrows}
\usepackage[utf8]{inputenc}
\numberwithin{equation}{section}
\usepackage{tikz-cd}
\usepackage[utf8]{inputenc}

\allowdisplaybreaks

\newtheorem{thm}{Theorem}[section]
\newtheorem*{theorem*}{Theorem}

\newtheorem{prop}[thm]{Proposition}
\newtheorem{lm}[thm]{Lemma}

\newtheorem{coro}[thm]{Corollary}

\newcommand\numberthis{\addtocounter{equation}{1}\tag{\theequation}}

\theoremstyle{definition}
\newtheorem{example}[thm]{Example}
\newtheorem{df}[thm]{Definition}
\newtheorem{remark}[thm]{Remark}

\newcommand{\fusion}[3]{{\binom{#3}{#1\;#2}}}

\renewcommand{\hat}[1]{\widehat{#1}}

\newcommand{\al}{\alpha}

\newcommand{\Hom}{{\rm Hom}\,}

\newcommand{\End}{{\rm End}\,}
\newcommand{\Res}{{\rm Res}\,}

\newcommand{\Z}{\mathbb{Z}}

\newcommand{\N}{\mathbb{N}}
\def\Res{{\rm Res}}
\def\wt{{\rm wt}}

\newcommand{\la}{\lambda}

\def\C{{\mathbb C}}
\def\P{{\mathbb P}}

\def\Z{{\mathbb Z}}

\def\N{{\mathbb N}}

\def\1{{\bf 1}}

\def \End{{\rm End}}
\def \Hom{{\rm Hom}}

\def \b{\beta}

\def \h{\mathfrak{h}}

\def \w{\omega}
\def \ra {\rightarrow}

\def\Cor{\mathrm{Cor}}
\def\Rad{\mathrm{Rad}}

\def\LA{\langle}
\def\RA{\rangle}
\def \ssq{\subseteq}

\def\Id{\mathrm{Id}}
\def \vac{\mathbf{1}}
\def\<{\langle}
\def\>{\rangle}

\def\ds{\dots}

\newcommand \WV {\begin{smallmatrix}W\\WV\end{smallmatrix}}

\begin {document}

\title{A proof of the fusion rules theorem}
	\author{Jianqi Liu}
\address{Department of Mathematics, University of California, Santa Cruz, CA 95064}
\email{jliu230@ucsc.edu}

\maketitle

\begin{abstract}
We prove that the space of intertwining operators associated with certain admissible modules over vertex operator algebras is isomorphic to a quotient of the vector space of conformal blocks on a three-pointed rational curve defined by the same data. This provides a new proof and alternative version of Frenkel and Zhu’s fusion rules theorem, in terms of the dimension of certain bimodules over Zhu’s algebra, without the assumption of rationality.
\end{abstract}

\section{Introduction}
The space of intertwining operators of vertex operator algebras (see \cite{B,FLM,FHL}) and its dimension, the so-called fusion rule in the physics literature \cite{MS,TK,TUY}, plays an essential role in studying the tensor product of modules over vertex operator algebras. In the semi-simple case, the fusion rule is the multiplicity of an irreducible module in a tensor product. For the affine Lie algebras or the associated affine vertex operator algebras \cite{FZ}, the fusion rules in case $\widehat{sl_{2}(\C)}$ were computed in \cite{TK}, and a general version was stated in \cite{TUY} without proof. 
In \cite{FZ}, Frenkel and Zhu proposed a formula (Theorem 1.5.2 in \cite{FZ}) to compute the fusion rules for arbitrary vertex operator algebras by using Zhu's algebra $A(V)$ defined in \cite{Z} and some of its (bi)modules. Given irreducible modules $M^{1},M^{2}$ and $M^3$ over a vertex operator algebra $V$, Frenkel and Zhu's fusion rules theorem claimed that the space of intertwining operators $I\fusion{M^1}{M^2}{M^3}$ can be identified with the vector space $(M^{3}(0)^{*}\otimes _{A(V)}A(M^{1})\otimes_{A(V)}M^{2}(0))^{*},$ where $A(M^{1})$ is a bimodule over the Zhu's algebra $A(V)$, and $M^{2}(0)$ and $M^{3}(0)$ are the bottom levels of the $V$-modules $M^{2}$ and $M^{3}$, which are modules over $A(V)$, see Section $1$ in \cite{FZ} for more details.

However, it was later realized by Li (see \cite{L}) that some additional conditions are needed in Frenkel and Zhu's fusion rules theorem. Li gave a counter-example in \cite{L} in the case of the universal Virasoro vertex operator algebra that shows that $I\fusion{M^1}{M^2}{M^3}$  is not isomorphic to $(M^{3}(0)^{*}\otimes _{A(V)}A(M^{1})\otimes_{A(V)}M^{2}(0))^{*}$ in general. Li also proposed in \cite{L} that the fusion rules theorem is true when $M^2$ and $M^3$ are the so-called generalized Verma modules constructed in \cite{DLM1}. In particular, it is true for the rational vertex operator algebras (see Section 2 in \cite{L} for more detailed discussions and the counter-example). 

In this paper, we give an alternative version of the fusion rules theorem for general vertex operator algebras. It can be stated as follows: 

\begin{thm}\label{thm1.1}
	Let $V$ be a CFT-type vertex operator algebra, and let $M^{1}$, $M^{2}$, and $M^{3}$ be $V$-modules with conformal weights $h_{1}$, $h_{2}$, and $h_{3}$, respectively. Assume $M^{2}(0)$ and $M^{3}(0)$ are irreducible $A(V)$-modules, then we have the following isomorphism of vector spaces:
$$I\fusion{M^1}{\bar{M}(M^{2}(0))}{\bar{M}(M^{3}(0)^{\ast})'}\cong	I\fusion{M^{1}}{\bar{M^2}}{\bar{M^3}}\cong (M^{3}(0)^{*}\otimes _{A(V)}B_{h}(M^{1})\otimes_{A(V)}M^{2}(0))^{*},$$
	where $h=h_{1}+h_{2}-h_{3}$, and $\bar{M^2}=\bar{M}/\Rad(\bar{M})$ and $\bar{M^3}'=\tilde{M}/\Rad{\tilde{M}}$ are quotient modules of the generalized Verma modules $\bar{M}(M^{2}(0))$ and $\bar{M}(M^{3}(0)^{\ast})$, respectively.
\end{thm}
In our version of the fusion rules theorem, we replaced the $A(V)$-bimodule $A(M^{1})$ by a newly defined $A(V)$-bimodule $B_{h}(M^{1})$, which is given by $B_{h}(M^{1})=M^{1}/span\{a\circ u,\ L(-1)v+(L(0)+h_2-h_3)v:a\in V,\ u,v\in M^{1}\}$. We will show that $B_h(M^1)$ is a quotient module of $A(M^1)$, and we will give examples to show that the vector spaces $(M^{3}(0)^{*}\otimes _{A(V)}B_{h}(M^{1})\otimes_{A(V)}M^{2}(0))^{*}$ and $(M^{3}(0)^{*}\otimes _{A(V)}A(M^{1})\otimes_{A(V)}M^{2}(0))^{*}$ are not isomorphic in general. We need to mod out the additional terms $L(-1)v+(L(0)v+h_{2}-h_{3})v$ in $A(M^{1})$ because otherwise, the $L(-1)$-derivation property of the intertwining operators cannot be correctly reflected. We will also give sufficient conditions for modules $\bar{M^2}$ and $\bar{M^3}$ to be irreducible. In particular, for a CFT-type rational vertex operator algebra $V$, the modules $\bar{M^2}$ and $\bar{M^3}$ are automatically irreducible, then the fusion rule $\dim I\fusion{M^1}{M^2}{M^3}$ for three irreducible $V$-modules is equal to the dimension of $(M^{3}(0)^{*}\otimes _{A(V)}B_{h}(M^{1})\otimes_{A(V)}M^{2}(0))^{*}$.

Our proof of Theorem \ref{thm1.1} is different than Li’s proof of Theorem 2.11 in \cite{L}. We prove Theorem \ref{thm1.1} based on a combination of ideas the ideas from \cite{TUY} and extensions made in \cite{Z0}, wherein a system of correlation functions is associated with every vector in the space of conformal blocks (see Theorem 6.2 in \cite{Z0}). Based on the properties of the following prototype system of $(n+3)$-point correlation functions on $\P^1_\C$:
\begin{equation}\label{1.2}(v_{3}', Y_{M^{3}}(a_{1},z_{1})\ds Y_{M^{3}}(a_{k},z_{k})I(v,w)Y_{M^{2}}(a_{k+1},z_{k+1})\ds Y_{M^2}(a_{n},z_{n})v_{2}),
\end{equation}
where $v_{3}'\in M^{3}(0)^{\ast}$, $v\in M^{1}$, $v_{2}\in M^{2}$, $a_{1},\ds,a_{n}\in V$, and $I$ is an intertwining operator of type $\fusion{M^1}{M^2}{M^3}$, we introduce the notion of space of correlation functions associated with $V$-modules $M^{1}$, $M^{2}$, and $M^{3}$, denoted by $\Cor\fusion{M^1}{M^2}{M^3}$. It is essentially a quotient of the vector space of three-point genus zero conformal blocks, the dual space to a certain quotient of the tensor product of $3$ admissible $V$-modules (see \cite{TUY,Z0}). Then we prove that $\Cor\fusion{M^1}{M^2}{M^3}$ is isomorphic to $I\fusion{M^1}{M^2}{M^3}$.

 In order to relate $\Cor\fusion{M^1}{M^2}{M^3}$ with the modules over $A(V)$, we introduce an auxiliary notion of the space of correlation functions associated with $M^1$, $M^{2}(0)$, and $M^{3}(0)$, denoted by $\Cor \fusion{M^1}{M^{2}(0)}{M^{3}(0)}$. This space can be viewed as the space $A(V)$-conformal blocks on the $3$-pointed rational curve $\P^1_\C$ defined from the representations of Zhu's algebra $A(V)$. The axioms we imposed on this space are based on the restriction of \eqref{1.2} onto the bottom levels $M^{2}(0)$ and $M^{3}(0)^{\ast}$. Then we use certain generating formulas satisfied by the correlation function \eqref{1.2} and prove that $\Cor\fusion{M^1}{M^{2}(0)}{M^{3}(0)}$ is isomorphic to both $\Cor \fusion{M^1}{\bar{M^2}}{\bar{M^3}} $ and $\Cor\fusion{M^1}{\bar{M}(M^{2}(0))}{\bar{M}(M^{3}(0)^{\ast})'}$ when $M^2(0)$ and $M^3(0)$ are irreducible modules over $A(V)$. However, unlike building $V$-modules from $A(V)$-modules (see Theorem 2.2.1 in \cite{Z}) based on the ordinary correlation functions $(v',Y(a_1,z_1)\dots Y(a_n,z_n)v)$, in our case, due to the appearance of intertwining operator $I(v,w)$ in \eqref{1.2}, the modules $\bar{M^2}$ and $\bar{M^3}$ constructed by \eqref{1.2} are not necessarily irreducible. This issue was first observed by Li in \cite{L}. The $V$-modules $\bar{M^2}$ and $\bar{M^3}$ are quotient modules of certain generalized Verma modules. They can be proved to be irreducible if a technical condition depends only on the (bi)modules over $A(V)$ is satisfied. 
 
 We then prove that $\Cor\fusion{M^1}{M^2(0)}{M^3(0)}$ is isomorphic to $(M^{3}(0)^{*}\otimes _{A(V)}B_{h}(M^{1})\otimes_{A(V)}M^{2}(0))^{*}$. Given a linear function $f$ on $M^{3}(0)^{*}\otimes _{A(V)}B_{h}(M^{1})\otimes_{A(V)}M^{2}(0)$, we shall use the recursive formulas satisfied by \eqref{1.2} and reconstruct a system of correlation functions in $\Cor\fusion{M^1}{M^2(0)}{M^3(0)}$. There is one recursive formula ((2.2.1) in \cite{Z}) of the correlation functions $S(v',(a_1,z_1)\dots (a_n,z_n)v)=(v',Y(a_1,z_1)\dots Y(a_n,z_n)v)$, where $v\in M(0)$ and $v'\in M(0)^\ast$, obtained by expanding the left-most term $Y(a_1,z_1)$. However, in our case, this formula alone is not enough to rebuild the correlation functions from $f$. The reason is again because of the appearance of $I(v,w)$ in the correlation functions, which makes expanding the left-most term $(v,w)$ in $S(v_{3}', (v,w)(a_1,z_1)\dots (a_n,z_n)v_2)$ unreasonable, as the action $v(n)a_i=\Res_{z} w^{n+h}I(v,w)a_i $ is not yet defined. We remedy this situation by introducing an additional recursive formula for the correlation functions \eqref{1.2} obtained by expanding the right-most term $Y(a_n,z_n)$ in $(v_{3}',I(v,w)Y(a_{1},z_{1})\ds Y(a_n,z_n)v_{2})$, where $v_{3}'\in M^3(0)^\ast$ and $v_{2}\in M^2(0)$, and we use both the recursive formulas to reconstruct the correlation functions from $f$. Then Theorem \ref{thm1.1} follows from the isomorphisms
 $I\fusion{M^{1}}{\bar{M^2}}{\bar{M^3}}\cong \Cor\fusion{M^1}{\bar{M^2}}{\bar{M^3}}\cong \Cor\fusion{M^1}{M^{2}(0)}{M^{3}(0)}\cong (M^{3}(0)^{*}\otimes _{A(V)}B_{h}(M^{1})\otimes_{A(V)}M^{2}(0))^{*}.$
 
 This paper is organized as follows: In Section $2$, we define $\Cor  \fusion{M^1}{M^2}{M^3}$ and prove that it is isomorphic to $I \fusion{M^1}{M^2}{M^3}$. In Section $3$, we define $\Cor\fusion{M^1}{M^{2}(0)}{M^{3}(0)}$ for irreducible $A(V)$-modules $M^2(0)$ and $M^3(0)$ and prove that $\Cor\fusion{M^1}{M^{2}(0)}{M^{3}(0)}$ is isomorphic to both $\Cor \fusion{M^1}{\bar{M^2}}{\bar{M^3}} $ and $\Cor\fusion{M^1}{\bar{M}(M^{2}(0))}{\bar{M}(M^{3}(0)^{\ast})'}$.  In section $4$, we define the $A(V)$-bimodule $B_{h}(M^1)$ and prove that $\Cor\fusion{M^1}{M^{2}(0)}{M^{3}(0)}$ is isomorphic to $ (M^{3}(0)^{*}\otimes _{A(V)}B_{h}(M^{1})\otimes_{A(V)}M^{2}(0))^{*},$ which finishes the proof of Theorem \ref{thm1.1}.  Then we verify this theorem on some particular examples, one of which shows that the counter-example given by Li in \cite{L} does not contradict Theorem \ref{thm1.1}. 
 
 We expect the readers are familiar with the concept of vertex operator algebras, modules over vertex operator algebras, and the $A(V)$-theory, see \cite{B,FHL,Z}. 



 \section{The space of correlation functions associated with $M^{1},M^{2}$, and $M^{3}$}
 We fix some notations that will be in force throughout this paper. We denote by $\C,\Z,$ and $\N$ the set of complex numbers, the set of integers, and the set of natural numbers, including $0$. All vector spaces are defined over $\C$.

Let $V=(V,Y,\vac,\w)$ be a vertex operator algebra (VOA) which is of the CFT-type: $V=\bigoplus_{n=0}^\infty V_n$, with $V_{0}=\C \vac$. A module $M$ over $V$ is an ordinary $V$-module: $M=\bigoplus_{n=0}^{\infty} M_{\la+n}$, where each $M_{\la+n}$ is an eigenspace of $L(0)$ with eigenvalue $\la+n$. Any $V$-module $M$ is $\N$-gradable (or admissible): $M=\bigoplus_{n=0}^{\infty}M(n)$, with $M(n)=M_{\la+n}$ for each $n$. We 
write $Y_{M}(a,z)=\sum_{n\in \Z} a(n)z^{-n-1},$  for all $a\in V$, and we write $  Y_{M}(\w,z)=\sum_{n\in \Z}L(n)z^{-n-2}.$ One can find more details about the definitions in \cite {DLM1,FHL,FLM,Z}.

When we use the integral sign $\int_{C} f(z)dz$, where $C$ is a simple closed contour of $z$, it means $\frac{1}{2 \pi i}\int_{C}f(z)dz$.

\subsection{The $(n+3)$-point correlation functions}  
Let $M^{1},M^{2}$, and $M^{3}$ be $V$-modules with conformal weights $h_{1},h_{2}$, and $h_{3}$,  respectively, and let $I\in I\fusion{M^{1}}{M^{2}}{M^{3}}$ be an intertwining operator. Recall that 
$I(v,w)=\sum_{n\in \Z} v(n)w^{-n-1}\cdot w^{-h},$
where $h=h_{1}+h_{2}-h_{3}$, and $v(n)=\Res_{w}I(v,w)w^{n+h}$. Moreover, $v(n)M^{2}(m)\subseteq M^{3}(\deg v-n-1+m)$ for all $n\in \Z$ and $m\in \N$, see \cite{FZ} for more details. Consider the power series 
\begin{equation}\label{n2.1}
\langle v_{3}', Y(a_{1},z_{1})\dots I(v,w)\ds Y(a_{n},z_{n})v_{2}\rangle w^{h}
\end{equation}
in $n+1$ complex variables $z_{1},\ds ,z_{n}, w$ with integer powers, where $a_{1},\ds ,a_{n}\in V$, $v\in M^{1}$, $v_{2}\in M^{2}$, and $v_{3}'\in M^{3'}$ which is the contragredient module of $M^{3}$ (cf. \cite{FHL}). We multiply the term $w^h$ to avoid the appearance of the logarithm when computing the integrations.

Recall that the power series \eqref{n2.1} converges in the domain $$\mathbb{D}=\{(z_{1},\dots,z_{n},w)\in \mathbb{C}^{n+1}||z_{1}|>|z_{2}|>\dots>|w|>\dots>|z_{n}|>0\}$$
to a rational function in $z_{1},...,z_{n}, w,\ z_{i}-z_{j}$ and $z_{k}-w$, where $1\leq i\neq j\leq n$ and $1\leq k\leq n$. We denote this rational function by: 
\begin{equation}\label{n2.2}
(v_{3}', Y(a_{1},z_{1})\ds I(v,w)\ds Y(a_{n},z_{n})v_{2}),
\end{equation}
also recall that the only possible poles of \eqref{n2.2} are at $z_{i}=0,\ w=0,\ z_{i}=z_{j}$ and $z_{k}=w$, see \cite{FHL} for more details. 

Moreover, it is also essentially proved in \cite{FHL} that the rational function \eqref{n2.2} is invariant under the permutation of the terms $Y(a_{1},z_{1}),\ds,Y(a_{n},z_{n})$, and $I(v,w)$. In other words, the power series \eqref{n2.1}
and the power series 
$\<v_{3}',Y(a_{i_{1}},z_{i_{1}})\ds I(v,w)\ds Y(a_{i_{n}},z_{i_{n}})v_{2}\>w^{h}$
have the same limit function \eqref{n2.2} on their corresponding domain of convergence. 
 
We use the symbol $S_I$ as in \cite{Z} to denote the limit function \eqref{n2.2}: 
 \begin{equation}\label{n2.3}
 S_I(v_{3}',(a_{1},z_{1})\ds(v,w)\ds(a_{n},z_{n})v_{2}):=(v_{3}', Y(a_{1},z_{1})\ds I(v,w)\ds Y(a_{n},z_{n})v_{2}).
 \end{equation}
 Then we have a system of linear maps $S_I=\{(S_I)_{V...M^1...V}^n\}_{n=0}^\infty$: 
  \begin{equation}\label{n2.4}
  \begin{aligned}
  (S_I)_{V...M^1...V}^n: M^{3'}\times V\times\ldots\times M^{1}&\times\ds V\times M^{2}\rightarrow \mathcal{F}(z_{1},
  \ds,z_{n},w),\\
  (v_{3}',a_{1},\ds,v,\ds,a_{n},v_{2})&\mapsto S_I(v_{3}',(a_{1},z_{1})\ds(v,w)\ds(a_{n},z_{n})v_{2}),
  \end{aligned}
  \end{equation} 
  where $\mathcal{F}(z_{1},\ds,z_{n},w)$ is the space of rational functions in $n+1$ variables $z_{1},z_{2},\ds,z_{n}, w$, with only possible poles at $z_{i}=0,\ w=0,\ z_{i}=z_{j},\ z_{k}=w$. For a fixed $n\in \N$, we have $(S_I)^n_{M^1V...V}=(S_I)^n_{VM^1...V}=\dots =(S_I)^n_{V...VM^1}$, since the terms $(a_1,z_1),\dots ,(a_n,z_n)$, and $(v,w)$ can be permuted within $S_I$ in \eqref{n2.3}. 
  
  We introduce the following notion that generalizes Definition 4.1.1 in \cite{Z}: 
  
\begin{df}\label{df2.1}
A system of linear maps $S=\{S^n_{V...M^1...V}\}_{n=0}^\infty$, 
\begin{align*}
  S_{V...M^1...V}^n: M^{3'}\times V\times\ldots\times M^{1}&\times\ds V\times M^{2}\rightarrow \mathcal{F}(z_{1},
\ds,z_{n},w),\\
(v_{3}',a_{1},\ds,v,\ds,a_{n},v_{2})&\mapsto S(v_{3}',(a_{1},z_{1})\ds(v,w)\ds(a_{n},z_{n})v_{2}),
\end{align*}
is said to satisfy the {\em genus-zero property associated with $M^1,M^2$, and $M^3$} if it satisfies
  	\begin{enumerate}
  		\item (The truncation property) For fixed $v\in M^{1}$ and $v_{2}\in M^{2}$, the Laurent series expansion of  $S(v_{3}',(v,w)v_{2})$ around $w=0$ has a uniform lower bound for $w$ independent of $v_{3}'\in M^{3'}$. i.e., $S(v_{3}',(v,w)v_2)=\sum_{n\leq N} a_n w^{-n-1}$ for all $v_{3}'\in M^{3'}$. 
  		\item (The locality) The terms $(a_1,z_1),\dots ,(a_n,z_n)$, and $(v,w)$ can be permuted arbitrarily within $S$. i.e., $S^n_{M^1V...V}=S^n_{VM^1...V}=\dots =S^n_{V...VM^1}$ for any fixed $n\in \N$.
  		\item (The vacuum property)
  	\begin{equation}\label{n2.5}	S(v_{3}',(\vac, z)(a_{1},z_{1})\ds(v,w)\ds(a_{n},z_{n})v_{2})=S(v_{3}',(a_{1},z_{1})\ds(v,w)\ds(a_{n},z_{n})v_{2}).
  	\end{equation}
  	\item (The $L(-1)$-derivation property)
  	\begin{equation}\label{n2.6}
\begin{aligned}
&S(v_{3}',(L(-1)a_{1},z_{1})\ds (a_n,z_n)(v,w)v_{2})=\frac{d}{dz_{1}}S(v_{3}',(a_{1},z_{1})\ds(a_n,z_n)(v,w)v_{2}),\\
&S(v_{3}',(L(-1)v,w)(a_{1},z_{1})\ds v_{2})w^{-h}=\frac{d}{dw}\left(S(v_{3}',(v,w)(a_{1},z_{1})\ds v_{2})w^{-h}\right).
\end{aligned}
\end{equation}
  		\item (The associativity)

  		\begin{align}\label{n2.7}
  		&\int_{C}S(v_{3}',(a_{1},z_{1})(v,w)\ds(a_{n},z_{n})v_{2})(z_{1}-w)^{k} d{z_{1}}=S(v_{3}',(a_{1}(k)v,w)\ds(a_{n},z_{n})v_{2}), \nonumber \\
  		&\int_{C}S(v_{3}',(a_{1},z_{1})(a_{2},z_{2})\ds(v,w)v_{2})(z_{1}-z_{2})^{k} d{z_{1}}=S(v_{3}',(a_{1}(k)a_{2},z_{2})\ds(v,w)v_{2}),
\end{align}
  		where in the first equation of \eqref{n2.7}, $C$ is a contour of $z_{1}$ surrounding $w$, with $z_{2},...,z_{n}$ outside of $C$; while in the second equation of \eqref{n2.7}, $C$ is a contour of $z_{1}$ surrounding $z_{2}$, with $z_{3},...,z_{n}, w$ outside of $C$.
  		\item (The Virasoro relation) 
  		Let $\w\in V$ be the Virasoro element, and let $x,x_{1},\ds,x_{m}$ be complex variables, denote the rational function
  		$$S(v_{3}',(\w,x_{1})\ds(\w,x_{m})(a_{1},z_{1})\ds(v,w)\ds(a_{n},z_{n})v_{2})$$
  		by $S$ for simplicity. Assume that $v_{3}',v, v_{2}, a_{1},\ds,a_{n}$ are highest weight vectors for the Virasoro algebra, then we have: 

  		\begin{align*}
  		&S(v_{3}',(\w,x)(\w,x_{1})\ds(\w,x_{m})(a_{1},z_{1})\ds(v,w)\ds(a_{n},z_{n})v_{2})\\
  		&=\sum_{k=1}^{n}\frac{x^{-1}z_{k}}{x-z_{k}}\frac{d}{dz_{k}}S+\sum_{k=1}^{n}\frac{\wt a_{k}}{(x-z_{k})^{2}}S+\frac{x^{-1}w}{x-w}w^{h}\frac{d}{dw}(S\cdot w^{-h})+\frac{\wt v}{(x-w)^{2}}S\\
  		&\ \ +\frac{\wt v_{2}}{x^{2}}S+\sum_{k=1}^{m}\frac{x^{-1}w_{k}}{x-x_{k}} \frac{d}{dx_{k}}S+\sum_{k=1}^{m}\frac{2}{(x-x_{k})^{2}}S\numberthis\label{n2.8}\\
  		&\ \ +\frac{c}{2}\sum_{k=1}^{m}\frac{1}{(x-x_{k})^{4}} S(v_{3}',(\w,x_{1})\ds\widehat{(\w,x_{k})}\ds(\w,x_{m})(a_{1},z_{1})\ds(v,w)\ds(a_{n},z_{n})v_{2})
  		\end{align*}

  		\item (The generating property for $M^2$) 
  		For any $a\in V$ and $m\in \Z$, we have:
  		\begin{equation}\label{n2.9}
  		\begin{aligned}
  		&S(v_{3}',(a_{1},z_{1})\ds(v,w)\ds(a_{n},z_{n})a(m)v_{2})\\
  		=&\int_{C} S(v_{3}',(a_{1},z_{1})\ds(v,w)\ds(a_{n},z_{n})(a,z)v_{2})z^{m}dz,
  		\end{aligned}
  		\end{equation}
  			where $C=C_{R}(0)$ is a contour of $z$ surrounding $0$ with $z_{1},\ds,z_{n},w$ lying outside.
  		\item (The generating property for $M^{3'}$)
  	 Denote $(e^{z^{-1}L(1)}(-z^{2})^{L(0)}a, z)$ by $(a,z)'$, then 
  		\begin{equation}\label{n2.10}
  		\begin{aligned}
  		&S(a(m)v_{3}', (a_{1},z_{1})\ldots(v,w)\ds(a_{n},z_{n})v_{2})\\
  		&=\int_{C'} S(v_{3}', (a,z)'(a_{1},z_{1})\ds(v,w)\ds(a_{n},z_{n})v_{2})z^{-m-2}dz,
  		\end{aligned}
  		\end{equation}
  		where $C'=C_{r}(0)$ is a contour of $z$ surrounding $0$ with $z_{1},\ldots,z_{n},w$ lying inside. 
  		\end{enumerate}
\end{df}

\begin{df}\label{df2.2}
	The vector space of the system of linear maps $S=\{S^n_{V...M^1...V}\}_{n=0}^\infty$ satisfying the genus-zero property associated with $M^1,M^2$, and $M^{3}$ is called the {\em space of correlation functions associated with $M^1,M^2$, and $M^{3}$}. We denote it by $\Cor\fusion{M^1}{M^2}{M^3}$.
\end{df}

\begin{prop}\label{prop2.3}
The system of functions $S_I$ given by \eqref{n2.3} and \eqref{n2.4} satisfies the genus-zero property associated with $M^1,M^2$, and $M^3$ in Definition \ref{df2.1}. Thus $S_I\in \Cor\fusion{M^1}{M^2}{M^3}$. 
\end{prop}
\begin{proof}
The properties (1) - (6) for $S_{I}$ follow immediately from the axioms satisfied by the intertwining operator $I$ and the vertex operator $Y$; see Section 5.6 in \cite{FHL} for more details.

To prove \eqref{n2.9}, we note that the Laurent series expansion of the rational function \eqref{n2.3} on the domain $|z|<|z_{i}|,|w|$ for all $i$ is 
$\sum_{m\in \Z} (v_{3}', Y(a_{1},z_{1})\ds I(v,w)\ds a(m)v_{2})z^{-m-1}.$
The coefficient of $z^{-m-1}$ in the Laurent series is also 
$$\int_{C} (v_{3}', Y(a_{1},z_{1})\ds I(v,w)\ds Y(a_{n},z_{n})Y(a,z)v_{2})z^{m}dz,$$
where $C=C_{R}(0)$ is a contour of $z$ surrounding $0$ with $z_{1},\ds ,z_{n}$ and $w$ lying outside. This proves \eqref{n2.9}. To prove \eqref{n2.10}, we denote the term $\sum_{j\geq 0}\frac{1}{j!}(-1)^{\wt a}(L(1)a^{j})(2\wt a-m-j-2)$ by $a'(m)$, then by the definition of contragredient module (see (5.2.4) in \cite{FHL}), the series
\begin{align*}
&\sum_{m\in \Z }(a(m)v_{3}', Y(a_{1},z_{1})\ds I(v,w)\ds Y(a_{n},z_{n})Y(a,z)v_{2})z^{-m-1}\\
&=\sum_{m\in \Z}(v_{3}', a'(m)Y(a_{1},z_{1})\ds I(v,w)\ds Y(a_{n},z_{n})v_{2})z^{-m-1}
\end{align*}
is the expansion of $(v_{3}', Y(e^{zL(1)}(-z^{-2})^{L(0)}a,z^{-1})Y(a_{1},z_{1})\ds I(v,w)\ds Y(a_{n},z_{n})v_{2})$
on the domain $|z^{-1}|>|z_{i}|, |w|$, or equivalently, $|z|<1/|z_{i}|, 1/|w|$, for $i=1,\ds,n$. By comparing the Laurent coefficient of $z^{-m-1}$, we have:
\begin{align*}
&	(a(m)v_{3}', Y(a_{1},z_{1})\ds I(v,w)\ds Y(a_{n},z_{n})Y(a,z)v_{2})\\
&=\int_{C_{R}(0)}(v_{3}', Y(e^{zL(1)}(-z^{-2})^{L(0)}a,z^{-1})Y(a_{1},z_{1})\ds I(v,w)\ds Y(a_{n},z_{n})v_{2}) z^{m}dz, \numberthis\label{n2.11} 
\end{align*}
where $R$ is small enough such that $R<1/|z_{i}|, 1/|w|$, for $i=1,\ds,n$. Change the variable $z\rightarrow 1/z$ in the integral \eqref{n2.11}. Note that the parametrization of $1/z$ is $(1/R)e^{-i\theta}$,  which gives us a clockwise orientation, and $d(1/z)=-(1/z^{2})dz$. Let $C'=C_{r}(0)$, with radius $r=1/R>|z_{i}|,|w|$ for $i=1,\ds,n$, equipped with the counterclockwise orientation. Then $z_{1},\ds ,z_{n},w$ are inside of $C'$, and 
\begin{align*}
\eqref{n2.11}&=-\int_{C'} (v_{3}', Y(e^{z^{-1}L(1)}(-z^{2})^{L(0)}a,z)Y(a_{1},z_{1})\ds I(v,w)\ds Y(a_{n},z_{n})v_{2})z^{-m}(-z^{-2})dz\\
&=\int_{C'} (v_{3}', Y(e^{z^{-1}L(1)}(-z^{2})^{L(0)}a,z)Y(a_{1},z_{1})\ds I(v,w)\ds Y(a_{n},z_{n})v_{2})z^{-m-2}dz\\
&=\int_{C'} S_I(v_{3}',(a,z)' (a_{1},z_{1})\ds (v,w)\ds (a_{n},z_{n})v_{2} )z^{-m-2}dz.
\end{align*}
This proves \eqref{n2.10}. 
\end{proof}

\begin{remark}\label{rmk3.2} Let $S\in \Cor\fusion{M^1}{M^2}{M^3}$. With the notations of Proposition \ref{prop2.3}, we have: 
\begin{align*}
&S(a'(m)v_{3}',(a_{1},z_{1})\ds (v,w)\ds (a_{n},z_{n})v_{2})\\
&=\sum_{j\geq 0}\frac{1}{j!}(-1)^{\wt a}\int_{C'}S( v_{3}',(e^{z^{-1}L(1)}(-z^{2})^{L(0)}(L(1)^{j}a),z)(a_{1},z_{1})\ds v_{2}) z^{-2\wt a+m+j}dz\\
&=\int_{C'}S( v_{3}',  (e^{z^{-1}L(1)}(-z^{2})^{L(0)}e^{zL(1)}(-z^{-2})^{L(0)}a,z)(a_{1},z_{1})\ds v_{2}) z^{m}dz\\
&=\int_{C'}S(v_{3}', (e^{z^{-1}L(1)}e^{-z^{-1}L(1)}a,z)(a_{1},z_{1})\ds v_{2}) z^{m}dz\\
&=\int_{C'} S(v_{3}', (a,z)(a_{1},z_{1})\ds (v,w)\ds (a_{n},z_{n})v_{2})z^{m}dz.
\end{align*}
Hence the generating property for $M^{3'}$ (8) in Definition \ref{df2.1} is equivalent to:

\begin{align}\label{n2.12}
&	S(a'(m)v_{3}',(a_{1},z_{1})\ds (v,w)\ds (a_{n},z_{n})v_{2})\\
&=\int_{C'} S(v_{3}', (a,z)(a_{1},z_{1})\ds (v,w)\ds (a_{n},z_{n})v_{2})z^{m}dz,\nonumber
\end{align}
where $a'(m)=\sum_{j\geq 0}\frac{1}{j!}(-1)^{\wt a}(L(1)^{j}a)(2\wt a-m-j-2)$ and $C'=C_{r}(0)$ as in (8). 
\end{remark}

As a consequence of Proposition \ref{prop2.3}, we have a well-defined linear map:
\begin{equation}\label{n2.13}
\al: I\fusion{M^{1}}{M^{2}}{M^{3}}\ra \Cor\fusion{M^{1}}{M^{2}}{M^{3}},\qquad  I\mapsto S_I,
\end{equation}
where $S_I$ is given by \eqref{n2.3} and \eqref{n2.4}.

\subsection{The space of correlation functions and the space of intertwining operators}
Although the genus-zero property associated with three $V$-modules in Definition \ref{df2.1} seems long and intrinsic, it is good enough to characterize an intertwining operator. In other words, we can construct an inverse of the map $\al$ in \eqref{n2.13}. 

Fix a system of correlation functions $S$ in $ \Cor\fusion{M^{1}}{M^{2}}{M^{3}}$, we construct an intertwining operator $I_S\in I\begin{psmallmatrix}
M^{3}\\M^{1}\ M^{2} \end{psmallmatrix}$ in the following way:

Let $v\in M^{1}$, define a linear map $v(n):M^{2}\rightarrow M^{3}$ by the formula:
\begin{equation}\label{n2.14}
\langle v_{3}',v(n)v_{2}\rangle:= \int_{C}S(v_{3}',(v,w)v_{2})w^{n}dw,
\end{equation}
where $C$ is a contour of $w$ surrounding $0$. Note that an element $u\in M^{3}$ is uniquely determined by the value $\LA v_{3}',u \RA$ for $v_{3}'\in M^{3'}$, so we have a well-defined element $v(n)v_2$ in $M^{3}$.  Define $I(v,w)$ by 
\begin{equation}\label{n2.15}
I_S(v,w):=\sum_{n\in \mathbb{Z}}v(n)w^{-n-1}\cdot w^{-h},
\end{equation}
where $h=h_{1}+h_{2}-h_{3}$. Then $I(v,w)\in \Hom(M^2,M^3)\{z\}$.

\begin{thm}\label{thm2.5}
The series $I_S(v,w)$ defined by \eqref{n2.14} and \eqref{n2.15} is an intertwining operator of type $\fusion{M^{1}}{M^{2}}{M^{3}}$.
\end{thm}
\begin{proof}
By Definition \ref{df2.1}, $S(v_{3}',(v,w)v_{2})$ is a rational function in $w$ with the only possible pole at $w=0$, and the term \eqref{n2.14} is the Laurent coefficient of $S(v_{3}',(v,w)v_{2})$. Thus the series $\langle x_{3}', I_S(v,w)x_{2}\rangle w^{h}$ is the Laurent series expansion of $S(x_{3}', (v,w)x_{2})$ around $w=0$ by \eqref{n2.15}. In particular, if we denote the limit of the Laurent series $\langle v_{3}', I(v,w)v_{2}\rangle w^{h}$ by $(v_{3}', I(v,w)v_{2})$, then we have the following equality of rational functions:
\begin{equation}\label{n2.16}
(v_{3}', I_S(v,w)v_{2})=S(v_{3}',(v,w)v_{2})
\end{equation}

Since $S$ satisfies the property (1) in Definition \ref{df2.1}, for $v\in M^{1}$ and $v_{2}\in M^{2}$, there exists $N\in \Z$ such that 
$\<v_{3}',I_S(v,w)v_{2}\>w^{h}=\sum_{n\leq N}\left(\int_{C}S(v_{3}',(v,w)v_{2})w^{n}dw\right)w^{-n-1},$ for all $v_{3}'\in M^{3'}$. Hence we have $v(n)v_{2}=0$ for $n\gg0$. 
By the locality of $S$, together with \eqref{n2.15}, we have: 
$$
\<v_{3}', I_S(L(-1)v,w)v_{2}\>=\frac{d}{dw}(S(v_{3}',(v,w)v_{2})w^{-h})=\frac{d}{dw}\<v_{3}', I_S(v,w)v_{2}\>.
$$
Hence $I_S(L(-1)v,w)=\frac{d}{dw}I_S(v,w)$. Moreover, we claim that the following equation holds:
\begin{align}\label{n2.17}
&\sum_{i=0}^{\infty}\binom{m}{i}(a(l+i)v)(m+n-i)v_{2}\\
=&\sum_{i=0}^{\infty}(-1)^{i}\binom{l}{i}a(m+l-i)v(n+i)v_{2}-\sum_{i=0}^{\infty}(-1)^{l+i}\binom{l}{i}v(n+l-i)a(m+i)v_{2},\nonumber
\end{align}
for all $m,n,l\in\Z$, $a\in V$, $v\in M^1,$ and $v_{2}\in M^{2}$. Note that \eqref{n2.17} is the component form of the Jacobi identity for the intertwining operator $I_S$ (see (1.2.9) in \cite{Z}).

Indeed,	by \eqref{n2.14} and the generating property \eqref{n2.12} of $S$, we have:

\begin{align*}
& \langle v_{3}', \sum_{i=0}^{\infty}(-1)^{i}\binom{l}{i}a(m+l-i)v(n+i)v_{2}\rangle\\
&=\sum_{i=0}^{\infty}(-1)^{i}\binom{l}{i} \int_{C_{1}'}S( a'(m+l-i)v_{3}', (v,w)v_{2})w^{n+i}dw \numberthis \label{n2.18}\\
&=\sum_{i=0}^{\infty}(-1)^{i}\binom{l}{i} \int_{C_{1}'}\int_{C_{2}'}S( v_{3}', (a,z)(v,w)v_{2})z^{m+l-i}w^{n+i}dw,
\end{align*}
where $C_{1}'$ is a contour of $w$, and $C_{2}'$ is a contour of $z$ which contains $C_{1}'$.
On the other hand, by \eqref{n2.14} and the generating property \eqref{n2.9} of $S$, we have:
\begin{align*}
&\langle v_{3}', \sum_{i=0}^{\infty}(-1)^{l+i}\binom{l}{i}v(n+l-i)a(m+i)v_{2} \rangle\\
&=\sum_{i=0}^{\infty}(-1)^{l+i}\binom{l}{i} \int_{C_{1}}S(v_{3}',(v,w)a(m+i)v_{2})w^{n+l-i}dw \numberthis \label{n2.19}\\
&= \sum_{i=0}^{\infty}(-1)^{l+i}\binom{l}{i} \int_{C_{1}}\int_{C_{2}}S(v_{3}',(v,w)(a,z)v_{2})z^{m+i}w^{n+l-i}dzdw,
\end{align*}
where $C_{1}$ and $C_{2}$ are contours in $w$ and $z$, respectively, and $C_2$ is contained in $C_{1}$. \par 
We adopt the notations in Proposition A.2.8 in \cite{FLM}. Choose the contours $C_{1},C_{2},C_{1}'$, and $C_{2}'$ in the following way: 
Let $C_{\alpha}^{z}$ be a circle in the variable $z$ centered at $0$, with radius $\alpha$, and $C_{\epsilon}^{1}(w_{2})$ be the circle of $w_{1}$ centered at $w_{2}$ with radius $\epsilon$. We may choose $\epsilon$ small enough so that $|w_{1}-w_{2}|<|w_{2}|$ for any $w_{1}$ lying on $C_{\epsilon}^{1}(w_{2})$. Choose $R,r,\rho>0$ so that $1>R>\rho>r$. Let $C_{1}'=C_{\rho}^{w}$, $C_{2}'=C_{R}^{z}$, $C_{1}=C_{\rho}^{w}$, and $C_{2}=C_{r}^{z}$. Then by \eqref{n2.14}, \eqref{n2.18}, and \eqref{n2.19}, together with (2) and (5) in Definition \ref{df2.1}, we have:

\begin{align*}
&\langle v_{3}', \sum_{i=0}^{\infty}(-1)^{i}\binom{l}{i}a(m+l-i)v(n+i)v_{2}-\sum_{i=0}^{\infty}(-1)^{l+i}\binom{l}{i}v(n+l-i)a(m+i)v_{2} \rangle\\
&=\sum_{i=0}^{\infty}(-1)^{i}\binom{l}{i} \int_{C_{\rho}^{w}}\int_{C_{R}^{z}}S( v_{3}', (a,z)(v,w)v_{2})z^{m+l-i}w^{n+i}dwdz\\
&\ \ \ - \sum_{i=0}^{\infty}(-1)^{l+i}\binom{l}{i} \int_{C_{\rho}^{w}}\int_{C_{r}^{z}}S(v_{3}',(v,w)(a,z)v_{2})z^{m+i}w^{n+l-i}dzdw\\
&=\int_{C_{\rho}^{w}}\int_{C_{R}^{z}}S( v_{3}', (a,z)(v,w)v_{2})\iota_{z,w}(z-w)^{l}z^{m}w^{n}dwdz\\
&\ \ \ - \int_{C_{\rho}^{w}}\int_{C_{r}^{z}}S(v_{3}',(v,w)(a,z)v_{2})\iota_{w,z}(z-w)^{l}z^{m}w^{n}dzdw \numberthis\label{n2.20}\\
&=\int_{C_{\rho}^{w}}\int_{C_{\epsilon}^{z}(w)}S(v_{3}',(a,z)(v,w)v_{2})(z-w)^{l}z^{m}w^{n}dzdw\\
&=\int_{C_{\rho}^{w}}\int_{C_{\epsilon}^{z}(w)}S(v_{3}',(a,z)(v,w)v_{2})(z-w)^{l}\iota_{w,z-w}(w+(z-w))^{m}w^{n}dzdw\\
&=\sum_{i\geq 0} \binom{m}{i}\int_{C_{\rho}^{w}}\int_{C_{\epsilon}^{z}(w)}S(v_{3}',(a,z)(v,w)v_{2})(z-w)^{l+i}w^{n+m-i}dzdw\\
&=\sum_{i\geq 0} \binom{m}{i}\int_{C_{\rho}^{w}}S(v_{3}',(a(l+i)v,w)v_{2})w^{m+n-i}\\
&=\sum_{i\geq 0} \binom{m}{i}\langle v_{3}',(a(l+i)v)(m+n-i)v_{2}\rangle .
\end{align*}

The graph of the contours appear in \eqref{n2.20} can be sketched as follows: 
\begin{center}
	\begin{tikzpicture} [scale=0.7]
	\draw  [thick] circle [radius=3] node[xshift=2.4cm]{$C_{\rho}^{w}$}
	(0,0)--node[right]{$\rho$}(-30:3cm);
	\draw[thick] circle [radius=4.0]   node[xshift=2.3cm, yshift=2.0cm]{$C_{R}^{z}$}
	--node[right]{$R$}(-120:4cm);
	\draw (0,3)[thick] circle [radius=0.75] node[below ,yshift=-0.4cm]{$C_{\epsilon}^{z}(w)$};
	\draw[thick] circle [radius=1] node[xshift=0.95cm]{$C_{r}^{z}$}--node[left]{$r$}(125:1cm);
	\end{tikzpicture}
\end{center}

Since $v_{3}'$ in \eqref{n2.20} can be choosen arbitraily, the Jacobi identity \eqref{n2.17} follows, and so $I_S$ given by \eqref{n2.15} is an intertwining operator of type $\fusion{M^{1}}{M^{2}}{M^3}$.
\end{proof}
\begin{coro}\label{coro2.6}
The vector space of intertwining operators $I\fusion{M^{1}}{M^{2}}{M^3}$ is isomorphic to the vector space $\Cor\fusion{M^{1}}{M^{2}}{M^3}$ in Definition \ref{df2.2}. 
\end{coro}
\begin{proof}
Theorem \ref{thm2.5} indicates that there exists a well-defined linear map:
\begin{equation}\label{n2.21}
\b: \Cor\fusion{M^{1}}{M^{2}}{M^3}\ra I \fusion{M^{1}}{M^{2}}{M^3}, \quad S\mapsto I_{S}.
\end{equation}
By \eqref{n2.3} and \eqref{n2.16}, it is clear that $\b$ is an inverse of the linear map $\al$ in \eqref{n2.13}. Hence $I\fusion{M^{1}}{M^{2}}{M^3}\cong \Cor\fusion{M^{1}}{M^{2}}{M^3}$ as vector spaces.
\end{proof}

\section{Extension of correlation functions from the bottom levels}
Let $M^{2}$ and $M^{3}$ be any $V$-modules with bottom levels $M^{2}(0)$ and $M^3(0)$, respectively.

Recall the bottom level $M(0)$ of any $\N$-gradable $V$-module $M=\bigoplus_{n=0}^{\infty}M(n)$ is a module over the Zhu's algebra $A(V)$ defined in \cite{Z} or the generalized Zhu's algebra $A_n(V)$ defined in \cite{DLM2} under the module action: $$[a].v=o(a)v=a(\wt a-1)v,$$ for all $[a]\in A(V)$ or $A_{n}(V)$, and $v\in M(0)$ (see Theorem 2.1.2 in \cite{Z}). 

In this section, we assume that the $A(V)$-modules $M^{2}(0)$ and $M^{3}(0)$ are irreducible.

\subsection{The space of correlation functions associated with $M^{1}$, $M^{2}(0)$, and $M^{3}(0)$}

Let $S\in \Cor \fusion{M^{1}}{M^{2}}{M^3}$, and let $I\in I\fusion{M^{1}}{M^{2}}{M^3}$ be its corresponding intertwining operator under the isomorphism $\b$ in \eqref{n2.21}. For each $n\in \N$, consider the restriction of $S$ onto the bottom levels $M^2(0)$ and $M^3(0)^\ast$:
\begin{equation}\label{n3.1}
S|_{M^{3}(0)^{\ast}\times...M^{1}...\times M^{2}(0)}: M^{3}(0)^{\ast}\times V\times \ds\times M^{1}\ds\times V\times M^{2}(0)\rightarrow \mathcal{F}(z_{1},\ds,z_{n},w).
\end{equation}
To simplify our notation, we use the same symbol $S$ to denote the restricted function \eqref{n3.1}. Clearly, $S$ in \eqref{n3.1} satisfies properties (1)-(6) in Definition \ref{df2.1}, with the elements $v_{3}'$ and $v_{2}$ in these properties belong to $M^{3}(0)^\ast$ and $M^{2}(0)$, respectively. Moreover, since $(v_{3}',I(v,w)v_{2})=S(v_{3}',(v,w)v_{2})$ by \eqref{n2.16}, and $v(n)M^{2}(m)\subseteq M^{3}(m+\deg v-n-1)$ for all $v\in M^1$ homogeneous, $n\in \Z$, and $m\in \N$ (see (1.5.4) in \cite{FZ}), then we have: 
\begin{equation}\label{n3.2}
S(v_{3}',(v,w)v_{2})=\<v_{3}', v(\deg v-1)(v_{2})\> w^{-\deg v}.
\end{equation}
We introduce the following intermediate notion based on the properties satisfied by the system of restricted correlation functions \eqref{n3.1}. 

\begin{df}\label{df3.1}
Let $M^{2}(0)$ and $M^{3}(0)$ be irreducible $A(V)$-modules. A system of linear maps $S=\{S^n_{V...M^1...V}\}_{n=0}^\infty$, 
\begin{align*}
S_{V...M^1...V}^n: M^{3}(0)^\ast\times V\times\ldots\times M^{1}&\times\ds V\times M^{2}(0)\rightarrow \mathcal{F}(z_{1},
\ds,z_{n},w),\\
(v_{3}',a_{1},\ds,v,\ds,a_{n},v_{2})&\mapsto S(v_{3}',(a_{1},z_{1})\ds(v,w)\ds(a_{n},z_{n})v_{2}),
\end{align*}
is said to satisfy the genus-zero property associated with $M^1$, $M^2(0)$, and $M^3(0)$ if it satisfies the following: 
\begin{enumerate}
	\item Properties $(2)-(6)$ in Definition \ref{df2.1}, with the elements $v_{3}'$ and $v_{2}$ in these properties belong to $M^{3}(0)^\ast$ and $M^{2}(0)$, respectively.
	\item There exists a linear functional $f:M^1\ra \Hom_\C(M^2(0),M^3(0)),v\mapsto f_v$, such that
	\begin{equation}\label{n3.3}
	S(v_{3}',(v,w)v_{2})=\<v_{3}',f_{v}(v_{2})\>w^{-\deg v},
	\end{equation}
	for all $v_{2}\in M^2(0)$ and $v_{3}'\in M^3(0)^\ast$.
	\item (The recursive formula for $M^3(0)^\ast$) For any $v_{3}'\in M^{3}(0)^{\ast}$, $v\in M^{1}$, $v_{2}\in M^{2}(0)$, and $a_{1},\ds ,a_{n}\in V$, 
	\begin{align*}
	&S(v_{3}',(a,z)(a_{1},z_{1})\ds (a_{n},z_{n})(v,w)v_{2})=S(v_{3}'o(a),(a_{1},z_{1})\ds(a_{n},z_{n})(v,w)v_{2})z^{-\wt a}\\
	&\quad +\sum_{k=1}^{n}\sum_{i\geq 0} F_{\wt a,i}(z,z_{k}) S(v_{3}',(a_{1},z_{1})\ds(a(i)a_{k},z_{k})\ds(a_{n},z_{n})(v,w)v_{2})\numberthis \label{n3.4}\\
	&\quad +\sum_{i\geq 0} F_{\wt a,i}(z,w) S(v_{3}', (a_{1},z_{1})\ds(a_{n},z_{n})(a(i)v,w)v_{2}),
	\end{align*}
	where $F_{\wt a,i}(z,w)$ is a rational function in $z,w$ given by:
	\begin{equation}\label{n3.5}
	\begin{aligned}
	\iota_{z,w}(F_{\wt a,i}(z,w))&=\sum_{j\geq 0} \binom{\wt a+j}{i}z^{-\wt a-j-1}w^{\wt a+j-i},\\
	F_{m,i}(z,w)&=\frac{z^{-m}}{i!}\bigg(\frac{d}{dw}\bigg)^{i}\frac{w^{m}}{z-w},
	\end{aligned}
	\end{equation} 
	for any $m\in \N$, and $v_{3}'o(a)$ is given by the natural right module action on $M^3(0)^\ast$.
	\item (The recursive formula for $M^2(0)$) For any $v_{3}'\in M^{3}(0)^{\ast}$, $v\in M^{1}$, $v_{2}\in M^{2}(0)$, and $a_{1},\ds,a_{n}\in V$, we have:
	\begin{align*}
	&S(v_{3}', (a_{1},z_{1})\ds(a_{n},z_{n})(v,w)(a,z)v_{2})=S(v_{3}', (a_{1},z_{1})\ds(a_{n},z_{n})(v,w)o(a)v_{2})z^{-\wt a}\\
	&\quad +\sum_{k=1}^{n}\sum_{i\geq 0} G_{\wt a,i}(z,w)S(v_{3}', (a_{1},z_{1})\ds(a(i)a_{k},z_{k})...(a_{n},z_{n})(v,w)v_{2})\numberthis\label{n3.6}\\
	&\quad +\sum_{i\geq 0} G_{\wt a,i}(z,w)S(v_{3}', (a_{1},z_{1})\ds(a_{n},z_{n})(a(i)v,w)(a,z)v_{2}),
	\end{align*}
	where $G_{\wt a,i}(z,w)$ is a rational function defined by
	\begin{equation}\label{n3.7}
	\begin{aligned}
	\iota_{w,z}(G_{\wt a,i}(z,w))&=-\sum_{j\geq 0}\binom{\wt a-2-j}{i}w^{\wt a-j-2-i} z^{-\wt a+1+j},\\
	G_{m,i}(z,w)&=\frac{z^{-m+1}}{i!}\bigg(\frac{d}{dw}\bigg)^{i}\bigg(\frac{w^{m-1}}{z-w}\bigg),
	\end{aligned}
	\end{equation}
	for any $m\in \N$. 
\end{enumerate}
The vector space of the system of functions satisfying the genus-zero property associated with $M^1$, $M^2(0)$, and $M^3(0)$ is denoted by $\Cor\fusion{M^1}{M^2(0)}{M^3(0)}$.
\end{df}

We observe that the rational functions $F$ and $G$ given by \eqref{n3.5} and \eqref{n3.7} satisfy the following relation:
$$
F_{m,i}(z,w)-G_{m,i}(z,w)
=\frac{z^{-m}}{i!}\left(\frac{d}{dw}\right)^{i} \bigg(\frac{w^{m}}{z-w}-\frac{zw^{m-1}}{z-w}\bigg)\\
=-\binom{m-1}{i}z^{-m}w^{m-1-i},
$$
for all $m\in \N$. In particular, we have 
\begin{equation}\label{n3.8}
F_{\wt a,i}(z,w)-G_{\wt a,i}(z,w)=-\binom{\wt a-1}{i}z^{-\wt a}w^{\wt a-1-i}.
\end{equation}
The equation \eqref{n3.8} will be used multiple times in Section 4 when we build a system of correlation functions $S$ from a linear map on a tensor product of $A(V)$-modules. 

\begin{prop}\label{prop3.2}
Let $S\in \Cor\fusion{M^1}{M^2}{M^3}$. Then the system of restricted functions $S$ in \eqref{n3.1} satisfies the genus-zero property associated with $M^1$, $M^2(0)$, and $M^3(0)$.
\end{prop}
\begin{proof}
By our discussion in the begining of this subsection, $S$ in \eqref{n3.1} satisfies (1) and (2) in Definition \ref{df3.1}, where the $f_v$ in \eqref{n3.3} is given by $f_v=v(\deg v-1)$, for all $v\in M^1$. The proof of \eqref{n3.4} is similar to the proof of Lemma 2.2.1 in \cite{Z}. We omit the details. To prove \eqref{n3.6}, we only consider the case when $n=0$ (the general case follows from a similar argument.) Note that $a(n)v_{2}=0$ if $\wt a-n-1<0$, it follows that
$\langle v_{3}', I(v,w)Y(a,z)v_2 \rangle=\langle v_{3}', I(v,w)o(a)v_{2}\rangle z^{-\wt a}+ \sum_{\wt a-n-1>0}\langle v_{3}',I(v,w)a(n)v_{2}\rangle z^{-n-1}.  $
By the definition of contragredient modules, we have $\langle v_{3}', a(n)u\rangle= \sum_{i\geq 0}\frac{1}{i!}(-1)^{i}\langle (L(i)a)(2\wt a-n-i-2)v_{3}',u\rangle,$ for any $n\in \Z$.
But $(L(i)a)(2\wt a-n-i-2)v_{3}'\in M^{3'}(-\wt a+n+1)=0$ when $\wt a-n-1>0$. Thus 
\begin{align*}
&\sum_{wta-n-1>0}\langle v_{3}',I(v,w)a(n)v_{2}\rangle z^{-n-1}=-\sum_{\wt a-n-1>0}\langle v_{3}', [a(n),I(v,w)]v_{2}\rangle z^{-n-1}\\
&=-\sum_{\wt a-n-1>0}\sum_{i\geq 0}\binom{n}{i}\langle v_{3}', I(a(i)v,w)v_{2}\rangle z^{-n-1}w^{n-i}\\
&=-\sum_{j\geq 0}\sum_{i\geq 0}\binom{\wt a-j-2}{i}z^{-\wt a+j+2-1}w^{\wt a-j-2-i}\langle v_{3}', I(a(i)v,w)v_{2}\rangle\\
&=\sum_{i\geq 0} \iota_{w,z}(G_{\wt a,i}(z,w))\langle v_{3}', I(a(i)v,w)v_{2}\rangle,
\end{align*}
where the last equality follows from \eqref{n3.7}. Hence we have:
$$\langle v_{3}', I(v,w)Y(a,z) \rangle=\langle v_{3}', I(v,w)o(a)v_{2}\rangle z^{-\wt a}+\sum_{i\geq 0} \iota_{w,z}(G_{\wt a,i}(z,w))\langle v_{3}', I(a(i)v,w)v_{2}\rangle$$
as power series. By taking the limit of this series, we obtain \eqref{n3.6} for $n=0$. 
\end{proof}
As a consequence of Proposition \ref{prop3.2}, we have a well-defined restriction map:
\begin{equation}\label{n3.9}
\varphi: \Cor\fusion{M^{1}}{M^{2}}{M^{3}}\ra \Cor\fusion{M^{1}}{M^{2}(0)}{M^{3}(0)}, \quad S\mapsto S|_{M^{3}(0)^{\ast}\times...M^{1}...\times M^{2}(0)},
\end{equation}
where $M^{2}$ and $M^{3}$ are any $V$-modules with bottom levels $M^{2}(0)$ and $M^{3}(0)$,

The following Lemma will be used in the next subsection:
\begin{lm}\label{lm3.3}
Let $S\in \Cor\fusion{M^{1}}{M^{2}(0)}{M^{3}(0)}$, and let $f:M^{1}\ra \Hom_{\C}(M^{2}(0),M^{3}(0)), v\mapsto f_{v}$ be the linear functional in Definition \ref{df3.1}. Suppose that $f_{v}=0$ for all $v\in M^{1}$. Then $S=0$. 
\end{lm}
\begin{proof}
We use induction on $n$ to show that $S(v_{3}', (a_{1},z_{1})\ds(a_{n},z_{n})(v,w)v_{2})=0$ for all $v_{3}'\in M^{3}(0)^{\ast}, v\in M^{1}$, $v_{2}\in M^{2}(0)$, and $a_{1},\ds,a_{n}\in V$. When $n=0$, by the assumption and \eqref{n3.3}, we have: $S(v_{3}',(v,w)v_{2})=\<v_{3}',f_{v}(v_{2})\>w^{-\deg v}=\<v_{3}',0\>w^{-\deg w}=0,$
for all $v_{3}'\in M^{3}(0)^{\ast}$, $v\in M^{1}$, and $v_{2}\in M^{2}(0)$. For $n>0$, by the recursive formula \eqref{n3.4}, we have 
\begin{align*}
&S(v_{3}', (a_{1},z_{1})\ds(a_{n},z_{n})(v,w)v_{2})=S(v_{3}'o(a_{1}), (a_{2},z_{2})\ds(a_{n},z_{n})(v,w))z^{-\wt a_{1}}\\
&\quad +\sum_{k=2}^{n}\sum_{i\geq 0} F_{\wt a_{1},i}(z_{1},z_{k}) S(v_{3}', (a_{2},z_{2})\ds(a_{1}(i)a_{k},z_{k})\ds(a_{n},z_{n})(v,w)v_{2})\\
&\quad +\sum_{i\geq 0}F_{\wt a_{1}, i}(z_{1},w) S(v_{3}', (a_{2},z_{2})\ds(a_{n},z_{n})(a_{1}(i)v,w)v_{2}).
\end{align*}
Since each term on the right-hand side has a smaller length, the right-hand side is equal to $0$ by the induction hypothesis, so we have $S(v_{3}', (a_{1},z_{1})\ds(a_{n},z_{n})(v,w)v_{2})=0$. 
\end{proof}

\subsection{Extension from the bottom levels} In this subsection, we will show that the restriction map $\varphi$ in \eqref{n3.9} has an inverse for certain $V$-modules $M^{2}$ and $M^{3}$, with (irreducible) bottom levels $M^2(0)$ and $M^3(0)$, respectively.

Recall that for any irreducible $A(V)$-module $U$, Dong, Li, and Mason constructed a generalized Verma module $\bar{M}(U)$ in \cite{DLM1}. By construction, $\bar{M}(U)=(U(\mathcal{L}(V))\otimes_{U(\mathcal{L}(V)_{\geq 0})}U)/U(\mathcal{L}(V))W,$ where \begin{equation}\label{n3.10}
\mathcal{L}(V)=V\otimes \C[t,t^{-1}]/(L(-1)\otimes 1+1\otimes \frac{d}{dt})(V\otimes \C[t,t^{-1}])
\end{equation} is the Lie algebra associated with the VOA $V$ (cf.\cite{B,DLM1}), and $W$ is the subspace of $U(\mathcal{L}(V))\otimes_{U(\mathcal{L}(V)_{\geq 0})}U$ spanned by the coefficients of the weak associativity equality, see Section 5 in \cite{DLM1} for more details.

$\bar{M}(U)$ is $\N$-gradable: $\bar{M}(U)=\bigoplus_{n=0}^{\infty} \bar{M}(n)$, with the bottom level $\bar{M}(U)(0)=U$. It satisfies a universal property in the sense that any $\N$-gradable $V$-module with bottom level $U$ is a quotient module of $\bar{M}(U)$ (Theorem 6.2 in \cite{DLM1}). Moreover, $\bar{M}(U)$ admits a unique maximal graded $\mathcal{L}(V)$-submodule $J$ subject to $J\cap U=0$, and $L(U)=\bar{M}(U)/J$ is an irreducible $V$-module (Theorem 6.3 in \cite{DLM2}).

In Section 2 of \cite{L}, Li gave an alternative definition of the generalized Verma module $\bar{F}(U)$ associated with $U$, namely, $\bar{F}(U)=(U(\mathcal{L}(V))\otimes_{U(\mathcal{L}(V)_{\geq 0})}U)/J(U),$
where $J(U)$ is the intersection of $\ker \al$, where $\al$ runs over all $\mathcal{L}(V)$-homomorphisms from $\bar{F}(U)$ to weak $V$-modules. Clearly, $\bar{M}(U)=\bar{F}(U)$ since they satisfy the same universal property.

Choose an element 
\begin{equation}\label{n3.11}
S:M^{3}(0)^{*}\times V\times \ds\times M^{1}\times \ds\times V\times M^{2}(0)\ra \mathcal{F}(z_{1},...,z_{n},w)
\end{equation}
in $\Cor\fusion{M^{1}}{M^{2}(0)}{M^{3}(0)}$. We will extend the first and the last input vector spaces from $M^{3}(0)^{\ast}$ and $M^{2}(0)$ to some $V$-modules $\tilde{M}/\Rad \tilde{M}$ and $\bar{M}/\Rad \bar{M}$, which are certain quotient modules of the generalized Verma modules $\bar{M}(M^{3}(0)^\ast)$ and $\bar{M}(M^2(0))$, respectively.

We first extend $M^{2}(0)$, and we will proceed like the proof of Theorem 2.2.1 in \cite{Z}. In our case, however, the extended $V$-module is {\em not} necessarily irreducible like the extended module in Theorem 2.2.1 \cite{Z} . 

Let $\bar{M}:=T(\mathcal{L}(V))\otimes_{\C} M^{2}(0)$, where $T(\mathcal{L}(V))$ is the tensor algebra of $\mathcal{L}(V)$. To simplify our notation, we omit the tensor symbol in an element of $\bar{M}$ and denote an element $\overline{b\otimes t^{n}}$ in $\mathcal{L}(V)$ by $(b,n)$, then an element in $\bar{M}$ can be written as:
\begin{equation}\label{n3.12}
(b_{1},i_{1})(b_{1},i_{2})\ds(b_{m},i_{m})v_{2}
\end{equation}
where $b_{i}\in V, \ i_{k}\in\mathbb{Z},\ v_{2}\in M^{2}(0)$, and $(b,i)$ linear in $b$. Denote the vector in \eqref{n3.12} by $x$. Extend $M^2(0)$ to $\bar{M}$ by repeatedly using the generating formula \eqref{n2.9}. i.e., we let:

\begin{align*}
&S:M^{3}(0)^{*}\times V\times \ds\times M^{1}\times \ds\times V\times \bar{M}\ra\mathcal{F}(z_{1},...,z_{n},w),\\
&S(v_{3}',(a_{1},z_{1})\ds(a_{n},v_{n})(v,w)x)\numberthis \label{n3.13}\\
&:=\int_{C_{1}}\ds\int_{C_{m}}S(v_{3}',(a_{1},z_{1})\ds(a_{n},z_{n})(v,w)(b_{1},w_{1})\ds(b_{m},w_{m})v_{2})w_{1}^{i_{1}}\ds w_{m}^{i_{m}}dw_{1}\ds dw_{m},
\end{align*}
where $C_{k}$ is a contour of $w_{k}$, $C_{k}$ contains $C_{k+1}$ for each $k$, $C_{m}$ contains $0$, and $z_{1},\ds,z_{n}, w$ are lying outside of $C_{1}$. For the well-definedness of $S$ in \eqref{n3.13}, by \eqref{n3.10}, we just need to show that $S$ in \eqref{n3.13} agrees on the elements:
$$(b_{1},i_{1})\ds (L(-1)b_{k},i_{k})\ds (b_{m},i_{m})v_{2},\quad \mathrm{and}\quad -i_{k}(b_{1},i_{1})\ds(b_{k},i_{k}-1)\ds(b_{m},i_{m})v_{2}.$$
Indeed, by the Definition \ref{df3.1}, $S$ in \eqref{n3.11} satisfies \eqref{n2.6}. Thus, 
\begin{align*}
&S(v_{3}',(a_{1},z_{1})...(a_{n},v_{n})(v,w)(b_{1},i_{1})...(L(-1)b_{k},i_{k})...(b_{m},i_{m})v_{2})\\
&=\int_{C_{1}}...\int_{C_{m}}\frac{d}{dw_{k}}S(v_{3}',(a_{1},z_{1})...(a_{n},z_{n})(v,w)...(b_{k},w_{k})...v_{2})w_{1}^{i_{1}}...w_{k}^{i_{k}}...w_{m}^{i_{m}}dw_{1}...dw_{m}\\
&=-\int_{C_{1}}...\int_{C_{m}}S(v_{3}',(a_{1},z_{1})...(a_{n},z_{n})(v,w)...(b_{k},w_{k})...v_{2})w_{1}^{i_{1}}...(i_{k})w_{k}^{i_{k}-1}...w_{m}^{i_{m}}dw_{1}...dw_{m}\\
&=S(v_{3}',(a_{1},z_{1})...(a_{n},v_{n})(v,w)(-i_{k})(b_{1},i_{1})...(b_{k},i_{k}-1)...(b_{m},i_{m})v_{2}).
\end{align*}

Introduce a gradation on $\bar{M}$ by letting
\begin{equation}\label{n3.14}
\deg ((b_{1},i_{1})(b_{1},i_{2})\ds(b_{m},i_{m})v_{2}):=\sum_{k=1}^{m} (\wt b_{k}-i_{k}-1),
\end{equation}
and denote the degree $n$ subspace by $\bar{M}(n)$. Then $\bar{M}=\bigoplus_{n\in \mathbb{Z}}\bar{M}(n)$, with $M^{2}(0)\subseteq \bar {M}(0)$.

Similar to (2.2.30) in \cite{Z}, we define the radical of $S$ on $\bar{M}$ by
\begin{equation}\label{n3.15}
\begin{aligned}
\Rad(S):=\{x\in \bar{M}|&S(v_{3}',(a_{1},z_{1})\dots(a_{n},z_{n})(v,w)x)=0,\\
& \forall  n\geq 0,\ a_{1},\dots a_n\in V,\ v\in M^{1},\ v_{3}\in M^{3}(0)^{*} \},
\end{aligned}
\end{equation}
then let $\Rad(\bar{M}):=\bigcap_{S}\Rad(S),$
where the intersection is taken over all $S\in \Cor\fusion{M^{1}}{M^{2}(0)}{M^{3}(0)}$. In fact, we can take the intersection over all nonzero $S$, since $\Rad(S)=\bar{M}$ if $S=0$. 

It is clear that the extended $S$ in \eqref{n3.13} factors through $\bar{M}/\Rad(\bar{M})$.
Next, we show that $\bar{M}/\Rad(\bar{M})$ carries a structure of $\N$-gradable $V$-module whose bottom level is $M^{2}(0)$.

\begin{lm}\label{lm3.4}
Let $W$ be the subspace of $\bar{M}$ spanned by the following elements:
\begin{align}\label{n3.16}
&\sum_{i=0}^{\infty}\binom{m}{i}(a(l+i)b,m+n-i)x\\
-\bigg(&\sum_{i=0}^{\infty}(-1)^{i}\binom{l}{i}(a,m+l-i)(b,n+i)x-\sum_{i=0}^{\infty}(-1)^{l+i}\binom{l}{i}(b,n+l-i)(a,m+i)x \bigg),\nonumber
\end{align}
where $a,b\in V,\ m,n,l\in \mathbb{Z},$ and $ x\in \bar{M}$. Then we have $W\subset \Rad(\bar{M})$.
\end{lm}
\begin{proof}
By the formula \eqref{n3.13}, it is easy to see that for the following element in $\bar{M}$:
$$x'=(b_{1},i_{1})\ds(b_{m},i_{m})x,$$
where $x= (c_{1},j_{1})\ds(c_{n},j_{n})v_{2}$ for some $b_{i},c_{j}\in V$ and $i_{k},j_{l}\in \Z$, we have:
\begin{equation}\label{n3.17}
\begin{aligned}
&S(v_{3}',(a_{1},z_{1})...(a_{n},v_{n})(v,w)x')\\
&=\int_{C_{1}}...\int_{C_{m}}S(v_{3}',(a_{1},z_{1})...(a_{n},z_{n})(v,w)(b_{1},w_{1})...(b_{m},w_{m})x)w_{1}^{i_{1}}...w_{m}^{i_{m}}dw_{1}...dw_{m},
\end{aligned}
\end{equation}
where $C_{k}$ is a contour of $w_{k}$, $C_{k+1}$ is inside of $C_{k}$ for each $k$, $C_{m}$ contains $0$, and $z_{1},\ds,z_{n},w$ are lying outside of $C_{1}$. Now we fix a nonzero element $S\in \Cor\fusion{M^{1}}{M^{2}(0)}{M^{3}(0)}$. 

Denote the element \eqref{n3.16} by $y$. We adopt the notations in Proposition A.2.8 in \cite{FLM} again. Let $C_{R}^{i}$ be the circle of $w_{i}, \ i=1,2$, centered at $0$ with radius $R$, and let $C_{\epsilon}^{1}(w_{2})$ be the circle of $w_{1}$ centered at $w_{2}$ with radius $\epsilon$. We may choose $\epsilon$ small enough so that $|w_{1}-w_{2}|<|w_{2}|$ for any $w_{1}$ lying on $C_{\epsilon}^{1}(w_{2})$. Choose $R,r,\rho>0$ so that $R>\rho>r$.
By \eqref{n3.17} and the locality (2) in Definition \ref{df2.1} of $S$, we have:
\begin{align*}
&S(v_{3}',(a_{1},z_{1})...(a_{n},z_{n})(v,w)y)\\
&=\int_{C_{\rho}^{2}}\sum_{i=0}^{\infty}\binom{m}{i}S(v_{3}',(a_{1},z_{1})...(a_{n},z_{n})(v,w)(a(l+i)b,w_{2})x)w_{2}^{m+n-i}dw_{2}\\
&\ \ -\int_{C_{R}^{1}}\int _{C_{\rho}^{2}}\sum_{i=0}^{\infty}(-1)^{i}\binom{l}{i}S(v_{3}',(a_{1},z_{1})...(a_{n},z_{n})(v,w)(a,w_{1})(b,w_{2})x)w_{1}^{m+l-i}w_{2}^{n+i}dw_{1}dw_{2}\\
&  \ \ +\int_{C_{\rho}^{2}}\int_{C_{r}^{1}}\sum_{i=0}^{\infty}(-1)^{l+i}\binom{l}{i}S(v_{3}',(a_{1},z_{1})...(a_{n},z_{n})(v,w)(b,w_{2})(a,w_{1})x)w_{1}^{m+i}w_{2}^{n+l-i}dw_{1}dw_{2}\\
&=\int_{C_{\rho}^{2}}\sum_{i=0}^{\infty}\binom{m}{i}S(v_{3}',(a_{1},z_{1})...(a_{n},z_{n})(v,w)(a(l+i)b,w_{2})x)w_{2}^{m+n-i}dw_{2}\\
& \ \ -\int_{C_{R}^{1}}\int _{C_{\rho}^{2}}S(v_{3}',(a_{1},z_{1})...(a_{n},z_{n})(v,w)(a,w_{1})(b,w_{2})x)\cdot \iota_{w_{1},w_{2}}((w_{1}-w_{2})^{l})w_{1}^{m}w_{2}^{n}dw_{1}dw_{2}\\
& \ \ +\int_{C_{\rho}^{2}}\int_{C_{r}^{1}}S(v_{3}',(a_{1},z_{1})...(a_{n},z_{n})(v,w)(b,w_{2})(a,w_{1})x)\cdot \iota_{w_{2},w_{1}}((-w_{2}+w_{1})^{l})w_{1}^{m}w_{2}^{n}dw_{1}dw_{2}\\
&=\int_{C_{\rho}^{2}}\sum_{i=0}^{\infty}\binom{m}{i}S(v_{3}',(a_{1},z_{1})...(a_{n},z_{n})(v,w)(a(l+i)b,w_{2})x)w_{2}^{m+n-i}dw_{2}\\
&\ \ -\int_{C_{\rho}^{2}}\int_{C_{\epsilon}^{1}(w_{2})}S(v_{3}',(a_{1},z_{1})...(a_{n},z_{n})(v,w)(a,w_{1})(b,w_{2})v_{2})(w_{1}-w_{2})^{l}w_{1}^{m}w_{2}^{n}dw_{1}dw_{2}.\\
&=\int_{C_{\rho}^{2}}\sum_{i=0}^{\infty}\binom{m}{i}S(v_{3}',(a_{1},z_{1})...(a_{n},z_{n})(v,w)(a(l+i)b,w_{2})x)w_{2}^{m+n-i}dw_{2}\\
&\ \ -\int_{C_{\rho}^{2}}\int_{C_{\epsilon}^{1}(w_{2})}\sum_{i=1}^{\infty}\binom{m}{i}S(v_{3}',(a_{1},z_{1})...(v,w)(a,w_{1})(b,w_{2})v_{2})(w_{1}-w_{2})^{l+i}w_{2}^{m+n-i}dw_{1}dw_{2}\\
&=0,
\end{align*}
for all $v_{3}'\in M^3(0)^\ast$, $a_{1},\dots a_n \in V$, and $v\in M^1$, where the last equality follows from the associativity (5) in Definition \ref{df2.1}. This shows $y\in \Rad (S)$. But $S$ is chosen arbitrarily. Hence we have $y\in \Rad(\bar{M})$. 
\end{proof}
The following facts are satisfied by $\Rad(\bar{M})$:
\begin{lm}\label{lm3.5}
(a) If $x\in \Rad(\bar{M})$, then $(b,i)x\in \Rad(\bar{M})$, for any $b\in V$ and $ i\in \mathbb{Z}$.\\
(b) $M^{2}(0)\cap \Rad(\bar{M})=0$. \\
(c)  $\bar{M}(n)\subset \Rad(\bar{M})$ for all $n<0$. 
\end{lm}
\begin{proof}
Since $\Rad(\bar{M})=\bigcap_{S} \Rad(S)$, we just need to show that (a), (b), and (c) hold for $\Rad(S)$, where $S\in \Cor\fusion {M^{1}}{M^{2}(0)}{M^{3}(0)}$ is nonzero.

$(a)$ Let $x\in \Rad(S)$, by \eqref{n3.13} and the definition \eqref{n3.15} of $\Rad(S)$, we have
$$S(v_{3}', (a_{1},z_{1})...(v,w)(b,i)x)=\int_{C}S(v_{3}', (a_{1},z_{1})...(v,w)(b,w_{1})x)w_{1}^{i}dw_{1}=\int_{C}0\cdot w_{1}^{i}dw_{1}=0,$$
where $C$ is a contour of $w_{1}$, with $z_{1},\ds,z_{n},w$ lying outside. Thus $(b,i)x\in \Rad(S)$.

$(b)$ Suppose there exists some $v_{2}\neq 0$ in $M^{2}(0)\cap \Rad(S)$, then by \eqref{n3.3} and the recursive formula \eqref{n3.6}, we have
\begin{align*}
0&=\iota_{w,z}(S(v_{3}',(a,z)(v,w)v_{2}))\\
&=S(v_{3}',(v,w)o(a)v_{2})z^{-\wt a}+\sum_{i\geq 0} \iota_{w,z}(G_{\wt a,i}(z,w))S(v_{3}',(a(i)v,w)v_{2})\numberthis\label{n3.18}\\
&=\<v_{3}',f_{v}(o(a)v_{2})\>z^{-\wt a}w^{-\deg w}-\sum_{i,j\geq 0}\binom{\wt a-2-j}{i}w^{\deg v-j-1}z^{-\wt a+1+j}\<v_{3}',f_{a(i)v} (v_{2})\>,
\end{align*}
for any $a\in V$, $v_{3}'\in M^{3}(0)^{\ast}$, and $v\in M^{1}$.
By comparing the coefficients of $z^{-\wt a}w^{-\deg w}$ on both sides of \eqref{n3.18}, we have $\<v_{3}',f_{v}(o(a)v_{2})\>=0$ for all $v_{3}\in M^{3}(0)^{\ast}$, $a\in V$, and $\ v\in M^{1}$. Then $f_{v}(M^2(0))=0$, since $M^{2}(0)$ is an irreducible $A(V)$-module, and $M^{2}(0)=A(V).v_{2}=span\{o(a)v_{2}| a\in V\}$. It follows that $f_{v}=0$ for all $v\in M^{1}$. By Lemma \ref{lm3.3}, we have $S=0$, which is a contradiction. 

$(c)$ Let $x=(b_{m},i_{m})\ds(b_{1},i_{1})v_{2}$, with $\sum_{k=1}^{m}(\wt b_{k}-i_{k}-1)<0$. We use induction on the length $m$ of $x$ to show that $x\in \Rad(S)$. For the base case, let $x=(b,t)v_{2}$ with $\wt b-t-1<0$, then by \eqref{n3.13} and \eqref{n3.6}, we have
\begin{align*}
&S(v_{3}',(a_{1},z_{1})...(v,w)x)=\int_{C}	S(v_{3}',(a_{1},z_{1})...(v,w)(b,z)v_{2})z^{t}dz\\
&=\int_{C} S(v_{3}',(a_{1},z_{1})...(a_{n},z_{n})(v,w)o(b)v_{2})z^{t-\wt b}dz\numberthis \label{n3.19}\\
&\ \ \ +\int_{C}\sum_{k=1}^{n}\sum_{i\geq 0}G_{\wt b,i}(z,z_{k})S(v_{3}',(a_{1},z_{1})...(b(i)a_{k},z_{k})...(v,w)v_{2})z^{t}dz\\
&\ \ \ +\int_{C} \sum_{i\geq 0} G_{\wt b,i}(z,w)S(v_{3}',(a_{1},z_{1})...(b(i)v,w)v_{2})z^{t}dz,
\end{align*}
where $C$ is a contour of $z$ surrounding $0$, with all other variables lying outside $C$.  In particular, we have $|z|<|z_{k}|$ for all $k$, and $|z|<|w|$. Then by \eqref{n3.7},
\begin{equation}\label{n3.20}
\int_{C}G_{\wt b,i}(z,z_{k})z^{t}dz=\int_{C}\frac{z^{-\wt b+1+t}}{i!}\bigg(\frac{d}{dz_{k}}\bigg)^{i}\bigg(\frac{z_{k}^{\wt b-1}}{z-z_{k}}\bigg)dz=0,
\end{equation}
since $-\wt b+1+t>0$, and $1/(z-z_{k})$ is a sum of nonnegative powers in $z$ for all $z$ lying on the contour $C$. We also have $\int_{C}z^{t-\wt b}dz=0$, since $t-\wt b>-1$. It follows that all the integrals on the right-hand side of \eqref{n3.19} are equal to $0$. This finishes the base case.

Now let $m>0$, and consider $x=(b_{m},i_{m})...(b_{1},i_{1})v_{2}\in \bar{M}$. We have:
\begin{align*}
&S(v_{3}',(a_{1},z_{1})...(v,w)x)\\
&=\int_{C_{m}}...\int_{C_{1}}S(v_{3}',(a_{1},z_{1})...(v,w)(b_{m},w_{m})...(b_{1},w_{1})v_{2})w_{m}^{i_{m}}...w_{1}^{i_{1}}dw_{1}...dw_{m}\\
&=\int_{C_{m}}...\int_{C_{1}}\underset{(1)}{S(v_{3}',(a_{1},z_{1})...(v,w)(b_{m},w_{m})...o(b_{1})v_{2})w_{m}^{i_{m}}...w_{1}^{-\wt b_{1}+i_{1}}dw_{1}...dw_{m}}\\
&\ \ +\int_{C_{m}}...\int_{C_{1}}\sum_{k=1}^{n}\sum_{i\geq 0}\underset{(2)}{G_{\wt b_{1},i}(w_{1},z_{k})S(v_{3}',...(b_{1}(i)a_{k},z_{k})...(v,w)...v_{2})w_{m}^{i_{m}}...w_{1}^{i_{1}}dw_{1}...dw_{m}}\\
&\ \ +\int_{C_{m}}...\int_{C_{1}}\sum_{i\geq 0}\underset{(3)}{G_{\wt b_{1},i}(w_{1},w)S(v_{3}',...(b_{1}(i)v,w)(b_{m},w_{m})...v_{2})w_{m}^{i_{m}}...w_{1}^{i_{1}}dw_{1}...dw_{m}}\\
&\ \ +\int_{C_{m}}...\int_{C_{1}}\sum_{l=2}^{m}\sum_{i\geq 0}\underset{(4)}{G_{\wt b_{1},i}(w_{1},w_{l})S(v_{3}',...(v,w)...(b_{1}(i)b_{l},w_{l})...v_{2})w_{m}^{i_{m}}...w_{1}^{i_{1}}dw_{1}...dw_{m}}\\
&=(1)+(2)+(3)+(4),
\end{align*}
where $C_{1}$ is a contour of $w_{1}$ surrounding $0$, with all other variables lying outside. We need to show that the sum of these integrals equals $0$. i.e.,  $(1)+(2)+(3)+(4)=0$.

Case 1. $\wt b_{1}-i_{1}-1<0$.

Similar to \eqref{n3.20}, we have 
$\int_{C_{1}}G_{\wt b_{1},i}(w_{1},z)w_{1}^{i_{1}}dw_{1}=0,$ for $z=z_{k}$, $w$ or $w_{l}$.	Thus we have $(2)=(3)=(4)=0$. We also have $(1)=0$ because $-\wt b_{1}+i_{1}>-1$.

Case 2. $\wt b_{1}-i_{1}-1>0$.

Then $-\wt b_{1}+i_{1}<-1$, which implies $(1)=0$. Moreover, by \eqref{n3.7} we have: 
\begin{align*}
\int_{C_{1}}G_{\wt b_{1},i}(w_{1},z)w_{1}^{i_{1}}dw_{1}	&=\textrm{Res}_{w_{1}=0}\bigg(-\sum_{j\geq 0} \binom{\wt b_{1}-2-j}{i}z^{\wt b_{1}-j-2-i}w_{1}^{-\wt b_{1}+1+j+i_{1}}\bigg)\\
&=-\binom{i_{1}}{i}z^{i_{1}-i}.\numberthis\label{n3.21}
\end{align*}
for $z=z_{k}$, $w$ or $w_{l}$. Apply \eqref{n3.21} to (2), (3), and (4), and we have:
\begin{align*}
(2)&=-\int_{C_{m}}...\int_{C_{2}}\sum_{k=1}^{n}\sum_{i\geq 0}\binom{i_{1}}{i}z_{k}^{i_{1}-i}S(v_{3}',...(b_{1}(i)a_{k},z_{k})...(v,w)(b_{m},w_{m})...(b_{2},w_{2})v_{2})\\
&=-\sum_{k=1}^{n}\sum_{i\geq 0} \binom{i_{1}}{i}z_{k}^{i_{1}-i}S(v_{3}',(a_{1},z_{1})...(b_{1}(i)a_{k},z_{k})...(a_{n},z_{n})(v,w)y),
\end{align*}
where $y=(b_{m},i_{m})\ds(b_{2},i_{2})v_{2}$. Note that $\deg y=\deg x-(\wt b_{1}-i_{1}-1)<0$, and the length of $y$ is $m-1$, then by the induction hypothesis we have $(2)=0$. Similarly, $(3)=0$.
\begin{align*}
(4)&=\int_{C_{m}}...\int_{C_{1}}\sum_{l=2}^{m}\sum_{i\geq 0}
\binom{i_{1}}{i}w_{l}^{i_{1}-i}S(v_{3}',...(v,w)...(b_{1}(i)b_{l},w_{l})...v_{2})w_{m}^{i_{m}}...w_{1}^{i_{1}}dw_{1}...dw_{m}\\
&=\sum_{l=2}^{m}\sum_{i\geq 0}\binom{i_{1}}{i}S(v_{3}',(a_{1},z_{1})...(a_{n},z_{n})(v,w)y_{l}),
\end{align*}  
where $y_{l}=(b_{m},i_{m})...(b_{1}(i)b_{l},i_{1}+i_{l}-i)...(b_{2},i_{2})v_{2}$. 
Note that $$\deg(b_{1}(i)b_{l},i_{1}+i_{l}-i)=\wt b_{1}+\wt b_{l}-i-1-i_{1}-i_{l}+i-1=\deg(b_{1},i_{1})+\deg(b_{l},i_{l}).$$
Thus, $\deg y_{l}=\sum_{k=1}^{m}\wt(b_{k},i_{k})=\deg x<0$, and the length of $y_{l}$ is $m-1$ for each $l$. Hence $(4)=0$ by the induction hypothesis.

Case 3. $\wt b_{1}-i_{1}-1=0$.

In this case, we have:
$\int_{C_{1}}G_{\wt b_{1},i}(w_{1},z)w_{1}^{i_{1}}dw_{1}=0$ in view of \eqref{n3.20}. Hence $(2)=(3)=(4)=0$. Moreover, since $-\wt b_{1}+i_{1}=-1$, we have:
\begin{align*}
(1)&=\int_{C_{m}}...\int_{C_{2}}S(v_{3}',(a_{1},z_{1})...(v,w)(b_{m},w_{m})...o(b_{1})v_{2})w_{m}^{i_{m}}...w_{2}^{i_{2}}dw_{2}...dw_{m}\\
&=S(v_{3}',(a_{1},z_{1})...(a_{n},z_{n})(v,w)y),
\end{align*}
where $y=(b_{m},i_{m})...(b_{2},i_{2})v_{2}$. Since $\deg y=\deg x<0$, and the length of $y$ is $m-1$, we have $(1)=0$ by the induction hypothesis. Now the proof of $(c)$ is complete.
\end{proof}

We define a vertex operator $Y_{\bar{M^2}}$ on the quotient space $\bar{M^2}=\bar{M}/\Rad(\bar{M})$ as follows:
\begin{equation}\label{n3.22}
Y_{\bar{M^2}}(a,z)(b_{1},i_{1})\ds(b_{m},i_{m})v_{2}:=\sum_{n\in\mathbb{Z}}(a,n)(b_{1},i_{1})\ds(b_{m},i_{m})v_{2}z^{-n-1},
\end{equation}
where $a\in V, \ (b_{1},i_{1})\ds(b_{m},i_{m})v_{2}\in \bar{M^2}$, and we use the same notation $(b_{1},i_{1})\ds (b_{m},i_{m})v_{2}$ for its image in the quotient space $\bar{M^2}$. We can express \eqref{n3.22} in the component form:
\begin{equation}\label{n3.23}
a(n)(b_{1},i_{1})\ds(b_{m},i_{m})v_{2}=(a,n)(b_{1},i_{1})\ds(b_{m},i_{m})v_{2},
\end{equation}
for all $a\in V,\ n\in \mathbb{Z}$, and $(b_{1},i_{1})\ds(b_{m},i_{m})v_{2}\in \bar{M}$. 

\begin{prop}\label{prop3.6}
$\bar{M^2}=\bar{M}/\Rad(\bar{M})$, together with $Y_{\bar{M^2}}: V\ra \End(\bar{M^2})[[z,z^{-1}]]$ given by \eqref{n3.22} and \eqref{n3.23}, is a weak $V$-module.
\end{prop}
\begin{proof}
By (a) of Lemma \ref{lm3.5}, we have $a(n)\Rad(\bar{M})\subseteq \Rad(\bar{M})$. Hence $Y_{\bar{M^2}}$ is well-defined. Let $x=(b_{1},i_{1})\ds (b_{m},i_{m})v_{2}\in \bar{M^2}$, we claim that $\mathbf{1}(-1)x=x$ and $\mathbf{1}(n)x=0$ for any $n\neq -1$. Indeed, for any $S\in \Cor\fusion{M^{1}}{M^{2}(0)}{M^{3}(0)}$, by the definition formula \eqref{n3.13}, the recursive formula \eqref{n3.6}, together with the fact that $\mathbf{1}(j)a=0$ for all $j\geq0$, $a\in V$, and $\vac(j)v=0$ for all $j\geq 0$, $v\in M^{1}$, we have:
\begin{align*}
&S(v_{3}',(a_{1},z_{1})...(v,w)\vac (n)x)\\
&= \int_{C_{0}}\int_{C_{m}}...\int_{C_{1}}S(v_{3}',(\mathbf{1},w_{0})(a_{1},z_{1})...(v,w)(b_{1},w_{1})...v_{2}) w_{0}^{n}w_{1}^{i_{1}}...w_{m}^{i_{m}}dw_{1}...dw_{m}dw_{0}\\
&=\int_{C_{0}}\int_{C_{m}}...\int_{C_{1}}S(v_{3}'o(\mathbf{1}),(a_{1},z_{1})...(v,w)(b_{1},w_{1})...v_{2})w_{0}^{n}w_{1}^{i_{1}}...w_{m}^{i_{m}}dw_{1}...dw_{m}dw_{0}\\
&=\delta_{n+1,0} \cdot S(v_{3}',(a_{1},z_{1})...(v,w)x),
\end{align*}
where the last equality follows from the fact that $\int_{C_{0}}w_{0}^{n}dw_{0}=\delta_{n+1,0}$. Thus, $(\vac(n)x-\delta_{n+1,0} x)\in \Rad(\bar{M})$, and so $\vac(n)x=\delta_{n+1,0} x$ in $\bar{M^2}$. Moreover, given homogeneous elements $x\in \bar{M}$ and $a\in V$, by \eqref{n3.14} and \eqref{n3.23}, $\deg (a(n).x)=\wt a-n-1+\deg x<0$ when $n>>0$. Then by part (c) of Lemma \ref{lm3.5}, we have $a(n)x=0$ in $\bar{M^2}$ when $n$ is large enough. Finally, by Lemma \ref{lm3.4} and \eqref{n3.23}, $(\bar{M^2},Y_{\bar{M^2}})$ satisfies the Jacobi identity. Hence it is a weak $V$-module.
\end{proof}

\begin{prop}\label{prop3.7}
$\bar{M^2}$ has a gradation $\bar{M^2}=\bigoplus_{n=0}^\infty\bar{M^2}(n)$, where $\bar{M^2}(n)$ is an eigenspace of $L(0)$ of eigenvalue $\la+n$ for each $n\in \N$, and $\bar{M^2}(0)=M^{2}(0)$. In particular, $\bar{M^2}$ is an ordinary $V$-module, and if $M^2(0)$ is the bottom level of some ordinary $V$-module $M^2$, with conformal weight $h_{2}$, then $\la=h_2$. 
\end{prop}
\begin{proof}
Let $\bar{M^2}(n)$ be the image of $\bar{M}(n)$ under the quotient map $\bar{M}\ra \bar{M^{2}}$. By Lemma \ref{lm3.5}, we have $\bar{M^2}=\sum_{n\geq 0}\bar{M^2}(n)$ and $ M^{2}(0)\subseteq \bar{M^2}(0)$. We claim that 
\begin{equation}\label{n3.24}
a(\wt a-1)v_{2}=o(a)v_{2},
\end{equation}
for all $v_{2}\in M^{2}(0)$ and homogeneous $a\in V$.
Indeed, we only need to show that $(a,\wt a-1)v_{2}-o(a)v_{2}\in \Rad(S)$, for all $S\in \Cor\fusion{M^{1}}{M^{2}(0)}{M^{3}(0)}$. By \eqref{n3.13} and \eqref{n3.6}, 
\begin{align*}
&S(v_{3}',(a_{1},z_{1})...(a_{n},z_{n})(v,w)(a,\wt a-1)v_{2})\\&=\int_{C}S(v_{3}',(a_{1},z_{1})...(a_{n},z_{n})(v,w)(a,w_{1})v_{2})w_{1}^{\wt a-1}dw_{1}\\
&=\int_{C}S(v_{3}',(a_{1},z_{1})...(a_{n},z_{n})(v,w)o(a)v_{2})w_{1}^{-\wt a}w_{1}^{\wt a-1}dw_{1}\\
& \ \ \ +\sum_{k=1}^{n}\sum_{i\geq 0}\int_{C}G_{\wt a,i}(w_{1},z_{k})S(v_{3}',(a_{1},z_{1})...(a(i)a_{k},z_{k})...(a_{n},z_{n})(v,w)v_{2})w_{1}^{\wt a-1}dw_{1}\\
&\ \ \ +\sum_{i\geq 0}\int_{C}G_{\wt a,i}(w_{1},w)S(v_{3}',(a_{1},z_{1})...(a_{n},z_{n})(a(i)v,w)v_{2})w_{1}^{\wt a-1}dw_{1},
\end{align*}
where $C$ is a contour of $w_{1}$ surrounding $0$, with all other variables lying outside of $C$. Since $|z_k|, |w|>|w_{1}|$ for all $k$, where $w_{1}$ is lying on $C$, then we have 
$$\int_{C} G_{\wt a,i}(w_{1},z)w_{1}^{\wt a-1}dw_{1}=\int_{C}w_{1}^{\wt a-1}\frac{w_{1}^{-\wt a+1}}{i!}\bigg(\frac{d}{dz}\bigg)^{i}\bigg(\frac{z^{\wt a-1}}{w_{1}-z}\bigg)dw_{1}=0,$$
for $z=z_k$ or $w$. Hence $(a,\wt a-1)v_{2}-o(a)v_{2}\in \Rad(S)$. This shows \eqref{n3.24}.

Since $L(0)=\omega(\wt \omega-1)$ on $\bar{M^2}$,  it follows from \eqref{n3.24} that $L(0)$ preserves $M^{2}(0)$. On the other hand, we have $[L(0), a(n)]=(\wt a-n-1)a(n)$ (see (4.2.2) in \cite{FHL}). Then by \eqref{n3.24} again, we have $[L(0),o(a)]v_{2}=[L(0),a(\wt a-1)]v_{2}=0$. Since $M^{2}(0)$ is an irreducible $A(V)$-module which is of countable dimension, then by the Schur's Lemma (Lemma 1.2.1 in \cite{Z}), there exists $\la\in \C$ such that $L(0)=\lambda \cdot \Id$ on $M^{2}(0)$. If $M^2(0)$ is the bottom level of $M^2$, with conformal weight $h_{2}$, then $L(0)=h_{2}\cdot \Id$ on $M^2(0)$, and so $h_{2}=\la$.

Now for any spanning element $x=(b_{1},i_{1})\ds(b_{m},i_{m})v_{2}=b_{1}(i_{1})\ds b_{m}(i_{m})v_{2}$ of $ \bar{M^2}(n)$, we have 
$L(0)x=(\sum_{k=1}^{m}(\wt b_{k}-i_k-1)+\lambda )x=(n+\lambda)x.$
Therefore, $\bar{M^2}(n)$ is an eigenspace of $L(0)$ of eigenvalue $n+\lambda$ for every $n\in \N$, and $\bar{M^2}=\bigoplus_{n=0}^\infty\bar{M^2}(n)$.

Finally, for any spanning element $x=b_{1}(i_{1})\ds b_{m}(i_{m})v_{2}$ of $ \bar{M^2}(0)$, it follows from \eqref{n3.24} and an easy induction that $x\in M^{2}(0)$, therefore $\bar{M^2}(0)=M^{2}(0)$.
\end{proof}

\begin{remark}Unlike the construction of $V$-modules from the correlation functions in Theorem 2.2.1 in \cite{Z}, in our case, it is unclear whether $\bar{M^2}=\bar{M}/\Rad(\bar{M})$ is an irreducible $V$-module. The reason is the following: 

Assume $N\leq \bar{M^2}$ is a submodule, by Proposition \ref{prop3.7} we have $N=\bigoplus_{n=0}^{\infty} N(n)$, with $N(n)=N\cap \bar{M^2}(n)$ for each $n$. If $N(0)\neq 0$, then clearly $N=\bar{M^2}$. So to show $\bar{M^2}$ is irreducible, we need to show that $N=0$ when $N(0)=0$. 

This is true for the module $\bar{M}/\Rad(\bar{M})$ constructed in Theorem 2.2.1 in \cite{Z}, wherein the correlation function $S(v',(a_{1},z_{1})\ds (a_{n},z_{n})N)$, with $v'\in M^2(0)$, is essentially the limit function of $\<v', Y(a_{1},z_{1})\dots Y(a_{n},z_{n})N\>$.
It is zero because $Y(a,z)N\subset N((z))$, and the bottom level of $N$ is $0$. Thus, $N\ssq \Rad(S)$, and so $N=0$ in $\bar{M}/\Rad(\bar{M})$. However, in our case, $S(v_{3}', (a_{1},z_{1})\dots (a_{n},z_{n})(v,w)N)$ with $v_{3}'\in M^{3}(0)^{\ast}$ is essentially the limit function of 
$\<v_{3}',I(v,w)Y(a_{1},z_{1})\dots Y(a_{n},z_{n})N\>w^{-h}.$
Although the components of $Y(a,z)$ still leave $N$ invariant, the intertwining operator $I(v,w)$ could send some element in $N$ to a nonzero element of $M^{3}(0)$. Hence we cannot conclude that $N\ssq \Rad (\bar{M})$ in general. 
\end{remark}
We give a sufficient condition under which $\bar{M^2}$ is irreducible.


\begin{lm}\label{lm3.9}
	Suppose $S\in \Cor\fusion{M^1}{M^2(0)}{M^3(0)}$ satisfies:
	\begin{equation}\label{n3.25}
	\sum_{i\geq 0}\binom{n}{i} \<v_{3}',f_{b(i)v}(v_{2})\>=0,
	\end{equation}
	for all $b\in V$, $n\in \Z$ such that $\wt b-n-1>0$, $v\in M^{1}$, $v_{3}'\in M^{3}(0)^{\ast}$, and $v_{2}\in M^{2}(0)$.
	Then $S(v_{3}',(v,w)y)=0$ for any $y\in M(m)$ with $m\geq 1$, $v_{3}'\in M^{3}(0)^{*}$, and $v\in M^{1}$.
\end{lm}
\begin{proof}
	It follows from an easy induction that $y$ can be written as a sum of the terms $(b_{m},n_{m})\ds(b_{1},n_{1})v_{2}$ for some $m\geq 1$ and $v_{2}\in M^{2}(0)$, with $\wt b_{j}-n_{j}-1>0$ for all $j$. 
	
	Let $y= (b_{m},n_{m})\ds (b_{1},n_{1})v_{2}$. We use induction on $m$ to show that $S(v_{3}',(v,w)y)=0$. For the base case $m=1$ and $y=(b,n)v_{2}$, with $\wt b-n-1>0$, by \eqref{n3.13}, \eqref{n3.3}, \eqref{n3.6}, \eqref{n3.7}, and the assumption \eqref{n3.25}, we have:
	\begin{align*}
	&S(v_{3}', (v,w)y)=\int_{C}S(v_{3}', (v,w)(b,z)v_{2})z^{n}dz\\
	&=\int_{C}S(v_{3}', (v,w)o(b)v_{2})z^{-\wt b+n}dz+\int_{C}\sum_{i\geq 0}G_{\wt b,i}(z,w)S(v_{3}',(b(i)v,w)v_{2})z^{n}dz\\
	&=0+ \sum_{i\geq 0}\int_{C} -\sum_{j
		\geq 0}\binom{\wt b-2-j}{i} w^{\wt b-j-2-i} z^{n-\wt b+1+j} S(v_{3}', (b(i)v,w)v_{2})dz\\
	&=-\sum_{i\geq 0}\binom{n}{i} \<v_{3}',f_{b(i)v}v_{2}\> w^{-\wt b-\deg v+1+n}=0.
	\end{align*}
	Now let $m>1$. Then by \eqref{n3.13} and \eqref{n3.6}, we have 
	\begin{align*}
	&S(v_{3}',(v,w)y)=\int_{C_{m}}...\int_{C_{1}}S(v_{3}',(v,w)(b_{m},z_{m})...(b_{1},z_{1})v_{2})z_{1}^{n_{1}}...z_{m}^{n_{m}}dz_{1}...dz_{m}\\
	&=\int_{C_{m}}...\int_{C_{1}}S(v_{3}',(v,w)(b_{m},z_{m})...(b_{2},z_{2})o(b_{1})v_{2})z_{1}^{-\wt b_{1}+n_{1}}...z_{m}^{n_{m}}dz_{1}...dz_{m}\\
	&\ +\int_{C_{m}}...\int_{C_{1}}\sum_{i\geq 0}G_{\wt b_{1},i}(z_{1},w)S(v_{3}',(b_{1}(i)v,w)(b_{m},z_{m})...(b_{2},z_{2})v_{2})z_{1}^{n_{1}}...z_{m}^{n_{m}}dz_{1}...dz_{m}\\
	&\ +\int_{C_{m}}...\int_{C_{1}}\sum_{k=2}^{m}\sum_{i\geq 0}G_{\wt b_{1},i}(z_{1},z_{k})S(v_{3}',(v,w)...(b_{1}(i)b_{k},z_{k})...(b_{2},z_{2})v_{2}) z_{1}^{n_{1}}...z_{m}^{n_{m}}dz_{1}...dz_{m}\\
	&=0+\int_{C_{m}}...\int_{C_{2}}\sum_{i\geq 0} \int_{C_{1}}-\sum_{j\geq 0}\binom{\wt b_{1}-2-j}{i} w^{\wt b_{1}-j-2-i} z_{1}^{n_{1}-\wt b_{1}+1+j} \\
	&\qquad \cdot S(v_{3}', (b_{1}(i)v,w)(b_{m},z_{m})...(b_{2},z_{2})v_{2})z_{2}^{n_{2}}...z_{m}^{n_{m}}dz_{2}...dz_{m}\\
	&\ +\int_{C_{m}}...\int_{C_{2}}\sum_{k=2}^{m}\sum_{i\geq 0} \int_{C_{1}}-\sum_{j\geq 0}\binom{\wt b_{1}-2-j}{i} z_{k}^{n_{k}+\wt b_{1}-j-2-i} z_{1}^{n_{1}-\wt b_{1}+1+j} \\
	&\qquad \cdot S(v_{3}', (v,w)(b_{m},z_{m})...(b_{1}(i)b_{k},z_{k})...(b_{2},z_{2})v_{2})z_{2}^{n_{2}}...\widehat{z_{k}^{n_{k}}}...z_{m}^{n_{m}}dz_{2}...dz_{m}\\
	&=-\int_{C_{m}}...\int_{C_{2}}\sum_{i\geq 0} \binom{n_{1}}{i} w^{n_{1}-i} S(v_{3}',(b_{1}(i)v,w)(b_{m},z_{m})...(b_{2},z_{2})v_{2})z_{2}^{n_{2}}...z_{k}^{n_{k}}dz_{2}...dz_{m}\\
	&\ -\int_{C_{m}}...\int_{C_{2}} \sum_{k=2}^{m} \sum_{i\geq 0} \binom{n_{1}}{i} S(v_{3}', (v,w)...(b_{1}(i)b_{k},z_{k})...v_{2}) z_{2}^{n_2}...z_{k}^{n_{1}+n_{k}-i}...z_{m}^{n_{m}}dz_{2}...dz_{m}\\
	&=-\sum_{i\geq 0} \binom{n_{1}}{i} w^{n_{1}-i}S(v_{3}', (b_{1}(i)v,w)(b_{m},n_{m})...(b_{2},n_{2})v_{2})\\
	&\ -\sum_{k=2}^{m}\sum_{i\geq 0} \binom{n_{1}}{i} S(v_{3}',(v,w)(b_{m},n_{m})...(b_{1}(i)b_{k},n_{1}+n_{k}-i)...(b_{2},n_{2})v_{2})\\
	&=0,
	\end{align*}
	where the last equality follows from the induction hypothesis, together with the fact that $\deg (b_{1}(i)b_{k}, n_{1}+n_{k}-i)=\wt b_{1}-n_{1}-1+\wt b_{k}-n_{k}-1>0$, for any $i\geq 0$.
\end{proof}

\begin{coro}\label{coro3.11}
	For any fixed $v\in M^{1}$ and $y\in \bar{M^{2}}=\bar{M}/\Rad(\bar{M})$, let $n\in \Z$ be an integer such that $n>\deg v+\deg y-1$. Then we have
	\begin{equation}\label{n3.26}
	\int_{C}S(v_{3}', (v,w)y)w^{n}dw=0,
	\end{equation}
	for all $v_{3}'\in M^{3}(0)$, where $C$ is a contour of $w$ surrounding $0$. In particular, for fixed $v\in M^{1}$ and $y\in \bar{M^{2}}$, the power series expansion of $S(v_{3}', (v,w)y)$ has a uniform lower bound for $w$ independent of $v_{3}'\in M^{3}(0)^{\ast}$.
\end{coro}
\begin{proof}
	It suffices to show \eqref{n3.26} for $y=(b_{m},n_{m})\ds(b_{1},n_{1})v_{2}$, where $v_{2}\in M^{2}(0)$, $m\geq 0$, and $\wt b_{j}-n_{j}-1>0$ for all $j$. Again, we use induction on $m$. When $m=0$, we have $y=v_{2}$ and $\deg y=0$. Then by \eqref{n3.3} and $-\deg v+n>-1$, we have: 
	$\int_{C} S(v_{3}', (v,w)v_{2})w^{n}dw=\int_{C} \<v_{3}', f_{v}(v_{2})\> w^{-\deg v+n}dw=0$. 
	Now let $m>0$, and let $n\in \Z $ be such that $n>\deg v+\deg y-1$. Since $-\wt b_{1}+n_{1}<-1$, by the calculations in Lemma \ref{lm3.9}, we have:
	\begin{align*}
	&\int_{C}S(v_{3}', (v,w)y)w^{n}dw=-\sum_{i\geq 0} \int_{C}\binom{n_{1}}{i} w^{n+n_{1}-i}\underset{(1)}{S(v_{3}', (b_{1}(i)v,w)(b_{m},n_{m})...(b_{2},n_{2})v_{2})}dw\\
	&\ \ \ -\sum_{k=2}^{m}\sum_{i\geq 0} \int_{C}\binom{n_{1}}{i}\underset{(2)}{w^{n}S(v_{3}',(v,w)(b_{m},n_{m})...(b_{1}(i)b_{k},n_{1}+n_{k}-i)...(b_{2},n_{2})v_{2})}dw\\
	&=(1)+(2).
	\end{align*}
	Since $n>\deg v+\deg y-1$, we have $n+n_{1}-i>\deg (b_{1}(i)v)+\sum_{j=2}^{m} (\wt b_{j}-n_{j}-1)-1$ for all $i\geq 0$. Then by the induction hypothesis, $(1)=0$ for all $v_{3}'\in M^{3}(0)^{\ast}$.
	On the other hand, since $\deg (b_{1}(i)b_{k}, n_{1}+n_{k}-i)=\wt b_{1}-n_{1}-1+\wt b_{k}-n_{k}-1$ for all $i\geq 0$, we have $(2)=0$ for all $v_{3}'\in M^{3}(0)^{\ast}$. Thus $\int_{C}S(v_{3}', (v,w)y)w^{n}dw=0$.
\end{proof}

\begin{prop}\label{prop3.11}
	Suppose every $S\in \Cor\fusion{M^1}{M^2(0)}{M^3(0)}$ satisfies the condition \eqref{n3.25}, then
	$\bar{M^2}=\bar{M}/\Rad(\bar{M})$ is an irreducible $V$-module with bottom level $M^{2}(0)$. In particular, $\bar{M^2}$ is isomorphic to $L(M^{2}(0))$, the unique irreducible $V$-module with bottom level is $M^{2}(0)$.
\end{prop}
\begin{proof}
	Note that for any $x\in M$, $S(v_{3}',(a_{1},z_{1})\ds (a_{n},z_{n})(v,w)x)$ is also a rational function in $z_{1},\ds,z_{n},w$ by \eqref{n3.13} and \eqref{n3.23}, and it has Laurent series expansion: 
	\begin{align*}
	&S(v_{3}',(a_{1},z_{1})...(a_{n},z_{n})(v,w)x)=S(v_{3}',(v,w)(a_{1},z_{1})...(a_{n},z_{n})x)\\
	&=\sum_{i_{1},...,i_{n}\in\mathbb{Z}} \left(\int_{C_{n}}...\int_{C_{1}}S(v_{3}',(v,w)(a_{n},z_{n})...(a_{1},z_{1})x)z_{1}^{i_{1}}...z_{n}^{i_{n}}dz_{1}...dz_{n}\right)z_{1}^{-i_{1}-1}...z_{n}^{-i_{n}-1}\\
	&=\sum_{i_{1},...,i_{n}\in\mathbb{Z}}S(v_{3}',(v,w)a_{n}(i_{n})...a_{1}(i_{1})x)z_{1}^{-i_{1}-1}...z_{n}^{-i_{n}-1}\numberthis \label{n3.27}
	\end{align*}
	on the domain 
	$\mathbb{D}=\{(z_{1},\ds,z_{n},w)||w|>|z_{n}|>\ds>|z_{1}|>0\}$.
	Let $N$ be a submodule of $\bar{M^2}$ such that $N(0)=0$, we need to show that $N=0$. Let $x\in N$, we have $y=a_{n}(i_{n})\ds a_{1}(i_{1})x\in N$, and if $y\neq 0$ then $\deg(y)>0$. By Lemma \ref{lm3.9}, we have $S(v_{3}',(v,w)y)=0$. Thus, the rational function $S(v_{3}',(a_{1},z_{1})\ds (a_{n},z_{n})(v,w)x)$ is equal to $0$ by \eqref{n3.27}. i.e., $x\in \Rad(S)$ for all $S\in \Cor\fusion{M^1}{M^2}{M^3}$. Thus $N=0$.
\end{proof}

In conclusion, given a $S\in \Cor\fusion{M^1}{M^2(0)}{M^3(0)}$, the extended $S$ in \eqref{n3.13} factors though an $\N$-gradable $V$-module $\bar{M^2}=\bar{M}/\Rad(\bar{M})$ whose bottom level is $M^2(0)$. It is irreducible if the condition \eqref{n3.25} is satisfied. Therefore, by \eqref{n3.13} and \eqref{n3.23}, we have a well-defined system of $(n+3)$-point correlation functions:
\begin{align*}
&S:M^{3}(0)^{*}\times V\times \ds\times M^{1}\times \ds\times V\times \bar{M^2}\rightarrow \mathcal{F}(z_{1},\ds,z_{n},w),\\
&S(v_{3}',(a_{1},z_{1})...(a_{n},z_{n})(v,w)b_{1}(i_{1})...b_{m}(i_{m})v_{2}) \numberthis \label{n3.28}\\
&=\int_{C_{1}}...\int_{C_{m}}S(v_{3}',(a_{1},z_{1})...(a_{n},z_{n})(v,w)(b_{1},w_{1})...(b_{m},w_{m})v_{2})w_{1}^{i_{1}}...w_{m}^{i_{m}}dw_{1}...dw_{m},
\end{align*}
for all $b_{1}(i_{1})\ds b_{m}(i_{m})v_{2}\in \bar{M^2}$, where $C_{k}$ is a contour of $w_{k}$, $C_{k}$ contains $C_{k+1}$ for all $k$, $C_{m}$ contains $0$, and $z_{1},\ds,z_{n}, w$ are outside of $C_{1}$. 

In particular, $S$ in \eqref{n3.28} satisfies the generating formula \eqref{n2.9} with $M^2=\bar{M^2}$, since the extended $S$ is defined by this formula. Moreover, by Corollary \ref{coro3.11} and the fact that the orginal $S$ in \eqref{n3.11} belongs to $\Cor\fusion{M^1}{M^2(0)}{M^3(0)}$, it is easy to see that the $S$ in \eqref{n3.28} also satisfies the properties $(1)-(6)$ in Definition \ref{df2.1}, with $v_{3}'\in M^3(0)^\ast$ and $v_{2}\in \bar{M^2}$. 

We adopt a similar method to extend the first input component of $S$ in \eqref{n3.28} from $M^{3}(0)^{\ast}$ to a $V$-module by using the other generating formula \eqref{n2.10}. First, we let
$$\tilde{M}:=T(\mathcal{L}(V))\otimes_{\C} M^{3}(0)^{\ast}.$$ 
Then $\tilde{M}$ is spanned by elements of the form:
$y=(b_{1},i_{1})\ds (b_{m},i_{m})v_{3}',$
where $b_{j}\in V$, $i_{j}\in \mathbb{Z}$ for $j=1,\ds ,m$, and $v_{3}'\in M^{3}(0)^{*}$. Next, we extend $S$ in \eqref{n3.28} by iterating the generating formula \eqref{n2.10}. i.e., we define:
\begin{align*}
&S: \tilde{M} \times V\times \ds\times M^{1}\times\ds\times V\times \bar{M^2}\ra \mathcal{F}(z_{1},\ds,z_{n},w)\\
&S((b_{1},i_{1})...(b_{m},i_{m})v_{3}',(a_{1},z_{1})...(a_{n},z_{n})(v,w)x_{2})\numberthis \label{n3.29}\\
&:=\int_{C_{1}}...\int_{C_{m}}S(v_{3}',(b_{m},w_{m})'...(b_{1}, w_{1})'(a_{1},z_{1})...(v,w)x_{2})w_{1}^{-i_{1}-2}...w_{m}^{-i_{m}-2}dw_{m}...dw_{1},
\end{align*}
where $(b,w)'=(e^{w^{-1}L(1)}(-w^{2})^{L(0)}b,w),$ $C_{k}$ is a contour of $w_{k}$ s.t. $C_{k}$ contains $C_{k-1}$  for each $k$, and $C_{1}$ contains all the variables $z_{1},\ds,z_{n}, w$. For $S$ in \eqref{n3.29}, we similarly define $$\Rad(S):=\{y\in \tilde{M}: S(y,(a_{1},z_{1})\ds(a_{n},z_{n})(v,w)x)=0, \forall a_{i}\in V, v\in M^{1}, x\in \bar{M^2}\},$$
and let $\Rad (\tilde{M}):=\bigcap \Rad(S)$, where the intersection is taken over all $S\in \Cor\fusion{M^{1}}{M^{2}(0)}{M^{3}(0)}$, with the extension given by \eqref{n3.29}. Clearly, $S$ factors though $\tilde{M}/\Rad(\tilde{M})$.

Similar to our previous argument, one can show that $ \bar{M^3}'=\tilde{M}/\Rad(\tilde{M})$ has a natural $\N$-gradable $V$-module structure $\bar{M^3}'=\bigoplus_{n=0}^{\infty} \bar{M^3}'(n)$, with $\bar{M^3}'(0)=M^{3}(0)^{\ast}$. Moreover, $\bar{M^3}'=\tilde{M}/\Rad(\tilde{M})$ is irreducible if the condition \ref{n3.25} is satisfied. Thus we have a well-defined system of correlation functions $S$:

\begin{align*}
&S: \bar{M^3}'\times V\times\ds\times M^{1}\times\ds\times V\times \bar{M^2}\rightarrow \mathcal{F}(z_{1},\ds,z_{n},w), \\
&S(b_{1}(i_{1})...b_{m}(i_{m})v_{3}',(a_{1},z_{1})...(a_{n},z_{n})(v,w)x_{2})\numberthis \label{n3.30}\\
&=\int_{C_{1}}...\int_{C_{m}}S(v_{3}',(b_{m},w_{m})'...(b_{1}, w_{1})'(a_{1},z_{1})...(v,w)x_{2})w_{1}^{-i_{1}-2}...w_{m}^{-i_{m}-2}dw_{m}...dw_{1},
\end{align*}
for all $b_{1}(i_{1})\ds b_{m}(i_{m})v_{3}'\in \bar{M^3}'$ and $x_{2}\in \bar{M^2}$.
Moreover, by Remark \ref{rmk3.2}, we also have: 
\begin{align*}
&S(b_{1}'(i_{1})...b_{m}'(i_{m})v_{3}',(a_{1},z_{1})...(a_{n},z_{n})(v,w)x_{2})\numberthis \label{n3.31}\\
&=\int_{C_{1}}...\int_{C_{m}}S(v_{3}',(b_{m},w_{m})...(b_{1}, w_{1})(a_{1},z_{1})...(a_{n},z_{n})(v,w)x_{2})w_{1}^{i_{1}}...w_{m}^{i_{m}}dw_{m}...dw_{1},
\end{align*}
where $b'(i)=\sum_{j\geq 0}\frac{1}{j!}(-1)^{\wt b}(L(1)^{j}b)(2\wt b-i-j-2)$, $C_{k}$ is a contour of $w_{k}$ such that $C_{k}$ contains $C_{k-1}$ for each $k$, and $z_{1},\ds,z_{n},w$ are inside of $C_{1}$. Since \eqref{n3.30} and \eqref{n3.31} are given by iterating the generating formula \eqref{n2.10}, it is clear that $S$ in \eqref{n3.30} also satisfies \eqref{n2.10} with $M^2=\bar{M^2}$ and $M^{3'}=\bar{M^3}'$. Denote the contragredient module of $\bar{M^3}'$ by $\bar{M^3}$.

\begin{thm}\label{thm3.12}
	The system of extended correlation functions $S$ in \eqref{n3.30} lies in $\Cor\fusion{M^{1}}{\bar{M^2}}{\bar{M^3}}$. Hence we have an isomorphism of vector spaces $\Cor\fusion{M^{1}}{M^{2}(0)}{M^{3}(0)}\cong \Cor\fusion{M^{1}}{\bar{M^2}}{\bar{M^3}}\cong I\fusion{M^{1}}{\bar{M^2}}{\bar{M^3}}$. 
\end{thm}
\begin{proof}
	We have already proven that $S$ satisfies $(7)$ and $(8)$ in Definition \ref{df2.1}, with $M^2=\bar{M^2}$ and $M^{3'}=\bar{M^3}'$. It remains to show that $S$ in \eqref{n3.30} satisfies the properties $(1)-(6)$ in Definition \ref{df2.1}, with $M^2=\bar{M^2}$ and $M^3=\bar{M^3}$. In fact, by the definition formulas \eqref{n3.28} and \eqref{n3.31}, together with the fact that the orginal $S$ in \eqref{n3.11} lies in $\Cor\fusion{M^{1}}{M^{2}(0)}{M^{3}(0)}$, the properties $(2)-(6)$ are straightforward. 
	
	To prove (1), we need an intermediate result first. We introduce the following notation: 
	\begin{equation}\label{n3.32}
	\begin{aligned}
	S(v_{3}', b_{1}(n_{1})...b_{m}(n_{m})(v,w)x_{2}):&=\int_{C_{m}}...\int_{C_{1}}S(v_{3}',(b_{1},z_{1})...(b_{m},z_{m})(v,w)x_{2})\\
	&\qquad \cdot  z_{1}^{n_{1}}...z_{m}^{n_{m}}dz_{1}...dz_{m},
	\end{aligned}
	\end{equation}
	where $m\geq 0$, $x_{2}\in \bar{M^2}$, $b_{k}\in V$, $n_{k}\in \Z$, $C_{k}$ is a contour of $z_{k}$ s.t. $C_{k}$ contains $C_{k+1}$ for all $k$, and $w$ is inside of $C_{m}$. Assume $\wt b_{1}-n_{1}-1<0$. We claim that:
	
	\begin{align*}
	&S(v_{3}', b_{1}(n_{1})...b_{m}(n_{m})(v,w)x_{2})\\
	&=\sum_{l=2}^{m}\sum_{i\geq 0} \binom{n_{1}}{i}S(v_{3}', b_{2}(n_{2})...(b_{1}(i)b_{l})(n_{1}+n_{l}-i)...b_{m}(n_{m})(v,w)x_{2})\\
	&\ \ \ +\sum_{i\geq 0} \binom{n_{1}}{i}S(v_{3}', b_{2}(n_{2})...b_{m}(n_{m})(b_{1}(i)v,w)x_{2})w^{n_{1}-i}\numberthis \label{n3.33}\\
	&\ \ \ + S(v_{3}', b_{2}(n_{2})...b_{m}(n_{m})(v,w)(b_{1}(n_{1})x_{2})).
	\end{align*}
	Let $x_{2}=c_{1}(i_{1})\ds c_{r}(i_{r})v_{2}$, for some $c_{j}\in V$, $i_{j}\in \Z$ for all $j$, and $v_{2}\in M^{2}(0)$. Note that $b_{1}(n_{1})v_{2}=0$ as $\wt b_{1}-n_{1}-1<0$. For $|z_{1}|>|w|$, by \eqref{n3.5} we have:
	$$\int_{C_{1}} F_{\wt b_{1},i}(z_{1},w)z_{1}^{n_{1}}dz_{1}=\sum_{j\geq 0} \int_{C_{1}}\binom{\wt b_{1}+j}{i} z_{1}^{n_{1}-\wt b_{1}-j-1}w^{\wt b_{1}+j-i+i_{t}}dz_{1}=\binom{n_{1}}{i} w^{n_{1}-i},$$
	where $C_{1}$ is a contour of $z_{1}$, with $w$ lying inside. Then by \eqref{n3.32}, \eqref{n3.28}, the recursive formula \eqref{n3.4}, together with the fact that $-\wt b_{1}+n_{1}>-1$, we have:
	\begin{align*}
	&S(v_{3}', b_{1}(n_{1})...b_{m}(n_{m})(v,w)x_{2})\\
	&=\int_{C_{m}}...\int_{C_{1}} \sum_{l=2}^{m}\sum_{i\geq 0}F_{\wt b_{1},i}(z_{1},z_{l})S(v_{3}', (b_{2},z_{2})...(b_{1}(i)b_{l}, z_{l})...(v,w)x_{2}) z_{1}^{n_{1}}...z_{m}^{n_{m}}dz_{1}...dz_{m}\\
	&\ \ \ + \int_{C_{m}}...\int_{C_{1}} \sum_{i\geq 0} F_{\wt b_{1},i}(z_{1},w) S(v_{3}', (b_{2},z_{2})...(b_{m},z_{m})(b_{1}(i)v,w)x_{2}) z_{1}^{n_{1}}...z_{m}^{n_{m}}dz_{1}...dz_{m}\\
	&\ \ \ +\int_{C_{m}}...\int_{C_{1}}\left(\int_{C'_{1}}...\int_{C_{r}'}\right. \sum_{i\geq 0} F_{\wt b_{1},i}(z_{1},w_{t})S(v_{3}',(b_{2},z_{2})...(b_{m},z_{m})(v,w)(c_{1},w_{1})...\\
	&\qquad (b_{1}(i)c_{t},w_{t})...(c_{r},w_{r})v_{2}) \cdot w_{1}^{i_{1}}...w_{r}^{i_{r}}dw_{r}...dw_{1})z_{1}^{n_{1}}...z_{m}^{n_{m}} dz_{1}...dz_{m}\\
	&=\int_{C_{m}}...\int_{C_{2}} \sum_{l=2}^{m}\sum_{i\geq 0}\binom{n_{1}}{i}S(v_{3}', (b_{2},z_{2})...(b_{1}(i)b_{l}, z_{l})...(b_{m},z_{m})(v,w)x_{2}) \\
	&\qquad \cdot z_{2}^{n_{2}}...z_{l}^{n_{1}-i+n_{l}} ...z_{m}^{n_{m}}dz_{2}...dz_{m}\\
	&\ \ \ + \int_{C_{m}}...\int_{C_{2}} \sum_{i\geq 0} \binom{n_{1}}{i}S(v_{3}', (b_{2},z_{2})...(b_{m},z_{m})(b_{1}(i)v,w)x_{2}) w^{n_{1}-i}z_{2}^{n_{2}}...z_{m}^{n_{m}}dz_{2}...dz_{m}\\
	&\ \ \ +\int_{C_{m}}...\int_{C_{2}}\sum_{i\geq 0} \binom{n_{1}}{i}S(v_{3}',(b_{2},z_{2})...(v,w)\left(c_{1}(i_{1})...(b_{1}(i)c_{t})(n_{1}-i+i_{t})...c_{r}(i_{r})v_{2}\right))\\
	&\qquad \cdot z_{2}^{n_{2}}...z_{m}^{n_{m}}dz_{2}...dz_{m}\\
	&=\sum_{l=2}^{m}\sum_{i\geq 0} \binom{n_{1}}{i}S(v_{3}', b_{2}(n_{2})...(b_{1}(i)b_{l})(n_{1}+n_{l}-i)...b_{m}(n_{m})(v,w)x_{2})\\
	&\ \ \ +\sum_{i\geq 0} \binom{n_{1}}{i}S(v_{3}', b_{2}(n_{2})...b_{m}(n_{m})(b_{1}(i)v,w)x_{2})w^{n_{1}-i}\\
	&\ \ \ + S(v_{3}', b_{2}(n_{2})...b_{m}(n_{m})(v,w)(b_{1}(n_{1})x_{2})).
	\end{align*}
	This proves \eqref{n3.33}. Now let $x_{3}' =b_{m}(n_{m})\ds b_{1}(n_{1})v_{3}'\in \bar{M^3}'$, with $\wt b_{i}-n_{i}-1>0$ for all $i$. We use induction on $m$ to show that
	\begin{equation}\label{n3.34}
	\int_{C} S(b_{m}(n_{m})...b_{1}(n_{1})v_{3}', (v,w)x_{2})w^{n}dw=0,
	\end{equation}
	for any fixed $v\in M^{1}$, $x_{2}\in \bar{M^2}$, and $n\in \Z$ such that $n>\deg v+\deg x_{2}-1$. The base case $m=0$ follows from the Corollary \ref{coro3.11}. Let $m>0$, then by \eqref{n3.30} and \eqref{n3.32}, we have:
	\begin{align*}
	&\int_{C} S(b_{m}(n_{m})...b_{1}(n_{1})v_{3}', (v,w)x_{2})w^{n}dw\\
	&=\int_{C} \int_{C_{m}}...\int_{C_{1}}S(v_{3}', (b_{1},z_{1})'...(b_{m},z_{m})'(v,w)x_{2})z_{1}^{-n_{1}-2}...z_{m}^{-n_{m}-2}w^{n}dz_{1}...dz_{m}dw\\
	&=\sum_{j_{1}\geq 0,...,j_{m}\geq 0}\frac{(-1)^{\wt b_{1}+...+\wt b_{m}}}{j_{1}!...j_{m}!} \int_{C}\int_{C_{m}}...\int_{C_{1}} S(v_{3}', (L(1)^{j_{1}}b_{1},z_{1})...(L(1)^{j_{m}}b_{m},z_{m})(v,w)x_{2})\\
	&\qquad \cdot z_{1}^{2\wt b_{1}-n_{1}-2-j_{1}} ...z_{m}^{2\wt b_{m}-n_{m}-2-j_{m}}w^{n}dz_{1}...dz_{m}dw.\\
	&=\sum_{j_{1}\geq 0,...,j_{m}\geq 0} \frac{(-1)^{\wt b_{1}+...+\wt b_{m}}}{j_{1}!...j_{m}!}\int_{C}S(v_{3}', (L(1)^{j_{1}}b_{1})(2\wt b_{1}-n_{1}-2-j_{1})...\numberthis \label{n3.35}\\
	&\qquad ...(L(1)^{j_{m}}b_{m})(2\wt b_{m}-n_{m}-2-j_{m})(v,w)x_{2}).
	\end{align*}
	It suffices to show that each summand in \eqref{n3.35} is $0$. For simplicity, we denote the term $(L(1)^{j_{i}}b_{i})(2\wt b_{i}-n_{i}-2-j_{i})$ by $c_{i}(r_{i})$ for each $i$, note that 
	$$\wt c_{1}(r_{1})=\wt (L(1)^{j_{1}}b_{1})(2\wt b_{1}-n_{1}-2-j_{1})=-\wt b_{1}+n_{1}+1<0.$$
	Then by \eqref{n3.33}, together with the definition formulas \eqref{n3.32} and \eqref{n3.31}, we have:
	\begin{align*}
	&\int_{C}S(v_{3}', c_{1}(r_{1})...c_{m}(r_{m})(v,w)x_{2})w^{n}dw\\
	&=\sum_{l=2}^{m}\sum_{i\geq 0} \binom{r_{1}}{i}\int_{C}S(v_{3}', c_{2}(r_{2})...(c_{1}(i)c_{l})(r_{1}+r_{l}-i)...c_{m}(r_{m})(v,w)x_{2})w^{n}dw\\
	&\ \ \ +\sum_{i\geq 0} \binom{r_{1}}{i}\int_{C}S(v_{3}', c_{2}(r_{2})...c_{m}(r_{m})(c_{1}(i)v,w)x_{2})w^{n+r_{1}-i}dw\\
	&\ \ \ + \int_{C}S(v_{3}', c_{2}(r_{2})...c_{m}(r_{m})(v,w)(c_{1}(r_{1})x_{2}))w^{n}dw\\
	&=\sum_{l=2}^{m}\sum_{i\geq 0}\binom{r_{1}}{i} \int_{C}\underset{(1)}{S(c_{m}'(r_{m})...(c_{1}(i)c_{l})'(r_{1}+r_{l}-i)...c_{2}'(r_{2})v_{3}', (v,w)x_{2})}w^{n}dw\\
	&\ \ \ +\sum_{i\geq 0} \binom{r_{1}}{i}\int_{C} \underset{(2)}{S(c_{m}'(r_{m})...c_{2}'(r_{2})v_{3}',(c_{1}(i)v,w)x_{2})}w^{n+r_{1}-i}dw\\
	&\ \ \ +\int_{C}\underset{(3)}{S(c_{m}'(r_{m})...c_{2}'(r_{2})v_{3}', (v,w)(c_{1}(r_{1})x_{2}))}w^{n}dw\\
	&=(1)+(2)+(3).
	\end{align*}
	Since $\wt c_{1}-r_{1}-1<0$ and $n>\deg v+\deg x_{2}-1$, we have
	\begin{align*}
	\deg (c_{1}(i)v)+\deg x_{2}-1&=\deg v+\deg x_{2}-1+\wt c_{1}-i-1<n+r_{1}-i,\\
	\deg v+\deg (c_{1}(r_{1})x_{2})-1&=\deg v+\deg x_{2}+\wt c_{1}-r_{1}-1-1<n,
	\end{align*}
	for all $i\geq 0$. Then by the induction hypothesis, we have $(1)=(2)=(3)=0$. This finishes the proof of \eqref{n3.34}. Hence $S$ in \eqref{n3.30} belongs to $ \Cor\fusion{M^1}{\bar{M^2}}{\bar{M^3}}$. 
\end{proof} 

So far in this subsection, by abuse of notations, we used the same symbol $S$ \eqref{n3.30} for the extension of a system of correlation functions $S$ in \eqref{n3.11}. We denote the extended $S$ in \eqref{n3.30} by $\psi(S)$ for the rest of this subsection. Then by the Theorem \ref{thm3.12}, we have a linear map:
\begin{equation}\label{n3.36}
\psi: \Cor\fusion{M^{1}}{M^{2}(0)}{M^{3}(0)}\ra \Cor\fusion{M^{1}}{\bar{M^2}}{\bar{M^2}}, \quad S\mapsto \psi(S),
\end{equation}
which is an inverse of the restriction map $\varphi$ in \eqref{n3.9}, with $M^2=\bar{M^2}$ and $M^3=\bar{M^3}$. 

\begin{coro}\label{coro3.13}
	Let $S\in \Cor\fusion{M^1}{M^{2}(0)}{M^{3}(0)}$. Then the linear functional $f$ in Definition \ref{df3.1} is given by $f_v=o(v)=v(\deg v-1)=\Res_{z}I(z,w) w^{\deg v-1+h}$,
	where $I\in I\fusion{M^1}{\bar{M^2}}{\bar{M^2}}$ is the intertwining operator corresponds to $\psi(S)$ in $\Cor\fusion{M^{1}}{\bar{M^2}}{\bar{M^2}}$. 
\end{coro}
\begin{proof}
	By \eqref{n3.3}, we have $S(v_{3}', (v,w)v_{2})=\<v_{3}',f_{v}(v_{2})\>w^{-\deg v}$, for all $v_{3}'\in M^{3}(0)^{\ast}$, $v_{2}\in M^{2}(0)$, and $v\in M^{1}$. On the other hand, by \eqref{n2.16}, 
	$$S(v_{3}', (v,w)v_{2})=\psi(S)(v_{3}', (v,w)v_{2})=(v_{3}',I(v,w)v_{2})=\<v_{3}', v(\deg v-1)v_{2}\> w^{-\deg v},$$
	since $v(m)M^{2}(0)\subseteq M^{3}(\deg v-m-1)$ for any $m\in \Z$. Thus, $f_{v}=v(\deg v-1)$.
\end{proof}

We finish this subsection by showing another property of the space of correlation functions associated with three modules. By \eqref{n3.28} and \eqref{n3.30}, the $\psi (S)$ in \eqref{n3.36} satisfies:
\begin{align*}
&\psi(S)(c_{1}(j_{1})...c_{m}(j_{m})v_{3}',(a_{1},z_{1})...(a_{p},z_{p})(v,w)b_{1}(i_{1})...b_{n}(i_{n})v_{2})\\&=\int_{C'_{1}}...\int_{C'_{m}} \int_{C_{n}}...\int_{C_{1}} S(v_{3}',(c_{m},w_{m})'...(c_{1},w_{1})'(a_{1},z_{1})...(v,w)(b_{1},x_{1})...(b_{n},x_{n})v_{2})\\
&\quad \cdot x_{1}^{i_{1}}...x_{n}^{i_{n}} w_{1}^{-j_{1}-2}...w_{m}^{-j_{m}-2}dx_{1}...dx_{n}dw_{m}...d{w_{1}},\numberthis\label{n3.37}
\end{align*}
where $v_{3}'\in M^{3}(0)^{\ast}$, $v_{2}\in M^{2}(0)$, $v\in M^{1}$, $a_{r},b_{s},c_{t}\in V$ for all $r,s,t$, $C'_{k}$ is a contour of $w_{k}$, $C_{l}$ is a contour of $x_{l}$ for all $k,l$, such that $C_{1}\subset\ds\subset C_{n}\subset C_{1}'\subset \ds\subset C_{m}'$ (we use the subset symbol to indicate one contour is inside of the other), and $z_{1},\ds,z_{n},w$ are outside of $C_{1}'$ but inside of $C_{n}$.

By Proposition \ref{prop3.7} and Theorem 6.2 in \cite{DLM1}, we have an epimorphism of $V$-modules $\pi: \bar{M}(M^{2}(0))\ra \bar{M^2}$, where $\bar{M}(M^{2}(0))$ is the generalized Verma module with bottom level $M^2(0)$. Similarly, there is an epimorphism $\pi: \bar{M}(M^{3}(0)^{\ast})\ra \bar{M^3}'$.
More generally, let $N^{2}$ and $N^{3}$ be any $V$-modules that are generated by their corresponding bottom levels, and assume that $N^{2}(0)=M^{2}(0)$ and $N^{3}(0)=M^{3}(0)$. Suppose there exist epimorphisms $\pi : N^{2}\ra \bar{M^2}$ and $\pi :N^{3'}\ra \bar{M^3}'$. 

If we write $\Res_{z} Y_{N}(b,z)z^{j}=b_{j}$ and $\Res_{z}Y_{\bar{M}}(b,z)z^{j}=b(j)$, then we have
$$\pi (c^{1}_{j_{1}}...c^{m}_{j_{m}}v'_{3})=c^{1}(j_{1})...c^{m}(j_{m})v_{3}',\quad \mathrm{and}\quad \pi( b^{1}_{i_{1}}...b^{n}_{i_{n}}v_{2})=b^{1}(i_{1})...b^{n}(i_{n})v_{2},$$
where $c^{k},b^{l}\in V$, $j_k,i_l\in \Z$ for all $k,l$, $v_{3}'\in M^{3}(0)^{\ast}$, and $v_{2}\in M^{2}(0)$. Thus we have a linear map: $\pi^{\ast} : \Cor\fusion{M^{1}}{\bar{M^2}}{\bar{M^3}}\ra \Cor\fusion{M^1}{N^2}{N^3}$ that is given by: 
\begin{align*}
&\pi^{\ast}(S)(c^{1}_{j_{1}}...c^{m}_{j_{m}}v'_{3},(a_{1},z_{1})...(a_{n},z_{n})(v,w)b^{1}_{i_{1}}...b^{n}_{i_{n}}v_{2})\numberthis \label{n3.38}\\
&= S(c^{1}(j_{1})...c^{m}(j_{m})v_{3}', (a_{1},z_{1})...(a_{n},z_{n})(v,w)b^{1}(i_{1})...b^{n}(i_{n})v_{2}).
\end{align*}
Compose $\psi$ and $\pi^\ast$, we have a linear map 
$\pi^{\ast}\psi: \Cor\fusion{M^{1}}{M^{2}(0)}{M^{3}(0)}\ra \Cor\fusion{M^{1}}{N^{2}}{N^{3}}.$
We claim that $\pi^{\ast} \psi$ is the inverse of the restriction map $\varphi: \Cor\fusion{M^{1}}{N^{2}}{N^{3}}\ra \Cor\fusion{M^{1}}{M^{2}(0)}{M^{3}(0)}$.

Indeed, for $S\in \Cor\fusion{M^{1}}{M^{2}(0)}{M^{3}(0)}$, by \eqref{n3.37} and \eqref{n3.38}, we have:
\begin{align*}
\varphi(\pi^{\ast} \psi)(S) (v_{3}',(a_{1},z_{1})...(a_{n},z_{n})(v,w)v_{2})
&=\psi (S)(\pi(v_{3}'),(a_{1},z_{1})...(a_{n},z_{n})(v,w)\pi(v_{2}))\\
&=S(v_{3}',(a_{1},z_{1})...(a_{n},z_{n})(v,w)v_{2}),
\end{align*}
where $v_{2}\in M^2(0)$ and $v_{3}'\in M^3(0)^\ast$. Hence $\varphi(\pi^{\ast} \psi)=1$. On the other hand, for $S\in \Cor\fusion{M^1}{N^2}{N^3}$, again by \eqref{n3.37} and \eqref{n3.38}, together with the fact that $S$ satisfies \eqref{n2.9} and \eqref{n2.10}, we have for any $c^{1}_{j_{1}}...c^{m}_{j_{m}}v'_{3}\in N^{3'}$,  $b^{1}_{i_{1}}...b^{n}_{i_{n}}v_{2}\in N^2$, $a_{1},...,a_{n}\in V$, and $v\in M^{1}$, 
\begin{align*}
&(\pi^{\ast}\psi)\varphi(S)(c^{1}_{j_{1}}...c^{m}_{j_{m}}v'_{3},(a_{1},z_{1})...(a_{n},z_{n})(v,w)b^{1}_{i_{1}}...b^{n}_{i_{n}}v_{2})\\
&=\psi(\varphi(S))(c^{1}(j_{1})...c^{m}(j_{m})v_{3}', (a_{1},z_{1})...(a_{n},z_{n})(v,w)b^{1}(i_{1})...b^{n}(i_{n})v_{2})\\
&=\int_{C'_{1}}...\int_{C'_{m}} \int_{C_{n}}...\int_{C_{1}} \varphi(S)(v_{3}',(c^{m},w_{m})'...(c^{1},w_{1})'(a_{1},z_{1})...(v,w)(b^{1},x_{1})...(b^{n},x_{n})v_{2})\\
&\quad \cdot x_{1}^{i_{1}}...x_{n}^{i_{n}} w_{1}^{-j_{1}-2}...w_{m}^{-j_{m}-2}dx_{1}...dx_{n}dw_{m}...d{w_{1}},\\
&=\int_{C'_{1}}...\int_{C'_{m}} \int_{C_{n}}...\int_{C_{1}} S(v_{3}',(c^{m},w_{m})'...(c^{1},w_{1})'(a_{1},z_{1})...(v,w)(b^{1},x_{1})...(b^{n},x_{n})v_{2})\\
&\quad \cdot x_{1}^{i_{1}}...x_{n}^{i_{n}} w_{1}^{-j_{1}-2}...w_{m}^{-j_{m}-2}dx_{1}...dx_{n}dw_{m}...d{w_{1}},\\
&=S(c^{1}_{j_{1}}...c^{m}_{j_{m}}v'_{3},(a_{1},z_{1})...(a_{n},z_{n})(v,w)b^{1}_{i_{1}}...b^{n}_{i_{n}}v_{2}).
\end{align*}
This shows $(\pi^{\ast}\psi)\varphi=1$, and so we have $\Cor\fusion {M^{1}}{N^2}{N^{3}}\cong \Cor\fusion{M^1}{M^{2}(0)}{M^{3}(0)}$. In particular, choose $N^{2}=\bar{M}(M^{2}(0))$ and $N^{3}=\bar{M}(M^{3}(0)^\ast)'$, then we have: 
\begin{equation}\label{n3.39}
\Cor\fusion{M^1}{\bar{M}(M^{2}(0))}{\bar{M}(M^{3}(0)^{\ast})'}\cong \Cor\fusion{M^1}{M^2(0)}{M^3(0)}\cong \Cor\fusion{M^1}{\bar{M^2}}{\bar{M^3}}
\end{equation}

Now by \eqref{n3.39}, Corollary \ref{coro2.6} and Theorem \ref{thm3.12}, we have the following theorem:
\begin{thm}\label{thm3.14}
	Let $M^{1}$ be a $V$-module, and let $M^{2}(0)$ and $M^{3}(0)$ be irreducible $A(V)$-modules, then we have the following isomorphism of vector spaces:
	\begin{equation}\label{n3.40}I\fusion{M^1}{\bar{M}(M^{2}(0))}{\bar{M}(M^{3}(0)^{\ast})'}\cong \Cor\fusion{M^1}{M^2(0)}{M^3(0)}\cong I\fusion{M^1}{\bar{M^2}}{\bar{M^3}}.\end{equation}
\end{thm}
If the VOA $V$ is rational, then the generalized Verma module $\bar{M}(U)$ is an irreducible $V$-module for any irreducible $A(V)$-module $U$. Thus, $\bar{M}(M^2(0))=\bar{M^2}=L(M^{2}(0))$, and $\bar{M}(M^{3}(0)^{\ast})'=\bar{M^3}=L(M^{3}(0))$. On the other hand, by Theorem 2.2.2 in \cite{Z}, if $M^2$ and $M^3$ are irreducible $V$-module, then $M^2(0)$ and $M^3(0)$ are irreducible $A(V)$-module. 
\begin{coro}\label{coro3.15}
	Let $V$ be an rational VOA, and let $M^{1}, M^{2}$, and $M^{3}$ be $V$-modules. Suppose $M^{2}$ and $M^{3}$ are irreducible, then we have $\Cor\fusion{M^1}{M^2(0)}{M^3(0)}\cong I\fusion{M^1}{M^2}{M^3}.$
\end{coro}
\begin{remark}\label{remark3.16}
	Let $W^{2}$ and $W^{3}$ be any $\N$-gradable $V$-module that are generated by their corresponding bottom levels, and assume that $W^{2}(0)=M^{2}(0)$ and $W^{3}(0)=M^{3}(0)$. Then there exist epimorphisms:
	$\pi :\bar{M}(M^{2}(0))\ra W^{2},$ and $\pi: \bar{M}(M^{3}(0)^{\ast})\ra W^{3'}.$
	Similar to \eqref{n3.38}, $\pi$ induces a linear map:
	$\pi^{\ast}:\Cor\fusion{M^1}{W^2}{W^3}\hookrightarrow \Cor\fusion  {M^1}{\bar{M}(M^{2}(0))}{\bar{M}(M^{3}(0)^{\ast})'},$
	which is injective since $\pi$ are surjective. Then by Corollary \ref{coro2.6}, \eqref{n3.39}, and \eqref{n3.40}, we have the following estimate for the fusion rule: 
	\begin{equation}\label{n3.41}
	\dim I\fusion{M^1}{W^2}{W^3}\leq \dim \Cor \fusion{M^1}{M^2(0)}{M^3(0)}.
	\end{equation}
\end{remark}

\section{$A(V)$-bimodules and the correlation function $S$}
In this section, we again assume that $M^2(0)$ and $M^3(0)$ are irreducible $A(V)$-modules. By Proposition \ref{prop3.7}, $L(0)=o(\omega)=h_{2}\cdot\Id$ on $M^2(0)$, and $L(0)=h_{3}\cdot \Id$ on $M^3(0)$, for some $h_{2},h_{3}\in \C$. Moreover, $h_{2}$ and $h_{3}$ are the conformal weights of $\bar{M^2}$ and $\bar{M^3}$, respectively.

We will show that $\Cor\fusion{M^{1}}{M^{2}(0)}{M^{3}(0)}$ can be identified with the vector space
$(M^{3}(0)^{\ast}\otimes _{A(V)} B_{h}(M^{1})\otimes _{A(V)}M^{2}(0))^{\ast},$
where $B_{h}(M^{1})$ is an $A(V)$-bimodule that is similar to the $A(V)$-bimodule $A_{0}(M^{1})$ constructed in \cite{HY}. 

However, there are counter-examples showing that this identification is false if one replaces $B_{h}(M^{1})$ by the $A(V)$-bimodule $A(M^{1})$ constructed in Theorem 1.5.1 in \cite{FZ} or  $A_{0}(M^{1})$ constructed in Section 4 of \cite{HY}. The reason is that the correct $L(-1)$-derivation property of the intertwining operators cannot be captured by $A(M^{1})$ nor $A_{0}(M^{1})$. We will see this by the end of this section. 
\subsection{The $A(V)$-bimodule $B_{\la}(W)$}
Let $W$ be a $V$-module with conformal weight $h'$. A sequence of  $A_{N}(V)$-bimodules $A_{N}(W)$ are constructed by Huang and Yang in Section 4 of \cite{HY}. In particular, the $A_{0}(V)=A(V)$-bimodule $A_{0}(W)$ is defined as follows: 

$A_{0}(W)=W/O_{0}(W)$, where $O_{0}(W)=span\{a\circ u,  L(-1)u+(L(0)-h')u:\forall a\in V, u\in W\}$. It is proved (see Theorem 4.7 in \cite{HY}) that $A_{0}(W)$ is an $A(V)$-bimodule under the left and right actions: 
$a\ast_{0} u=\Res_{z} Y_{W}(a,z)u\frac{(1+z)^{\wt a}}{z}$ and $ v\ast_{0} a=\Res_{z} Y_{WV}^{W}(u,z)a\frac{(1+z)^{\deg u}}{z},$
where $Y_{WV}^{W}$ is defined by the skew-symmetry formula (5.1.5) in \cite{FHL}: 
\begin{equation}\label{4.1'}
Y_{WV}^{W}(u,z)a=e^{zL(-1)}Y_{W}(a,-z)u.
\end{equation}

Now let $\la\in \C$ be a fixed complex number, we construct another $A(V)$-bimodule $B_{\la}(W)$ that is similar to $A_{0}(W)$ in the following way:

\begin{df}\label{df4.1'}
	For homogeneous elements $u\in W$ and $a\in V$, define:
	\begin{equation}\label{4.2'}
	u\circ_{\begin{smallmatrix}W\\WV\end{smallmatrix}} a:=\Res_{z}\left(Y_{WV}^{W}(u,z)a\frac{(1+z)^{\deg u+\la}}{z^{2}}\right),
	\end{equation}
	then extend $\circ$ bilinearly to $\circ: W\times V\ra W$. 
	Let $O_{WV}^{W}(W)$ be the vector space spanned by elements \eqref{4.2'} for all $a\in V$ and $u\in W$, and let $B_{\la}(W):=W/(O(W)+O_{WV}^{W}(W)),$
	where $O(W)=span\{a\circ u=\Res_{z}\left(Y_{W}(a,z)u \frac{(1+z)^{\wt a}}{z^{2}} \right): \forall a\in V,u\in W\}$.
\end{df}
\begin{lm}\label{lm 4.2'}
	Let $u\in W$ and $a\in V$ by homogeneous elements, and $m\geq n\geq 0$. Then
	\begin{equation}\label{4.3'}
	\Res_{z} Y_{WV}^{W}(u,z)a\frac{(1+z)^{\deg u+\la+ n}}{z^{2+m}}\in O_{WV}^{W}(W). 
	\end{equation}
\end{lm}
\begin{proof}
	Since $Y_{WV}^{W}(L(-1)u,z)=\frac{d}{dz}Y_{WV}^{W}(u,z)$, the proof of \eqref{4.3'} is almost the same as the proof of Lemma 2.1.2 in \cite{Z}, we omit the details. 
\end{proof}
Recall that the module actions of $A(V)$ on its bimodule $A(W)$ are given by:
$$
b\ast v=\Res_{z}\left(Y_{W}(b,z)v \frac{(1+z)^{\wt b}}{z}\right),\quad \mathrm{and}\quad v\ast b=\Res_{z}\left(Y_{W}(b,z)v\frac{(1+z)^{\wt b-1}}{z}\right),
$$
where $b\in V$ is homogeneous and $v\in W$ (see Definition 1.5.2 in \cite{FZ}). 
\begin{lm}\label{lm 1.2}
	$b\ast O_{WV}^{W}(W)\ssq O_{WV}^{W}(W)$ and $O_{WV}^{W}(W)\ast b\ssq O_{WV}^{W}(W)$, for all $b\in V$.
\end{lm}
\begin{proof}
	Let $u\in W$ and $ b\in V$ be homogeneous, and let $a\in V$. By Definition \ref{df4.1'}, Lemma \ref{lm 4.2'}, and the Jacobi identity of the intertwining operator $Y_{WV}^{W}$, we have:
	\begin{align*}
	&b\ast (u\circ_{\WV}a)\equiv \Res_{z_{1}} Y_{W}(b,z_{1}) \frac{(1+z_{1})^{\wt b}}{z_{1}}\Res_{z_{2}} Y_{WV}^{W}(u,z_{2})a\frac{(1+z_{2})^{\deg u+\la}}{z_{2}^{2}}\\
	&\qquad - \Res_{z_{2}} Y_{WV}^{W}(u,z_{2}) \frac{(1+z_{2})^{\deg u+\la}}{z_{2}^{2}}\Res_{z_{1}} Y_{V}(u,z_{1})a\frac{(1+z_{1})^{\wt b}}{z_{1}} \pmod{O_{WV}^{W}(W)}\\
	&=\Res_{z_{0}}\Res_{z_{2}}Y_{WV}^{W}(Y_{W}(b,z_{0}),z_{2})a\frac{(1+z_{2}+z_{0})^{\wt b}}{z_{2}+z_{0}}\cdot \frac{(1+z_{2})^{\deg u+\la}}{z_{2}^{2}}\\
	&=\Res_{z_{0}}\Res_{z_{2}}\sum_{i=0}^{\wt b}\sum_{j\geq 0}Y_{WV}^{W}(Y_{W}(b,z_{0})u,z_{2})a \binom{\wt b}{i}(-1)^{j} z_{0}^{i+j} \frac{(1+z_{2})^{\deg u+\la+ \wt b-i}}{z_{2}^{2+j+1}}\\
	&=\sum_{i=0}^{\wt b}\sum_{j\geq 0}\binom{\wt b}{i} \Res_{z_2}Y_{WV}^{W}(b_{i+j}u,z_{2})a\frac{(1+z_{2})^{\deg(b_{i+j}u)+\la+(j+1)}}{z_{2}^{2+(j+1)}}\\
	&\equiv 0\pmod{O_{WV}^{W}(W)},
	\end{align*}
	where the last congruence follows from Lemma \ref{lm 4.2'}. By a similar computation, we have:
	\begin{align*}
	(u\circ_{\WV}a)\ast b&\equiv \sum_{i=0}^{\wt b-1}\sum_{j\geq 0}\binom{\wt b-1}{i}\Res_{z_2}Y_{WV}^{W}(b_{i+j}u,z_{2})a\frac{(1+z_{2})^{\deg (b_{i+j}u)+\la+j}}{z_{2}^{2+(j+1)}}\\
	&\equiv 0\pmod{O_{WV}^{W}(W)}. 
	\end{align*}
	Hence we have $b\ast O_{WV}^{W}(W)\ssq O_{WV}^{W}(W)$, and $O_{WV}^{W}(W)\ast b\ssq O_{WV}^{W}(W)$. 
\end{proof}

By Lemma \ref{lm 1.2} and Theorem 1.5.1 in \cite{FZ}, $B_{\la}(W)=W/(O(W)+O_{WV}^{W}(W))$ has an $A(V)$-bimodule structure with respect to $b\ast v$ and $v\ast b$. Moreover, $B_{\la}(W)$ is a quotient module of $A(W)$. In particular, we have the following formula holds in $B_{\la}(W)$: 
\begin{equation}\label{4.4'}
a\ast u-u\ast a\equiv \sum_{j\geq 0} \binom{\wt a-1}{j} a(j)u  \pmod{O_{WV}^{W}(W)+O(W)}, 
\end{equation}
where $a\in V$ homogeneous, and $u\in W$. Let
\begin{equation}\label{4.5'}
O_{\la}(W):=span\{a\circ u, \ L(-1)u+(L(0)-h'+\la)u:\forall a\in V,u\in W\}\subset W.
\end{equation}
\begin{lm}\label{lm4.4'}
	For any $u\in W$, we have $L(-1)u+(L(0)-h'+\la)u\in O^{W}_{WV}(W)$.
\end{lm}
\begin{proof}
	Let $u\in W$ be homogeneous. Since $\deg u=(L(0)-h')u$, we have: 
	\begin{align*}
	u\circ _{\WV} \mathbf{1}&=\Res_{z}e^{zL(-1)}Y_{W}(\mathbf{1},-z)u\frac{(1+z)^{\deg u+\la}}{z^{2}}\\
	&=\Res_{z} \sum_{j\geq 0} \frac{z^{j}}{j!}L(-1)^{j} \sum_{i=0}^{\deg u+\la} \binom{\deg u+\la}{i} z^{i-2}\\
	&=\binom{\deg u+\la}{0} L(-1)u+ \binom{\deg u+\la}{1} L(-1)^{0}u\\
	&=(L(-1)+L(0)-h'+\la)u.
	\end{align*}
	Hence $(L(-1)+(L(0)-h'+\la))u\in O_{WV}^{W}(W)$.
\end{proof}
\begin{lm}\label{lm4.5'}
	We have $O(W)+O_{WV}^{W}(W)=O_{\la}(W)$. In particular, $B_{\la}(W)=W/O_{\la}(W)$. 
\end{lm}
\begin{proof}
	By Lemma \ref{lm4.4'}, we only need to show that $O_{WV}^{W}(W)\ssq O_{\la}(W)$. Similar to the proof of Lemma 2.1.3 in \cite{Z}, for any homogeneous $u\in W$ and $a\in V$, we have:
	$Y_{WV}^{W}(u,z)a\equiv (1+z)^{-\deg u-\la-\wt a} Y_{W}\left(a,\frac{-z}{1+z}\right)u \pmod{O_{\la}(W)}.$
	It follows that 
	\begin{align*}
	u\circ_{\WV} a &=\Res_{z} Y_{WV}^{W}(u,z)a \frac{(1+z)^{\deg u+\la}}{z^{2}} \\
	&\equiv \Res_{z} Y_{W}\left(a,\frac{-z}{1+z}\right)u \frac{(1+z)^{-\wt a}}{z^{2}}\pmod{O_{\la}(W)}\\
	&\equiv -\Res_{w} Y_{W}(a,w)u \frac{(1+w)^{\wt a}}{w^{2}} \pmod{O_{\la}(W)}.
	\end{align*}
	Hence $u\circ_{\WV}a\equiv -a\circ u\pmod{O_{\la}(W)}$, and so $O^{W}_{WV}(W)+O(W)=O_{\la}(W)$. 
\end{proof}

Now let $W=M^{1}$, and $\la=h=h_{1}+h_{2}-h_{3}$. Then by \eqref{4.5'} and Lemma \ref{lm4.5'}, $B_{h}(M^{1})=M^{1}/ O_{h}(M^{1})$, where $O_{h}(M^{1})=span\{ a\circ u,\ L(-1)u+(L(0)+h_{2}-h_{3})u: \forall a\in V, u\in M^{1}\}.$

\begin{lm}\label{lm4.6'}
	Let $I\in I\fusion{M^{1}}{\bar{M^2}}{\bar{M^3}}$, then the linear map $$o: M^{1}\ra \Hom_{
		\C}(M^{2}(0),M^{3}(0)),\ o(v)=v(\deg v-1)=\Res _{z} I(v,z) z^{\deg v-1+h}$$ factors through $B_{h}(M^{1})=M^{1}/O_{h}(M^{1})$. 
\end{lm}
\begin{proof}
	By Lemma \ref{lm4.5'}, we need to show that $o(O_{h}(M^{1}))=0$. By Lemma 1.5.2 in \cite{FZ}, we already have $o(a\circ u)=0$ for all $a\in V$ and $u\in M^{1}$. Furthermore,  by the $L(-1)$-derivation property of $I$, we have: 
	\begin{align*}
	o(L(-1)v)&=\Res _{z} I(L(-1)v,z)z^{\deg v+1-1+h}=\Res_{z} \left(\frac{d}{dz} I(v,z)\right)  z^{\deg v+h}\\
	&=\Res_{z} I(v,z) (-\deg v-h) z^{\deg v+h-1}=-((L(0)-h_{1}+h)v)(\deg v-1) \\
	&=-o((L(0)+h_{2}-h_{3})v).
	\end{align*}
	Hence $o(O_{h}(M^{1}))=0$, and so $o$ factors through $B_{h}(M^{1})$. 
\end{proof}
\begin{prop}\label{prop4.7'}
	There exists an injective linear map:
	\begin{align}\label{4.6'}
	\nu: \Cor\fusion{M^1}{M^{2}(0)}{M^{3}(0)}&\ra (M^{3}(0)^{\ast}\otimes _{A(V)} B_{h}(M^{1})\otimes _{A(V)}M^{2}(0))^{\ast}\\
	S&\mapsto f_{S},\qquad f_{S}(v_{3}'\otimes v\otimes v_{2}):=\<v_{3}', f_{v}(v_{2})\>,\nonumber
	\end{align}
	where we use the same symbol $v$ for its image in $B_{h}(M^{1})$.
\end{prop}
\begin{proof} First, we have $f_{v}=o(v)$ by Corollary \ref{coro3.13}, where $o(v)=\Res_{w}I(v,w)w^{\deg v-1+h}$, and $I\in I\fusion{M^1}{\bar{M^2}}{\bar{M^3}}$ is the intertwining operator corresponds to $\psi(S)$ in $\Cor\fusion{M^1}{\bar{M^2}}{\bar{M^3}}$, see Theorem \ref{thm3.14}. Moreover, it follows from Lemma \ref{lm4.6'} that $o(O_{h}(M^{1}))=0$. Hence $\nu$ is well-defined.
	The injectivity of $\nu$ follows from Lemma \ref{lm3.3}.
\end{proof}

\begin{remark}Although our definition for $B_{h}(M^{1})$ is similar to the $A(V)$-bimodule $A_{0}(M^{1})$ constructed by Huang and Yang in \cite{HY}, they are not isomorphic as $A(V)$-bimodules. We will give a counter-example in the next subsection. 
\end{remark}

Our goal next is to construct an inverse map of $\nu$ in \eqref{4.6'}. Given a $f\in (M^{3}(0)^{\ast}\otimes _{A(V)} B_{h}(M^{1})\otimes _{A(V)}M^{2}(0))^{\ast}$, we need to construct a  corresponding system of correlation functions $S$ in $ \Cor\fusion{M^1}{M^{2}(0)}{M^{3}(0)}$. Our strategy is to use the recursive formulas \eqref{n3.4} and \eqref{n3.6} and construct the system of functions $S$ inductively. The key is to show the locality ((2) in Definition \ref{df2.1}) in each step, which can be achieved by the properties of the $A(V)$-bimodule $B_{h}(M^1)$, together with the formula \eqref{n3.8}. 

\subsection{The construction of $4$-point and $5$-point functions}
From now on, we fix a linear function $f$ on the vector space $M^{3}(0)^{*}\otimes _{A(V)}B_{h}(M^{1})\otimes _{A(V)}M^{2}(0).$

\begin{df}\label{df4.9}
	Define $S_{M}: M^{3}(0)^{*}\times M^{1}\times M^{2}(0)\rightarrow \mathcal{F}(w)$ by 
	\begin{equation}\label{n4.7}
	S_{M}(v_{3}',(v,w)v_{2}):=f(v_{3}'\otimes v\otimes v_{2})w^{-\deg v},
	\end{equation}
	where on the right-hand side we use the same symbol $v$ for its image $v+O(M^{1})$ in $B_h(M^1)$.
	
	Define $S_{VM}^{L}: M^{3}(0)^{*}\times V\times M^{1}\times M^{2}(0)\rightarrow \mathcal{F}(z,w)$ by 
	\begin{equation}\label{n4.8}
	\begin{aligned}
	S_{VM}^{L}(v_{3}',(a,z)(v,w)v_{2}):&=S_{M}(v_{3}'o(a),(v,w)v_{2})z^{-\wt a}\\
	&\ \ \ +\sum_{i\geq 0}F_{\wt a, i} (z,w)S_{M}(v_{3}',(a(i)v,w)v_{2}).
	\end{aligned}
	\end{equation}
	Finally, define $S_{MV}^{R}: M^{3}(0)^{*}\times M^{1}\times V\times M^{2}(0)\rightarrow \mathcal{F}(z,w)$ by 
	\begin{equation}\label{n4.9}
	\begin{aligned}
	S_{MV}^{R}(v_{3}',(v,w)(a,z)v_{2}):&=S_{M}(v_{3}',(v,w)o(a)v_{2})z^{-\wt a}\\
	&\ \ \ +\sum_{i\geq 0}G_{\wt a, i} (z,w)S_{M}(v_{3}',(a(i)v,w)v_{2}).
	\end{aligned}
	\end{equation}
\end{df} 
The upper index $L$ (resp.$R$) in the 4-point functions $S$ indicates that we use the expansion formula for the left (resp. right) most term, namely, \eqref{n3.4} (resp.\eqref{n3.6}) to construct the new $S$ from the $3$-point function. We will denote the 3-point function $S_{M}$ by $S$.
\begin{prop}\label{prop4.10}
	As rational functions in $\mathcal{F}(z,w)$, we have: $$S_{VM}^{L}(v_{3}',(a,z)(v,w)v_{2})=S_{MV}^{R}(v_{3}',(v,w)(a,z)v_{2}).$$
\end{prop}
\begin{proof}
	By Definition \ref{df4.9}, \eqref{n3.8}, and the property of $M^{3}(0)^{*}\otimes _{A(V)}B_{h}(M^{1})\otimes _{A(V)}M^{2}(0)$, 
	\begin{align*}
	&S_{VM}^{L}(v_{3}',(a,z)(v,w)v_{2})-S_{MV}^{R}(v_{3}',(v,w)(a,z)v_{2})\\
	&=f(v_{3}'o(a)\otimes v\otimes v_{2})w^{-\deg v}z^{-\wt a}-f(v_{3}'\otimes v\otimes o(a)v_{2})w^{-\deg v}z^{-\wt a}\\
	&\ \ \ +\sum_{i\geq 0} (F_{\wt a,i}(z,w)-G_{\wt a,i}(z,w))S_M(v_{3}',(a(i)v,w)v_{2})\\
	&=f(v_{3}'\otimes a\ast v\otimes v_{2})w^{-\deg v}z^{-\wt a}-f(v_{3}'\otimes v\ast a\otimes v_{2})w^{-\deg v}z^{-\wt a}\\
	&\ \ \ -\sum_{i\geq 0}\binom{\wt a-1}{i}f(v_{3}'\otimes a(i)v\otimes v_{2})w^{-\deg v-\wt a+i+1}z^{-\wt a}w^{\wt a-1-i}\\
	&=f(v_{3}'\otimes(a\ast v-v\ast a)\otimes v_{2})w^{-\deg v}z^{-\wt a}-\sum_{i\geq 0}\binom{\wt a-1}{i}f(v_{3}'\otimes a(i)v\otimes v_{2})z^{-\wt a}w^{-\deg v}.
	\end{align*}
	By \eqref{4.4'}, we also have $a\ast v-v\ast a=\sum_{i\geq 0}\binom{\wt a-1}{i}a(i)v$ holds in the $A(V)$-bimodule $B_{h}(M^{1})$. Hence  $S_{VM}^{L}(v_{3}',(a,z)(v,w)v_{2})-S_{MV}^{R}(v_{3}',(v,w)(a,z)v_{2})=0$.
\end{proof}

By Proposition \ref{prop4.10}, the 4-point functions $S_{VM}^{L}$ and $S_{MV}^{R}$ in definition \ref{df4.9} give rise to one single 4-point function $S$ that satisfies
\begin{equation}\label{n4.10}
S(v_{3}',(a,z)(v,w)v_{2})=S(v_{3}',(v,w)(a,z)v_{2}),
\end{equation}
and this function can be defined either by \eqref{n4.8} or \eqref{n4.9}. 

We adopt a similar method to construct $5$-point functions. As long as the term $(v,w)$ does not appear at the left-most place, we use the formula \eqref{n3.4} to construct $S$ from the $4$-point function; if $(v,w)$ appears at the left-most place, we use \eqref{n3.6} to construct $S$.

\begin{df}\label{df4.11} Define the $5$-point functions with the upper index $L$,
	$$S_{VMV}^{L}(v_{3}',(a_{1},z_{1})(v,w)(a_{2},z_{2})v_{2}),\quad \mathrm{and} \quad S_{VVM}^{L}(v_{3}',(a_{1},z_{1})(a_{2},z_{2})(v,w)v_{2}),$$ 
	by expanding $(a_{1},z_{1})$ from the left, which is given by the common formula:
	\begin{equation}\label{n4.11}
	\begin{aligned}
	S(v_{3}'o(a_1),(v,w)(a_{2},z_{2})v_{2})z_{1}^{-\wt a_{1}}&+\sum_{j\geq 0}F_{\wt a_{1},j}(z_{1},w)S(v_{3}',(a_{1}(j)v,w)(a_{2},z_{2})v_{2})\\
	&+\sum_{j\geq 0}F_{\wt a_{1},j}(z_{1},z_{2})S(v_{3}',(v,w)(a_{1}(j)a_{2},z_{2})v_{2}).
	\end{aligned}
	\end{equation}
	Define the $5$-point functions with upper index $R$,
	$$S_{VMV}^{R}(v_{3}',(a_{2},z_{2})(v,w)(a_{1},z_{1})v_{2}), \quad \mathrm{and}\quad S_{MVV}^{R}(v_{3}',(v,w)(a_{2},z_{2})(a_{1},z_{1})v_{2}),$$
	by expanding $(a_{1},z_{1})$ from the right, which is given by the common formula:
	\begin{equation}\label{n4.12}
	\begin{aligned}
	S(v_{3}',(a_{2},z_{2})(v,w)o(a_{1})v_{2})z_{1}^{-\wt a_{1}}&+\sum_{j\geq 0}G_{\wt a_{1},j}(z_{1},w)S(v_{3}',(a_{2},z_{2})(a_{1}(j)v,w)v_{2})\\
	&+\sum_{j\geq 0}G_{\wt a_{1},j}(z_{1},z_{2})S(v_{3}',(a_{1}(j)a_{2},z_{2})(v,w)v_{2}).
	\end{aligned}
	\end{equation}
\end{df}
The function $S$ in \eqref{n4.11} and \eqref{n4.12} is the (common) $4$-point function in Definition \ref{df4.9}. By \eqref{n4.10}, it makes sense to define $S^L_{VMV}$ and $S^L_{VVM}$ by the same formula, same for $S^R_{VMV}$ and $S^R_{MVV}$. 
We will show that all the $5$-point functions in Definition \ref{df4.11} are the same. First, we
observe that the term $S_{VMV}(v_{3}',(a_{1},z_{1})(v,w)(a_{2},z_{2})v_{2})$ has the following two expressions:
$ S_{VMV}^{L}(v_{3}',(a_{1},z_{1})(v,w)(a_{2},z_{2})v_{2})$ and $S_{VMV}^{R}(v_{3}',(a_{1},z_{1})(v,w)(a_{2},z_{2})v_{2}). $

\begin{prop}\label{prop4.12}
	If \eqref{n4.11}=\eqref{n4.12}, then we have: 
	$$S_{VMV}^{L}(v_{3}',(a_{1},z_{1})(v,w)(a_{2},z_{2})v_{2})=S_{VMV}^{R}(v_{3}',(a_{1},z_{1})(v,w)(a_{2},z_{2})v_{2}).$$
\end{prop} 
\begin{proof}
	Note that \eqref{n4.11} is a generalization of the function (2.2.6) in \cite{Z}. By a similar calculation, it is easy to see that the formula (2.2.11) in \cite{Z} also holds for our case. i.e., we can swap the terms $(a_{1},z_{1})$ and $(a_{2},z_{2})$ in $S^{L}_{VVM}$: 
	\begin{equation}\label{n4.13}
	S^{L}_{VVM}(v_{3}',(a_{1},z_{1})(a_{2},z_{2})(v,w)v_{2})=S^{L}_{VVM}(v_{3}',(a_{2},z_{2})(a_{1},z_{1})(v,w)v_{2}).
	\end{equation}
	By the assumption that \eqref{n4.11}=\eqref{n4.12}, Definition \ref{df4.11}, and \eqref{n4.13}, we have: 
	\begin{align*}
	&S_{VMV}^{L}(v_{3}',(a_{1},z_{1})(v,w)(a_{2},z_{2})v_{2})=S_{VVM}^{L}(v_{3}',(a_{1},z_{1})(a_{2},z_{2})(v,w)v_{2})\\
	&=S^{L}_{VVM}(v_{3}',(a_{2},z_{2})(a_{1},z_{1})(v,w)v_{2})=S^{L}_{VMV}(v_{3}',(a_{2},z_{2})(v,w)(a_{1},z_{1})v_{2})\\
	&=S^{R}_{VMV}(v_{3}',(a_{1},z_{1})(v,w)(a_{2},z_{2})v_{2}),
	\end{align*}
	where the last equality follows from the assumption that \eqref{n4.11}=\eqref{n4.12}. 
\end{proof}
Next, we show that \eqref{n4.11}=\eqref{n4.12}. We use symbols $(1),(2)$, and $(3)$ to denote the difference of the three summands in  the term $\eqref{n4.11}-\eqref{n4.12}$: 

\begin{equation}\tag*{(1)}\label{(1)}
S(v_{3}'o(a_{1}),(v,w)(a_{2},z_{2})v_{2})z_{1}^{-\wt a_{1}}-S(v_{3}',(a_{2},z_{2})(v,w)o(a_{1})v_{2})z_{1}^{-\wt a_{1}}.
\end{equation}
\begin{equation}\tag*{(2)}\label{(2)}
\sum_{j\geq 0}(F_{\wt a_{1},j}(z_{1},w)-G_{\wt a_{1},j}(z_{1},w))S(v_{3}',(a_{1}(j)v,w)(a_{2},z_{2})v_{2}).
\end{equation}
\begin{equation}\tag*{(3)}\label{(3)}
\sum_{j\geq 0}(F_{\wt a_{1},j}(z_{1},z_{2})-G_{\wt a_{1},j}(z_{1},z_{2}))S(v_{3}',(v,w)(a_{1}(j)a_{2},z_{2})v_{2}).
\end{equation}
So we need to show that \ref{(1)}+\ref{(2)}+\ref{(3)}=0.

By \eqref{n4.10}, we may use the formula \eqref{n4.8} and expand both terms in \ref{(1)} with respect to $(a_{2},z_{2})$ from the left. Then \ref{(1)} can be expressed as:

\begin{align*}
&S(v_{3}'o(a_{1}),(v,w)(a_{2},z_{2})v_{2})z_{1}^{-\wt a_{1}}-S(v_{3}',(a_{2},z_{2})(v,w)o(a_{1})v_{2})z_{1}^{-
	\wt a_{1}}\\
&=S(v_{3}'o(a_{1})o(a_{2}),(v,w)v_{2})z_{1}^{-\wt a_{1}}z_{2}^{-\wt a_{2}}+\sum_{i\geq 0}F_{\wt a_{2},i}(z_{2},w)S(v_{3}'o(a_{1}),(a_{2}(i)v,w)v_{2})z_{1}^{-\wt a_{1}}\\
&\ \ -S(v_{3}'o(a_{2}),(v,w)o(a_{1})v_{2})z_{1}^{-\wt a_{1}}z_{2}^{-\wt a_{2}}+\sum_{i\geq 0}F_{\wt a_{2},i}(z_{2},w)S(v_{3}',(a_{2}(i)v,w)o(a_{1})v_{2})z_{1}^{-\wt a_{1}}\\
&=\underset{(11)}{ f(v_{3}'\otimes a_{1}\ast a_{2}\ast v\otimes v_{2})w^{-\deg v}z_{1}^{-\wt a_{1}}z_{2}^{-\wt a_{2}}}-\underset{(12)}{ f(v_{3}'\otimes a_{2}\ast v\ast a_{1}\otimes v_{2})w^{-\deg v}z_{1}^{-\wt a_{1}}z_{2}^{-\wt a_{2}}}\\
&\ \ +\sum_{i\geq 0}\underset{(13)}{F_{\wt a_{2},i}(z_{2},w)f(v_{3}'\otimes (a_{1}\ast (a_{2}(i)v))-(a_{2}(i)v)\ast a_{1})\otimes v_{2}) w^{-\wt a_{2}-\deg v+i+1}z_{1}^{-
		\wt a_{1}}}\\
&=(11)+(12)+(13).
\end{align*}
For the term \ref{(2)}, we use the formula \eqref{n4.8} agian and expand each summand in \ref{(2)} with respect to $(a_{2},z_{2})$ from the left. Then by \eqref{n3.8}, \ref{(2)} can be expressed as:
\begin{align*}
\eqref{(2)}&=\sum_{j\geq 0}F_{\wt a_{1},j}(z_{1},w)S(v_{3}'o(a_{2}),(a_{1}(j)v,w)v_{2})z_{2}^{-\wt a_{2}}\\
&\ \ +\sum_{j\geq 0}\sum_{i\geq 0} F_{\wt a_{1},j}(z_{1},w)F_{\wt a_{2},i}(z_{2},w)S(v_{3}',(a_{2}(i)a_{1}(j)v,w)v_{2})\\
&\ \ -\sum_{j\geq 0}G_{\wt a_{1},j}(z_{1},w)S(v_{3}'o(a_{2}),(a_{1}(j)v,w)v_{2})z_{2}^{-\wt a_{2}}\\
&\ \ -\sum_{j\geq 0}\sum_{\i\geq 0}G_{\wt a_{1},j}(z_{1},w)F_{\wt a_{2},i}(z_{2},w)S(v_{3}',(a_{2}(i)a_{1}(j)v,w)v_{2})\\
&=\sum_{j\geq 0}-\binom{\wt a_{1}-1}{j}\underset{(21)}{S(v_{3}'o(a_{2}),(a_{1}(j)v,w)v_{2})z_{1}^{-
		\wt a_{1}}z_{2}^{-\wt a_{2}}w^{\wt a_{1}-1-j}}\\
&\ \ +\sum_{j\geq 0}\sum_{i\geq 0}-\binom{\wt a_{1}-1}{j}\underset{(22)}{z_{1}^{-\wt a_{1}}w^{
		\wt a_{1}-1-j}F_{\wt a_{2},w}(z_{2},w)S(v_{3}',(a_{2}(i)a_{1}(j)v,w)v_{2})}\\
&=(21)+(22).
\end{align*}

Finally, for the term \ref{(3)}, we expand each of its summand with respect to $(a_{1}(j)a_{2},z_{2})$ from the left, so \ref{(3)} can be expressed as:

\begin{align*}
&\eqref{(3)}=\sum_{j\geq 0} F_{\wt a_{1},j}(z_{1},z_{2})S(v_{3}'o(a_{1}(j)a_{2}),(v,w)v_{2})z_{2}^{-\wt a_{1}-\wt a_{2}+j+1}\\
&\ \ +\sum_{j\geq 0}\sum_{i\geq 0}F_{\wt a_{1},j}(z_{1},z_{2})F_{\wt a_{1}+\wt a_{2}-j-1,i}(z_{2},w)S(v_{3}',((a_{1}(j)a_{2})(i)v,w)v_{2})\\
&\ \ -\sum_{j\geq 0} G_{\wt a_{1},j}(z_{1},z_{2})S(v_{3}'o(a_{1}(j)a_{2}),(v,w)v_{2})z_{2}^{-\wt a_{1}-\wt a_{2}+j+1}\\
&\ \ +\sum_{j\geq 0}\sum_{i\geq 0} G_{\wt a_{1},j}(z_{1},z_{2})F_{\wt a_{1}+\wt a_{2}-j-1,i}(z_{2},w)S(v_{3}', ((a_{1}(j)a_{2})(i)v,w)v_{2})\\
&=\sum_{j\geq 0}-\binom{\wt a_{1}-1}{j}\underset{(31)}{z_{1}^{-\wt a_{1}}z_{2}^{\wt a_{1}-1-j}S(v_{3}'o(a_{1}(j)a_{2}),(v,w)v_{2})z_{2}^{-\wt a_{1}-\wt a_{2}+j+1}}\\
&\ \ +\sum_{j\geq 0}\sum_{i\geq 0}-\binom{\wt a_{1}-1}{j}\underset{(32)}{z_{1}^{-\wt a_{1}}z_{2}^{\wt a_{1}-1-j}F_{\wt a_{1}+\wt a_{2}-j-1,i}(z_{2},w)S(v_{3}',(a_{1}(j)a_{2})(i)(v,w)v_{2})}\\
&=(31)+(32).
\end{align*}
We need to show that  $(11)+(12)+(13)+(21)+(22)+(31)+(32)=0$.
In fact, since 
$a\ast v-v\ast a=\textrm{Res}_{z}Y(a,z)v(1+z)^{\wt a-1}=\sum_{j \geq 0}\binom{\wt a-1}{j}a(j)v$ in $B_{h}(M^{1})$, see \eqref{4.4'}, and 
$a_{1}\ast a_{2}-a_{2}\ast a_{1}=\sum_{j\geq 0}\binom{\wt a_{1}-1}{j} a_{1}(j)a_{2}$ in $A(V)$, we can rewrite $(21)$ and $(31)$ as:
\begin{align*}
(21)&=-\sum_{j\geq 0}\binom{\wt a_{1}-1}{j}w^{-\wt a_{1}-\deg v+j+1}z_{1}^{\wt a_{1}}z_{2}^{\wt a_{2}}w^{\wt a_{1}-j-1}f(v_{3}'o(a_{2})\otimes a_{1}(j)v\otimes v_{2})\\
&=-w^{-\deg v}z_{1}^{-\wt a_{1}}z_{2}^{\wt a_{2}}f(v_{3}'\otimes(a_{2}\ast a_{1}\ast v-a_{2}\ast v\ast a_{1})\otimes v_{2});\\
(31)&=-\sum_{j\geq 0}\binom{\wt a_{1}-1}{j}z_{1}^{-\wt a_{1}}z_{2}^{-\wt a_{2}}w^{-\deg v}f(v_{3}'o(a_{1}(j)a_{2})\otimes v\otimes v_{2})\\
&=-z_{1}^{-\wt a_{1}}z_{2}^{-\wt a_{2}}w^{-\deg v}f(v_{3}'\otimes(a_{1}\ast a_{2}\ast v-a_{2}\ast a_{1}\ast v)\otimes v_{2}).
\end{align*}
Then by the bimodule property of $B_{h}(M^{1})$, we have:
\begin{align*}
&\ \ (11)+(12)+(21)+(31)\\
&=f(v_{3}'\otimes a_{1}\ast a_{2}\ast v\otimes v_{2})w^{-\deg v}z_{1}^{-\wt a_{1}}z_{2}^{-\wt a_{2}}- f(v_{3}'\otimes a_{2}\ast v\ast a_{1}\otimes v_{2})w^{-\deg v}z_{1}^{-\wt a_{1}}z_{2}^{-\wt a_{2}}\\
&\ \ -w^{-\deg v}z_{1}^{-\wt a_{1}}z_{2}^{\wt a_{2}}f(v_{3}'\otimes(a_{2}\ast a_{1}\ast v-a_{2}\ast v\ast a_{1})\otimes v_{2})\\
&\ \ -z_{1}^{-\wt a_{1}}z_{2}^{-\wt a_{2}}w^{-\deg v}f(v_{3}'\otimes(a_{1}\ast a_{2}\ast v-a_{2}\ast a_{1}\ast v)\otimes v_{2})\\
&=0.
\end{align*}
It remains to show that $(13)+(22)+(32)=0$.
\begin{lm}\label{lm4.5}
	Let $M$ be a $V$ module, and let $a_{1},a_{2}\in V$, $v\in M$, and $n\in \N$. We have:
	\begin{equation}\label{n4.14}
	\begin{aligned}
	&\sum_{i,j\geq 0}\binom{\wt a_{1}-1}{j}\binom{\wt a_{2}+n}{i}(a_{1}(j)a_{2}(i)v-a_{2}(i)a_{1}(j)v)\\
	&=\sum_{i,j\geq 0} \binom{\wt a_{1}-1}{j}\binom{\wt a_{1}+\wt a_{2}-j-1+n}{i}(a_{1}(j)a_{2})(i)v
	\end{aligned}
	\end{equation}
\end{lm}
\begin{proof}
	Choose complex variables $z_{1},z_{2}$ in the domain $|z_{1}|<1,\ |z_{2}|<1,\ |z_{1}-z_{2}|<|1+z_{2}|$.\\
	By the Jacobi identity in the residue form, the left-hand side of \eqref{n4.14} can be written as:
	\begin{align*}
	&\textrm{Res}_{z_{1},z_{2}}\sum_{i,j\geq 0}\binom{\wt a_{1}-1}{j}\binom{\wt a_{2}+n}{i}z_{1}^{j}z_{2}^{i}(Y(a_{1},z_{1})Y(a_{2},z_{2})v-Y(a_{2},z_{2})Y(a_{1},z_{1})v)\\
	&=\textrm{Res}_{z_{1},z_{2}}(1+z_{1})^{\wt a_{1}-1}(1+z_{2})^{\wt a_{2}+n}(Y(a_{1},z_{1})Y(a_{2},z_{2})v-Y(a_{2},z_{2})Y(a_{1},z_{1})v)\\
	&=\textrm{Res}_{z_{2}}\textrm{Res}_{z_{1}-z_{2}}(1+z_{2}+(z_{1}-z_{2}))^{\wt a_{1}-1}(1+z_{2})^{\wt a_{2}+n}Y(Y(a_{1},z_{1}-z_{2})a_{2},z_{2})v\\
	&=\textrm{Res}_{z_{2}}\textrm{Res}_{z_{1}-z_{2}}\sum_{j\geq 0}\binom{\wt a_{1}-1}{j}(1+z_{2})^{\wt a_{1}-1-j+\wt a_{2}+n}(z_{1}-z_{2})^{j}Y(Y(a_{1},z_{1}-z_{2})a_{2},z_{2})v\\
	&=\sum_{i,j\geq 0}\binom{\wt a_{1}-1}{j}\binom{\wt a_{1}+\wt a_{2}-j-1+n}{i}(a_{1}(j)a_{2})(i)v,
	\end{align*}
	which is the right-hand side of \eqref{n4.14}.
\end{proof}
We use the formula \eqref{4.4'} again and rewrite (13) as:
\begin{align*}
(13)=\sum_{i,j\geq 0}\binom{\wt a_{1}-1}{j}z_{1}^{-\wt a_{1}}w^{-\wt a_{2}-\deg v+i+1}F_{\wt a_{2},i}(z_{2},w)f(v_{3}'\otimes a_{1}(j)a_{2}(i)v\otimes v_{2}).
\end{align*}
Since the map $\iota_{z_{2},w}$ is injective (see Section 3 in \cite{FHL}), we only need to show that  $\iota_{z_{2},w}((13)+(22)+(32))=0$. By \eqref{n3.5}, $\iota_{z_{2},w}(F_{\wt a_{2},i}(z_{2},w))$ can be written as:
\begin{align*}
\iota_{z_{2},w}(F_{\wt a_{2},i}(z_{2},w))=\sum_{n\geq 0}\binom{\wt a_{2}+n}{i}w^{\wt a_{2}+n-i}z_{2}^{-\wt a_{2}-n-1}
\end{align*}
To simplify our notation, we denote $z_{1}^{\wt a_{1}}w^{-\deg v+n+1}z_{2}^{-\wt a_{2}-n-1}$ by $\gamma$. By Lemma \ref{lm4.5},
\begin{align*}
&\ \ \iota_{z_{2},w}(13)+\iota_{z_{2},w}(22)\\
&=\sum_{i,j\geq 0}\binom{\wt a_{1}-1}{j}z_{1}^{\wt a_{1}}w^{-\wt a_{2}-\deg v+i+1}\left(\sum_{n\geq 0}\binom{\wt a_{2}+n}{i}w^{\wt a_{2}+n-i}z_{2}^{-\wt a_{2}-n-1}\right)\\
&\ \ \cdot (f(v_{3}'\otimes a_{1}(j)a_{2}(i)v\otimes v_{2})-f(v_{3}'\otimes a_{2}(i)a_{1}(j)v\otimes v_{2}))\\
&=\sum_{i,j,n\geq 0}\binom{\wt a_{1}-1}{j}\binom{\wt a_{2}+n}{i}\gamma \cdot f(v_{3}'\otimes (a_{1}(j)a_{2}(i)v-a_{2}(i)a_{1}(j)v)\otimes v_{2})\\
&=\sum_{i,j,n\geq 0}\binom{\wt a_{1}-1}{j}\binom{\wt a_{1}+\wt a_{2}+n-j-1}{i}\gamma \cdot f(v_{3}'\otimes (a_{1}(j)a_{2})(i)v\otimes v_{2})\\
&=-\iota_{z_{2},w}(32).
\end{align*}
Now the proof of \eqref{n4.11}=\eqref{n4.12} is complete. 

Therefore, the 5 point functions in Definition \ref{df4.11} give rise to one single $5$-point function $S$ that satisfies:
\begin{align*}
S(v_{3}',(a_{1},z_{1})(a_{2},z_{2})(v,w)v_{2})&=S(v_{3}',(a_{2},z_{2})(a_{1},z_{1})(v,w)v_{2})\\
=S(v_{3}',(a_{1},z_{1})(v,w)(a_{2},z_{2})v_{2})&=S(v_{3}',(a_{2},z_{2})(v,w)(a_{1},z_{1})v_{2})\numberthis \label{n4.15}\\
=S(v_{3}',(v,w)(a_{1},z_{1})(a_{2},z_{2})v_{2})&=S(v_{3}',(v,w)(a_{2},z_{2})(a_{1},z_{1})v_{2}).
\end{align*}
In particular, the $5$-point function $S$ satisfies the locality in Definition \ref{df2.1}, with $v_{3}'\in M^3(0)^\ast$ and $v_2\in M^2(0)$. Moreover, $S(v_{3}',(a_{1},z_{1})(a_{2},z_{2})(v,w)v_{2})$ also satisfies both of the recursive formula \eqref{n3.4} and \eqref{n3.6} by its definition.

\subsection{Construction of $(n+3)$-point functions}
We construct the general $(n+3)$-point function $S$ using induction on $n$. We have finished the base cases $n=1,2$ in the previous subsection.
Now assume the $(n+2)$-point function: 
$$S:M^{3}(0)^{*}\times V\times\ds \times M^{1}\times\ds\times V\times M^{2}(0) \rightarrow \mathcal{F}(z_{1},\ds,z_{n-1},w)$$
exist and satisfy the following two properties: Let $\{(b_1,w_1),(b_2,w_2),\ds,(b_n,w_n)\}$ be the same set as $\{(a_1,z_1),\ds,(a_{n-1},z_{n-1}), (v,w)\}$. 
The first property is the locality:
\begin{equation}\tag{{I}}\label{I}
S(v_{3}',(a_{1},z_{1})(a_{2},z_{2})...(a_{n-1},z_{n-1})(v,w)v_{2})=S(v_{3}',(b_{1},w_{1})(b_{2},w_{2})...(b_{n},w_{n})v_{2}),
\end{equation}
that is, the terms $(a_1,z_1)$,$(a_2,z_2),\ds ,(a_{n-1},z_{n-1})$, and $(v,w)$ can be permutated arbitrarily within $S$. Denote by $S^{L}$ (resp. $S^{R}$) the expansion of the $(n+1)$-point function $S$ with respect to the left (resp. right)-most term using \eqref{n3.4} (resp. \eqref{n3.6}). The second property is that
\begin{equation}\tag{{II}}\label{II}
\begin{aligned}
S(v_{3}',(b_{1},w_{1})(b_{2},w_{2})...(b_{n},w_{n})v_{2})
&=S^{L}(v_{3}',(b_{1},w_{1})(b_{2},w_{2})...(b_{n},w_{n})v_{2})\\
&=S^{R}(v_{3}',(b_{1},w_{1})(b_{2},w_{2})...(b_{n},w_{n})v_{2}),
\end{aligned}
\end{equation}
where $(b_{1},w_{1})$ in $S^{L}$ is not $(v,w)$, and $(b_{n},w_{n})$ in $S^{R}$ is not $(v,w)$.

Note that properties \ref{I} and \ref{II} are satisfied by the $4$-point and $5$-point functions (see \eqref{n4.10} and \eqref{n4.15}.) We construct $(n+3)$-point functions as follows: 
\begin{df}\label{df4.14}
	Assume the number of $V$ in the sub-indices of $S^{L}_{VV...M^{1}...V}$ and $S^{R}_{V...M^{1}...VV}$ are both equal to $n$, the sub-index $M^{1}$ in $S^{L}$ is not at the first place, and the sub-index $M^{1}$ in $S^{R}$ is not at the last place. We define $S^{L}_{VV...M^{1}...V}$ by
	\begin{align*}
	&S^{L}_{VV...M^{1}...V}(v_{3}',(a_{1},z_{1})...(v,w)...v_{2}):=S(v_{3}'o(a_{1}),(a_{2},z_{2})...(a_{n},z_{n})(v,w)v_{2})z_{1}^{-\wt a_{1}}\\
	&\ +\sum_{k=2}^{n}\sum_{j\geq 0} F_{\wt a_{1},j}(z_{1},z_{k}) S(v_{3}',(a_{2},z_{2})...(a_{1}(j)a_{k},z_{k})...(a_n,z_n)(v,w)v_{2})\numberthis\label{n4.16}\\
	&\ +\sum_{j\geq 0}F_{\wt a_{1},j}(z_{1},w) S(v_{3}',(a_{2},z_{2})...(a_{n},z_{n})(a_{1}(j)v,w)v_{2}),
	\end{align*}
	and define $S^{R}_{V...M^{1}...VV}$ by 
	\begin{align*}
	&S^{R}_{V...M^{1}...VV}(v_{3}',...(v,w)...(a_{1},z_{1})v_{2}):=S(v_{3}',(a_{2},z_{2})...(a_{n},z_{n})(v,w)o(a_{1})v_{2})z_{1}^{-\wt a_{1}}\\
	&\ +\sum_{k=2}^{n}\sum_{j\geq 0} G_{\wt a_{1},j}(z_{1},z_{k}) S(v_{3}',(a_{2},z_{2})...(a_{1}(j)a_{k},z_{k})...(a_n,z_n)(v,w)v_{2})\numberthis\label{n4.17}\\
	&\ +\sum_{j\geq 0}G_{\wt a_{1},j}(z_{1},w) S(v_{3}',(a_{2},z_{2})...(a_{n},z_{n})(a_{1}(j)v,w)v_{2}),
	\end{align*}
	where the $S$ on right-hand sides of \eqref{n4.16} and \eqref{n4.17} is the $(n+2)$-point function.
\end{df}

The definition above indicates that $S^{L}_{VMV...V}=S^{L}_{VVM...V}=\ds=S^{L}_{VV...VM}$, which is reasonable because the $(n+2)$-point function $S$ on the right-hand side of \eqref{n4.16} satisfies the locality property \eqref{I}. For a similar reason, we can also expect that $S^{R}_{MV...VV}=S^{R}_{VM...VV}=\dots=S^{R}_{V...VMV}$. We need to show that 
\begin{equation}\label{n4.18}
\begin{aligned}
&S^{L}_{V...M...V}(v_{3}',(a_{1},z_{1})...(v,w)...(a_{2},z_{2})v_{2})\\
&=S^{R}_{V...M...V}(v_{3}',(a_{1},z_{1})...(v,w)...(a_{2},z_{2})v_{2}),
\end{aligned}
\end{equation}
for all $S^{L}_{VV...M...V}$ and $ S^{R}_{V...M...VV}$.

Indeed, as we mentioned in Proposition \ref{prop4.10}, since \eqref{n4.16} is the generalization of (2.2.6) in \cite{Z}, by a similar argument as the proof of (2.2.11) in \cite{Z}, we have: 
\begin{equation}\label{n4.19}
\begin{aligned}
&S^{L}_{VV...M...V}(v_{3}',(a_{1},z_{1})(a_{2},z_{2})...(v,w)...v_{2})
\\&=S^{L}_{VV...M...V}(v_{3}',(a_{2},z_{2})(a_{1},z_{1})...(v,w)...v_{2}).
\end{aligned}
\end{equation}
\begin{prop}\label{prop4.15}
	If $S^{L}_{VV...M...V}(v_{3}',(a_{1},z_{1})...v_{2})=S^{R}_{V...M...VV}(v_{3}',...(a_{1},z_{1})v_{2})$, i.e. if the right-hand side of \eqref{n4.16} is equal to the right-hand side of\eqref{n4.17}, then \eqref{n4.18} holds.
\end{prop}
\begin{proof}
	The proof is similar to the proof of Proposition \ref{prop4.12}. By \eqref{n4.19} and the assumption,
	\begin{align*}
	&S^{L}_{V...M...V}(v_{3}',(a_{1},z_{1})...(v,w)...(a_{2},z_{2})v_{2})=S^{L}_{VV...M...V}(v_{3}',(a_{1},z_{1})(a_{2},z_{2})...(v,w)...v_{2})\\
	&=S^{L}_{VV...M...V}(v_{3}',(a_{2},z_{2})(a_{1},z_{1})...(v,w)...v_{2})=S^{R}_{V...M...VV}(v_{3}',(a_{1},z_{1})...(v,w)...(a_{2},z_{2})v_{2})
	\end{align*}
	as asserted.
\end{proof}
Now we are left to show that: 
\begin{equation}\label{n4.20} S^{L}_{VV...M...V}(v_{3}',(a_{1},z_{1})...(v,w)...v_{2})=S^{R}_{V...M...VV}(v_{3}',...(v,w)...(a_{1},z_{1})v_{2}).
\end{equation} 
Similar to the previous subsection, we use the symbols $(1),(2),(3)$, and $(4)$ to denote the following summands on the right-hand side of $\eqref{n4.16}-\eqref{n4.17}$:
\begin{equation}\tag*{(1)}\label{1}
S(v_{3}'o(a_{1}),(a_{2},z_{2})...(v,w)v_{2})z^{-\wt a_{1}}
-S(v_{3}',(a_{2},z_{2})...(v,w)o(a_{1})v_{2})z^{-\wt a_{1}}.
\end{equation}
\begin{equation}\tag*{(2)}\label{2}
\sum_{j\geq 0}(F_{\wt a_{1},j}(z_{1},z_{2})- G_{\wt a_{1},j}(z_{1},z_{2}))S(v_{3}',(a_{1}(j)a_{2},z_{2})...(a_n,z_n)(v,w)v_{2}).
\end{equation}
\begin{equation}\tag*{(3)}\label{3}
\sum_{k=3}^{n}\sum_{j\geq 0}(F_{\wt a_{1},j}(z_{1},z_{k})-G_{\wt a_{1},j}(z_{1},z_{k}))S(v_{3}',(a_{2},z_{2})...(a_{1}(j)a_{k},z_{k})...(v,w)v_{2}).
\end{equation}
\begin{equation}\tag*{(4)}\label{4}
\sum_{j\geq 0}((F_{\wt a_{1},j}(z_{1},w)-G_{\wt a_{1},j}(z_{1},w))S(v_{3}',(a_{2},z_{2})...(a_n,z_n)(a_{1}(j)v,w)v_{2}).
\end{equation}
Then we need to show that (1)+(2)+(3)+(4)=0.\par
Our strategy is to apply the expansion formula \eqref{n3.4} and expand each summand of $\ref{1}-\ref{4}$ with respect to the left-most term. Then we add them all up and show that the sum equals $0$. (Since we are using the recursive formula \eqref{n3.4} twice and the $3$-point function cannot be expanded, the construction of the 5-point function in the previous subsection is necessary for our induction process.)

Start with \ref{1}, 
note that $S(v_{3}'o(a_{1}),(a_{2},z_{2})...(a_{n},z_{n})(v,w)v_{2})z^{-\wt a_{1}}$ can be written as:
\begin{align*}\tag{$\ast$}
&\ \ S(v_{3}'o(a_{1})o(a_{2}),(a_{3},z_{3})...(a_{n},z_{n})(v,w)v_{2})z_{1}^{-\wt a_{1}}z_{2}^{-\wt a_{2}}\\
&+\sum_{t=3}^{n}\sum_{i\geq 0}F_{\wt a_{2},i}(z_{2},z_{t})S(v_{3}'o(a_{1}),(a_{3},z_{3})...(a_{2}(i)a_{t},z_{t})...(a_{n},z_{n})(v,w)v_{2})z_{1}^{-\wt a_{1}}\\
&+\sum_{i\geq 0}F_{\wt a_{2},i}(z_{2},w)S(v_{3}'o(a_{1}),(a_{3},z_{3})...(a_{n},z_{n})(a_{2}(i)v,w)v_{2})z_{1}^{-\wt a_{1}},
\end{align*}
and $S(v_{3}',(a_{2},z_{2})...(a_{n},z_{n})(v,w)o(a_{1})v_{2})z_{1}^{-\wt a_{1}}$ can be written as
\begin{align*}
&\ \ \tag{$\ast\ast$}S(v_{3}'o(a_{2}),(a_{3},z_{3})...(a_{n},z_{n})(v,w)o(a_{1})v_{2})z_{1}^{-\wt a_{1}}z_{2}^{-\wt a_{2}}\\ 
&+\sum_{t=3}^{n}\sum_{i\geq 0}F_{\wt a_{2},i}(z_{2},z_{t})S(v_{3}',(a_{3},z_{3})...(a_{2}(i)a_{t},z_{t})...(a_{n},z_{n})(v,w)o(a_{1})v_{2})z_{1}^{-\wt a_{1}}\\
&+\sum_{i\geq 0}F_{\wt a_{2},i}(z_{2},w)S(v_{3}',(a_{3},z_{3})...(a_{n},z_{n})(a_{2}(i)v,w)o(a_{1})v_{2})z_{1}^{-\wt a_{1}}.
\end{align*}
We denote the first, second, and third corresponding terms in $(\ast)-(\ast\ast)$ by $(11)$, $(12)$, and $(13)$, respectively. In particular, (11) is
\begin{equation}\tag*{(11)}\label{11}
\begin{aligned}
& S(v_{3}'o(a_{1})o(a_{2}),(a_{3},z_{3})...(a_{n},z_{n})(v,w)v_{2})z_{1}^{-
	\wt a_{1}}z_{2}^{-\wt a_{2}}\\
-&S(v_{3}'o(a_{2}),(a_{3},z_{3})...(a_{n},z_{n})(v,w)o(a_{1})v_{2})z_{1}^{-\wt a_{1}}z_{2}^{-\wt a_{2}}.
\end{aligned}
\end{equation}
\begin{lm}\label{lm4.8}
	As $(n+1)$-point function, we have:
	\begin{align*}
	&S(v_{3}'o(a_{1}),(a_{3},z_{3})...(a_{n},z_{n})(v,w)v_{2})-S(v_{3}',(a_{3},z_{3})...(a_{n},z_{n})(v,w)o(a_{1})v_{2})\\
	&=\sum_{k=3}^{n}\sum_{j\geq 0}\binom{\wt a_{1}-1}{j}z_{k}^{\wt a_{1}-j-1}S(v_{3}',(a_{3},z_{3})...(a_{1}(j)a_{k},z_{k})...(a_{n},z_{n})(v,w)v_{2})\numberthis\label{n4.21}\\
	&+\sum_{j\geq 0} \sum_{j\geq 0}\binom{\wt a_{1}-1}{j}w^{\wt a_{1}-j-1}S(v_{3}',(a_{3},z_{3})...(a_{n},z_{n})(a_{1}(j)v,w)v_{2})
	\end{align*}
\end{lm}
\begin{proof}
	By the induction hypothesis for the $(n+2)$-point functions and \eqref{n3.8}, we have: 
	\begin{align*}
	0&=S(v_{3}',(a_{1},z_{1})(a_{3},z_{3})...(a_{n},z_{n})(v,w)v_{2})-S(v_{3}',(a_{3},z_{3})...(a_{n},z_{n})(a_{1},z_{1})(v,w)v_{2})\\
	&=S(v_{3}'o(a_{1})(a_{3},z_{3})...(a_{n},z_{n})(v,w)v_{2})z_{1}^{-wta_{1}}-S(v_{3}'(a_{3},z_{3})...(a_{n},z_{n})(v,w)o(a_{1})v_{2})z_{1}^{-\wt a_{1}}\\
	&\ \ +\sum_{k=3}^{n}\sum_{j\geq 0}(F_{\wt a_{1},j}(z_{1},z_{k})-G_{\wt a_{1},j}(z_{1},z_{k}))S(v_{3}',(a_{3},z_{3})...(a_{1}(j)a_{k},z_{k})...(a_{n},z_{n})(v,w)v_{2})\\
	&\ \ +\sum_{j\geq 0} (F_{\wt a_{1},j}(z_{1},w)-G_{\wt a_{1},j}(z_{1},w))S(v_{3}',(a_{3},z_{3})...(a_{n},z_{n})(a_{1}(j)v,w)v_{2})\\
	&=S(v_{3}'o(a_{1})(a_{3},z_{3})...(a_{n},z_{n})(v,w)v_{2})z_{1}^{-\wt a_{1}}-S(v_{3}'(a_{3},z_{3})...(a_{n},z_{n})(v,w)o(a_{1})v_{2})z_{1}^{-\wt a_{1}}\\
	&\ \ +\sum_{k=3}^{n}\sum_{j\geq 0}-\binom{\wt a_{1}-1}{j}z_{k}^{\wt a_{1}-j-1}S(v_{3}',(a_{3},z_{3})...(a_{1}(j)a_{k},z_{k})...(a_{n},z_{n})(v,w)v_{2})\\
	&\ \ +\sum_{j\geq 0} -\binom{\wt a_{1}-1}{j}w^{\wt a_{1}-j-1}S(v_{3}',(a_{3},z_{3})...(a_{n},z_{n})(a_{1}(j)v,w)v_{2}).
	\end{align*}
	This proves \eqref{n4.21}. 
\end{proof}

It follows from the Lemma \ref{lm4.8} that (12) and (13) can be written as: 
\begin{align*}
(12)&=\sum_{t=3}^{n}\sum_{k=3,k\neq t}^{n}\sum_{i,j\geq 0} F_{\wt a_{2},i}(z_{2},z_{t})\binom{\wt a_{1}-1}{j}z_{1}^{-\wt a_{1}}z_{k}^{\wt a_{1}-1-j}\\
&\qquad\cdot \underset{(121)}{S(v_{3}',(a_{3},z_{3})...(a_{1}(j)a_{k},z_{k})...(a_{2}(i)a_{t},z_{t})...(a_{n},z_{n})(v,w)v_{2})}\\
&\ \ + \sum_{t=3}^{n}\sum_{i,j\geq 0}F_{\wt a_{2},i}(z_{2},z_{t})\binom{\wt a_{1}-1}{j}z_{1}^{-\wt a_{1}}z_{t}^{\wt a_{1}-1-j}\\
&\qquad \cdot \underset{(122)}{S(v_{3}',(a_{3},z_{3})...(a_{1}(j)a_{2}(i)a_{t},z_{t})...(a_{n},z_{n})(v,w)v_{2})}\\
&\ \ +\sum_{t=3}^{n}\sum_{i,j\geq 0}F_{\wt a_{2},i}(z_{2},w)\binom{\wt a_{1}-1}{j}z_{1}^{-\wt a_{1}}w^{\wt a_{1}-1-j}\\
&\qquad \cdot  \underset{(123)}{S(v_{3}',(a_{3},z_{3})...(a_{2}(i)a_{t},z_{t})...(a_{n},z_{n})(a_{1}(j)v,w)v_{2})}\\
&=(121)+(122)+(123),\\
(13)&=\sum_{k=3}^{n}\sum_{i,j\geq 0}F_{\wt a_{2},i}(z_{2},z_{k})\binom{\wt a_{1}-1}{j}z_{1}^{-wta_{1}}z_{k}^{\wt a_{1}-1-j}\\ 
&\qquad \cdot  \underset{(131)}{S(v_{3}',(a_{3},z_{3})...(a_{1}(j)a_{k},z_{k})...(a_{n},z_{n})(a_{2}(i)v,w)v_{2})}\\
&\ \ +\sum_{i,j\geq 0}F_{\wt a_{2},i}(z_{2},w)\binom{\wt a_{1}-1}{j}z_{1}^{-\wt a_{1}}w^{\wt a_{1}-1-j}\\ 
&\qquad \cdot \underset{(132)}{ S(v_{3}',(a_{3},z_{3})...(a_{n},z_{n})(a_{1}(j)a_{2}(i)v,w)v_{2})}\\
&=(131)+(132).
\end{align*}
Then \ref{1}=(11)+(121)+(122)+(123)+(131)+(132).

Now we expand \ref{2}, \ref{3}, and \ref{4} with respect to their corresponding left-most terms. By \eqref{n3.8}, they can be expressed as follows:

\begin{align*}
(2)=&\sum_{j\geq 0}-\binom{\wt a_{1}-1}{j} z_{1}^{-\wt a_{1}} \underset{(21)}{z_{2}^{-\wt a_{2}} S(v_{3}'o(a_{1}(j)a_{2}),(a_{3},z_{3})...(a_{n},z_{n})(v,w)v_{2})}\\
&+\ \ \sum_{k=3}^{n}\sum_{i,j\geq 0}-\binom{\wt a_{1}-1}{j} z_{1}^{-\wt a_{1}}z_{2}^{\wt a_{1}-1-j}F_{\wt a_{1}+\wt a_{2}-j-1,i}(z_{2},z_{k}) \\
&\qquad\cdot  \underset{(22)}{S(v_{3}',(a_{3},z_{3})...((a_{1}(j)a_{2})(i)a_{k},z_{k})...(a_{n},z_{n})(v,w)v_{2})}\\
&+\ \ \sum_{i,j\geq 0}-\binom{\wt a_{1}-1}{j}z_{1}^{-\wt a_{1}} z_{2}^{-\wt a_{1}-1-j}F_{\wt a_{1}+\wt a_{2}-j-1,i}(z_{2},w)\\
&\qquad\cdot  \underset{(23)}{S(v_{3}',(a_{3},z_{3})...(a_{n},z_{n})((a_{1}(j)a_{2})(i)v,w)v_{2})}\\
&=(21)+(22)+(23).\\
(3)=&\sum_{k=3}^{n}\sum_{j\geq 0}-\binom{wta_{1}-1}{j}\underset{(31)}{z_{1}^{-\wt a_{1}}z_{k}^{\wt a_{1}-1-j} S(v_{3}'o(a_{2}),(a_{3},z_{3})...(a_{1}(j)a_{k},z_{k})...v_{2})z_{2}^{-\wt a_{2}}}\\
&+\sum_{k=3}^{n}\sum_{j,i\geq 0}-\binom{\wt a_{1}-1}{j}z_{1}^{-\wt a_{1}}z_{k}^{\wt a_{1}-1-j}F_{\wt a_{2},i}(z_{2},w)\\
&\qquad \cdot  \underset{(32)}{S(v_{3}',(a_{3},z_{3})...(a_{1}(j)a_{k},z_{k})...(a_{n},z_{n})(a_{2}(i)v,w)v_{2})}\\
&+\sum_{k=3}^{n}\sum_{t=3,t\neq k}^{n}\sum_{j,i\geq 0}-\binom{\wt a_{1}-1}{j}z_{1}^{-\wt a_{1}}z_{k}^{\wt a_{1}-1-j}F_{\wt a_{2},i}(z_{2},z_{t})\\
&\qquad\cdot  \underset{(33)}{S(v_{3}',(a_{3},z_{3})...(a_{2}(i)a_{t},z_{t})...(a_{1}(j)a_{k},z_{k})...(a_{n},z_{n})(v,w)v_{2})}\\
&+\sum_{k=3}^{n}\sum_{j,i\geq 0}-\binom{\wt a_{1}-1}{j}z_{1}^{-\wt a_{1}}z_{k}^{\wt a_{1}-1-j}F_{\wt a_{2},i}(z_{2},z_{k})\\
&\qquad\cdot  \underset{(34)}{S(v_{3}',(a_{3},z_{3})...(a_{2}(i)a_{1}(j)a_{k},z_{k})...(a_{n},z_{n})(v,w)v_{2})}\\
&=(31)+(32)+(33)+(34).\\
(4)&=\sum_{j\geq 0}-\binom{\wt a_{1}-1}{j}\underset{(41)}{z_{1}^{-
		\wt a_{1}}w^{
		\wt a_{1}-1-j}S(v_{3}'o(a_{2}),(a_{3},z_{3})...(a_{n},z_{n})(a_{1}(j)v,w)v_{2})z_{2}^{-\wt a_{2}}}\\
&+\sum_{t=3}^{n}\sum_{j,i\geq 0}-\binom{\wt a_{1}-1}{j}z_{1}^{-\wt a_{1}}w^{\wt a_{1}-1-j}F_{\wt a_{2},i}(z_{2},z_{k})\\
&\qquad \cdot \underset{(42)}{S(v_{3}',(a_{3},z_{3})...(a_{2}(i)a_{k},z_{k})...(a_{n},z_{n})(a_{1}(j)v,w)v_{2})}\\
&+\sum_{j,i\geq 0}-\binom{\wt a_{1}-1}{j}z_{1}^{-\wt a_{1}}w^{\wt a_{1}-1-j}F_{\wt a_{2},i}(z_{2},w)\\
&\qquad \cdot \underset{(43)}{S(v_{3}',(a_{3},z_{3})...(a_{n},z_{n})(a_{2}(i)a_{1}(j)v,w)v_{2})}\\
&=(41)+(42)+(43).
\end{align*}

By Lemma \ref{lm4.5} and the formula \eqref{n3.5} of $\iota_{z_{2},z_{t}}F_{n,i}(z_{2},z_{t})$, we have: 
\begin{align*}
&\sum_{i,j\geq 0}\binom{\wt a_{1}-1}{j} F_{\wt a_{2},i}(z_{2},z_{t})a_{1}(j)a_{2}(i)a_{t}+\sum_{i,j\geq 0} -\binom{\wt a_{1}-1}{j}F_{\wt a_{2},i}(z_{2},z_{t})a_{1}(j)a_{2}(i)a_{t}\\
&+\sum_{i,j\geq 0}-\binom{\wt a_{1}-1}{j} F_{\wt a_{1}+
	\wt a_{2}-j-1,i}(z_{2},z_{t})(a_{1}(j)a_{2})(i)a_{t}\numberthis \label{n4.22}\\
&=0,
\end{align*}
and the same equation holds if we replace $z_{t}$ with $w$ and $a_{i}$ with $v$. Using \eqref{n4.22}, we have the cancelations $(122)+(22)+(34)=0,$ and $(132)+(23)+(43)=0.$
Moreover, it follows directly from the expressions of the terms $(123),(42), (121),(33),(131)$, and $(32)$ that
$$(123)+(42)=0,\quad (121)+(33)=0,\ \ \textrm{and}\ \  (131)+(32)=0.$$
Now it remains to show \ref{11}+(21)+(31)+(41)=0, or equivalently, 
\begin{align*}
&S(v_{3}'o(a_{1})o(a_{2}),(a_{3},z_{3})...(a_{n},z_{n})(v,w)v_{2})-S(v_{3}'o(a_{2}),(a_{3},z_{3})...(a_{n},z_{n})(v,w)o(a_{1})v_{2})\\
&=\sum_{j\geq 0}\binom{\wt a_{1}-1}{j}  S(v_{3}'o(a_{1}(j)a_{2}),(a_{3},z_{3})...(a_{n},z_{n})(v,w)v_{2})\\
&+\sum_{k=3}^{n}\sum_{j\geq 0}\binom{\wt a_{1}-1}{j}z_{k}^{\wt a_{1}-1-j} S(v_{3}'o(a_{2}),(a_{3},z_{3})...(a_{1}(j)a_{k},z_{k})...(v,w)v_{2})\numberthis \label{n4.23}\\
&+\sum_{j\geq 0}\binom{\wt a_{1}-1}{j}w^{\wt a_{1}-1-j}S(v_{3}'o(a_{2}),(a_{3},z_{3})...(a_{n},z_{n})(a_{1}(j)v,w)v_{2}),
\end{align*}
but this is a consequence of Lemma \ref{lm4.8}. In fact,
\begin{align*}
&L.H.S. \ of \ \eqref{n4.23}\\
&=S(v_{3}'o(a_{1})o(a_{2}),(a_{3},z_{3})...(a_{n},z_{n})(v,w)v_{2})-S(v_{3}'o(a_{2})o(a_{1}),(a_{3},z_{3})...(a_{n},z_{n})(v,w)v_{2})\\
&\ +S(v_{3}'o(a_{2})o(a_{1}),(a_{3},z_{3})...(a_{n},z_{n})(v,w)v_{2})-S(v_{3}'o(a_{2}),(a_{3},z_{3})...(a_{n},z_{n})(v,w)o(a_{1})v_{2}).
\end{align*}
Since $S$ is linear in the place $M^{3}(0)^{*}$, we have
\begin{align*}
&\ \ \ S(v_{3}'o(a_{1})o(a_{2}),(a_{3},z_{3})...(a_{n},z_{n})(v,w)v_{2})-S(v_{3}'o(a_{2})o(a_{1}),(a_{3},z_{3})...(a_{n},z_{n})(v,w)v_{2})\\
&=S(v_{3}'[o(a_{1}),o(a_{2})],(a_{3},z_{3})...(a_{n},z_{n})(v,w)v_{2})\\
&=\sum_{j\geq 0}\binom{\wt a_{1}-1}{j}S(v_{3}'o(a_{1}(j)a_{2}), (a_{3},z_{3})...(a_{n},z_{n})(v,w)v_{2}),
\end{align*}
which is the first term on the right-hand side of \eqref{n4.23}. Moreover, by Lemma \ref{lm4.8}, 
\begin{align*}
&S(v_{3}'o(a_{2})o(a_{1}),(a_{3},z_{3})...(a_{n},z_{n})(v,w)v_{2})-S(v_{3}'o(a_{2}),(a_{3},z_{3})...(a_{n},z_{n})(v,w)o(a_{1})v_{2})\\
&=\sum_{k=3}^{n}\sum_{j\geq 0}\binom{\wt a_{1}-1}{j}z_{k}^{\wt a_{1}-1-j} S(v_{3}'o(a_{2}),(a_{3},z_{3})...(a_{1}(j)a_{k},z_{k})...(v,w)v_{2})\\
&+\sum_{j\geq 0}\binom{\wt a_{1}-1}{j}w^{\wt a_{1}-1-j}S(v_{3}'o(a_{2}),(a_{3},z_{3})...(a_{n},z_{n})(a_{1}(j)v,w)v_{2}),
\end{align*}
which gives us the last two summands on the right-hand side of \eqref{n4.23}. This proves \eqref{n4.23}. Hence $\ref{1}+\ref{2}+\ref{3}+\ref{4}=0$, and so \eqref{n4.20} holds. 

Then by Proposition \ref{prop4.15}, all the $(n+3)$-point functions $S_{VV...M...V}^{L}$ and $S_{V...M...VV}^{R}$ defined by \eqref{n4.16} and \eqref{n4.17} give rise to one single $(n+3)$-point function:
\begin{equation}\label{n4.24}
S: M^{3}(0)^{*}\times V\times\ds\times M^{1}\times\ds\times V\times M^{2}(0) \rightarrow \mathcal{F}(z_{1},\ds,z_{n},w),
\end{equation}
where $M^1$ can be placed anywhere in between the first and the last place of $V$. Moreover, by Definition \ref{df4.14} and \eqref{n4.18}, $S$ in \eqref{n4.24} satisfies the locality \ref{I} and the expansion property \ref{II}, with $n$ replaced by $n+1$.  Therefore, the induction step is complete. 

\begin{thm}\label{thm4.17}
	The system of $(n+3)$-point functions $S$ we constructed by Definitions \ref{df4.9}, \ref{df4.11}, and \ref{df4.14} in this subsection lies in $ \Cor\fusion{M^1}{M^{2}(0)}{M^{3}(0)}$.
\end{thm}
\begin{proof}
	Since S is constructed inductively by the recursive formulas \eqref{n3.4} and \eqref{n3.6} in view of Defintions \ref{df4.9}, \ref{df4.11}, and \ref{df4.14}, it obviously satisfies \eqref{n3.4} and \eqref{n3.6}. By \eqref{n4.7}, we have  $S(v_{3}',(v,w)v_{2})=f(v_{3}'\otimes v\otimes v_{2})w^{-\deg v}$, for any $v_{3}'\in M^3(0)^\ast, v\in M^{1},$ and $v_{2}\in M^{2}(0)$. By the Hom-tensor duality, we have a well-defined element $f_{v}\in \Hom_\C(M^2(0),M^3(0))$ such that  $\<v_{3}',f_{v}(v_2)\>=f(v_{3}'\otimes v\otimes v_{2})$ for each $v\in M^1$. Hence $S$ satisfies \eqref{n3.3}.
	
	In view of Definition \ref{df3.1}, it remains to show that $S$ satisfies $(2)-(6)$ in Definition \ref{df2.1} for $v_2\in M^2(0)$ and $v_{3}'\in M^3(0)^\ast$. Indeed, the locality follows from \eqref{I}, and by \eqref{n4.16}, 
	\begin{align*}
	&S(v_{3}', (\vac,z)(a_{1},z_{1})...(a_{n},z_{n})(v,w)v_{2})=S(v_{3}'o(\vac), (a_{1},z_{1})...(a_{n},z_{n})(v,w)v_{2}) z^{-\wt \vac}\\
	&\ \ \ +\sum_{k=1}^{n} \sum_{j\geq 0}F_{\wt \vac, j}(z,z_{j})S(v_{3}', (a_{1},z_{1})...(\vac(j)a_{k},z_{k})...(a_{n},z_{n})(v,w)v_{2})\\
	&\ \ \ +\sum_{j\geq 0} F_{\wt \vac, j}(z,w) S(v_{3}', (a_{1},z_{1})...(a_{n},z_{n})(\vac (j)v,w)v_{2})\\
	&=S(v_{3}', (a_{1},z_{1})...(a_{n},z_{n})(v,w)v_{2}),
	\end{align*}
	since $\vac(j) a_{k}=\vac(j)v=0$ when $j\geq 0$, and $o(\vac)=\Id$.
	
	Again because $S$ in \eqref{n4.24} satisfies \eqref{n4.16}, it is easy to verify the following associativity formulas by a similar argument to the proof of (2.2.9) in \cite{Z}: 
	\begin{equation}\label{4.23}
	\begin{aligned}
	&\int_{C}S(v_{3}',(a_{1},z_{1})(v,w)...(a_{n},z_{n})v_{2})(z_{1}-w)^{n} d{z_{1}}=S(v_{3}',(a_{1}(k)v,w)...(a_{n},z_{n})v_{2}),\\
	&\int_{C}S(v_{3}',(a_{1},z_{1})(a_{2},z_{2})...(v,w)v_{2})(z_{1}-z_{2})^{n} d{z_{1}}=S(v_{3}',(a_{1}(k)a_{2},z_{2})...(v,w)v_{2}),
	\end{aligned}
	\end{equation}
	where in the first equation of \eqref{4.23}, $C$ is a contour of $z_{1}$ surrounding $w$ with $z_{2},...,z_{n}$ outside of $C$; while in the second equation of \eqref{4.23}, $C$ is a contour of $z_{1}$ surrounding $z_{2}$ with $z_{3},...,z_{n}, w$ outside of $C$. We also have:
	\begin{equation}\label{n4.26}
	\begin{aligned}
	&S(v_{3}',(L(-1)a_{1},z_{1})...(a_{n},z_{n})(v,w)v_{2})=\frac{d}{dz_{1}}S(v_{3}',(a_{1},z_{1})...(a_{n},z_{n})(v,w)v_{2}),\\
	&S(v_{3}',(L(-1)v,w)(a_{1},z_{1})...(a_n,z_n)v_{2})w^{-h}=\frac{d}{dw}(S(v_{3}',(v,w)(a_{1},z_{1})...v_{2})w^{-h}).
	\end{aligned}
	\end{equation}
	The first equation in \eqref{n4.26} is similar to (2.2.8) in \cite{Z}. We omit the details of the proof. To show the second equation in \eqref{n4.26}, we use induction on $n$. When $n=0$, by \eqref{4.5'} and Lemma \ref{lm4.5'}, we have: $L(-1)v+(L(0)+h_{2}-h_{3})v\equiv 0 \mod O_{h}(M^{1})$ for all $v\in M^{1}$. Then 
	\begin{align*}
	&S(v_{3}',(L(-1)v,w)v_{2}) w^{-h}=f(v_{3}'\otimes L(-1)v\otimes v_{2})w^{-\deg v-1-h}\\
	&=-f(v_{3}'\otimes (L(0)+h_{2}-h_{3})v\otimes v_{2}) w^{-\deg v-1-h}=f(v_{3}'\otimes v\otimes v_{2}) \frac{d}{dw} (w^{-\deg v-h})\numberthis \label{n4.27}\\
	&=\frac{d}{dw}(S(v_{3}', (v,w)v_{2})w^{-h}).
	\end{align*}
	Now assume the second equation of \eqref{n4.26} holds for the $(n+2)$-point function, then by the properties \eqref{I} and \eqref{II} of $S$, we have:
	\begin{align*} &S(v_{3}',(L(-1)v,w)(a_{1},z_{1})...(a_{n},z_{n})v_{2})w^{-h}=S^{L}(v_{3}',(a_{1},z_{1})...(a_{n},z_{n})(L(-1)v,w)v_{2})w^{-h}\\
	&=S(v_{3}'o(a_{1}),(a_{2},z_{2})...(a_{n},z_{n})(L(-1)v,w)v_{2})z_{1}^{-\wt a_{1}}w^{-h}\\
	&\ \ + \sum_{k=2}^{n}\sum_{j\geq 0}F_{\wt a_{1},j}(z_{1},z_{k})S(v_{3}',(a_{2},z_{2})...(a_{1}(j)a_{k},z_{k})...(a_{n},z_{n})(L(-1)v,w)v_{2})w^{-h}\numberthis \label{n4.28}\\
	&\ \ +\sum_{j\geq 0}F_{\wt a_{1},j}(z_{1},w) S(v_{3}',(a_{2},z_{2})...(a_{n},z_{n})(a_{1}(j)L(-1)v,w)v_{2})w^{-h}.
	\end{align*} 
	Note that we can apply the induction hypothesis to the first two terms of \eqref{n4.28}. Moreover, by the $L(-1)$-bracket formula (4.2.1) in \cite{FHL}, we have: 
	$$a_{1}(j)L(-1)v_{2}=L(-1)a_{1}(j)v_{2}-[L(-1),a_{1}(j)]v_{2}=L(-1)a_{1}(j)v_{2}+ja_{1}(j-1)v_{2}.$$ It follows from the induction hypothesis and \eqref{n3.5} that 
	\begin{align*}
	&\sum_{j\geq 0}F_{\wt a_{1},j}(z_{1},w) S(v_{3}',(a_{2},z_{2})...(a_{n},z_{n})(a_{1}(j)L(-1)v,w)v_{2})w^{-h}\\
	&=\sum_{j\geq 0}F_{\wt a_{1},j}(z_{1},w) \frac{d}{dw}(S(v_{3}',(a_{2},z_{2})...(a_{n},z_{n})(a_{1}(j)v,w)v_{2})w^{-h})\\
	&\ \ +\sum_{j\geq 1}\frac{z_{1}^{-\wt a_{1}}}{(j-1)!}\bigg(\frac{d}{dw}\bigg)^{j}\bigg(\frac{w^{\wt a_{1}}}{z_{1}-w}\bigg) S(v_{3}',(a_{2},z_{2})...(a_{n},z_{n})(a_{1}(j-1)v,w)v_{2})w^{-h}\\
	&=\frac{d}{dw}\sum_{j\geq 0}F_{\wt a_{1},j}(z_{1},w) S(v_{3}',(a_{2},z_{2})...(a_{n},z_{n})(a_{1}(j)v,w)v_{2})w^{-h}.
	\end{align*}This proves \eqref{n4.26}. Finally, let $v_{3}'\in M^{3}(0)^{\ast},v\in M^{1},v_{2}\in M^{2}(0)$, and $a_{1},\ds,a_{n}\in V$ be highest weight vectors of the Virasoro algebra. By property \eqref{I} and \eqref{n4.26} of $S$, we have:
	\begin{align*}
	&S(v_{3}', (\w,x)(\w,x_{1})...(\w,x_{m})(a_{1},z_{1})...(a_{n},z_{n})(v,w)v_{2})\\
	&= S(v_{3'}, (\w,x_{1})...(a_{n},z_{n})(v,w)o(\w)v_{2}) x^{-2}\\
	&\ \ +\sum_{k=1}^{m} \sum_{j\geq 0} G_{2,j}(x,x_{k})S(v_{3}', (\w,x_{1})...(\w_{j}\w,x_{k})...(a_{n},z_{n})(v,w)v_{2})\\
	&\ \ +\sum_{k=1}^{n} \sum_{j\geq 0} G_{2,j}(x,z_{k}) S(v_{3}', (\w,x_{1})...(\w_{j}a_{k},z_{k})...(a_{n},z_{n})(v,w)v_{2})\\
	&\ \ +\sum_{j\geq 0} G_{2,j}(x,w)S(v_{3}', (\w,x_{1})...(a_{n},z_{n})(\w_{j}v,w)v_{2}).
	\end{align*}
	By the definition formula \eqref{n3.7}, it is easy to verify that: 
	$$
	G_{2,0}(x,z)=\frac{x^{-1}z}{x-z},\quad G_{2,1}(x,z)=\frac{1}{(x-z)^{2}}, \quad
	G_{2,3}(x,z)=\frac{1}{(x-z)^{4}}.
	$$
	Then by using the properties of the Virasoro element $\w$ (see Section 2.3 in \cite{FHL}), we have: 
	\begin{align*}
	&S(v_{3}',(\w,x)(\w,x_{1})...(\w,x_{m})(a_{1},z_{1})...(v,w)...(a_{n},z_{n})v_{2})\\
	&=\sum_{k=1}^{n}\frac{x^{-1}z_{k}}{x-z_{k}}\frac{d}{dz_{k}}S+\sum_{k=1}^{n}\frac{\wt a_{k}}{(x-z_{k})^{2}}S+\frac{x^{-1}w}{x-w}w^{h}\frac{d}{dw}(S\cdot w^{-h})+\frac{\wt v}{(x-w)^{2}}S\\
	&\ \ +\frac{ h_{2}}{x^{2}}S+\sum_{k=1}^{m}\frac{x^{-1}wx_{k}}{x-x_{k}} \frac{d}{dx_{k}}S+\sum_{k=1}^{m}\frac{2}{(x-x_{k})^{2}}S\\
	&\ \ +\frac{c}{2}\sum_{k=1}^{m}\frac{1}{(x-x_{k})^{4}} S(v_{3}',(\w,x_{1})...\widehat{(\w,x_{k})}...(\w,x_{m})(a_{1},z_{1})...(v,w)...(a_{n},z_{n})v_{2}),
	\end{align*}
	where $S=S(v_{3}',(\w,x_{1})...(\w,x_{m})(a_{1},z_{1})...(a_{n},z_{n})(v,w)v_{2})$. This shows that the $S$ in \eqref{n3.24} also satisfies \eqref{n2.8}, with $v_{3}'\in M^3(0)^\ast$ and $v_2\in M^2(0)$. Therefore, $S\in \Cor\fusion {M^{1}}{M^{2}(0)}{M^{3}(0)}$.
\end{proof}

\begin{remark}
	By equation \eqref{n4.27}, we see that it is necessary to have the equality $L(-1)v+(L(0)+h_{2}-h_{3})v=0$ hold in the bimodule $B_{h}(M^{1})$ to show the $L(-1)$-derivation property \eqref{n4.26} of $S$. However, in general, such equality does not hold in the bimodule $A(M^{1})$ in \cite{FZ} by its construction. This is the reason why $I\fusion{M^{1}}{M^{2}}{M^{3}}$ is not isomorphic to $(M^{3}(0)^{\ast}\otimes _{A(V)}A(M^{1})\otimes _{A(V)}M^{2}(0))^{\ast}$ in general.
\end{remark}
Theorem \ref{thm4.17} indicates that we have a well-defined linear map:
\begin{equation}\label{n4.29}
\mu: (M^{3}(0)^{\ast}\otimes_{A(V)}B_{h}(M^{1})\otimes _{A(V)} M^{2}(0))^{\ast}\ra \Cor\fusion{M^{1}}{M^{2}(0)}{M^{3}(0)},\quad 
f\mapsto S_{f},
\end{equation}
where $S_{f}$ is the $S$ we constructed in this subsection by Defintions \ref{df4.9}, \ref{df4.11}, and \ref{df4.14}.

Since we have $S_{f}(v_{3}',(v,w)v_{2})=f(v_{3}'\otimes v\otimes v_{2}) w^{-\deg v}$ by \eqref{n4.7}, and $f_{S_{f}}(v_{3}'\otimes v\otimes v_{2})w^{-\deg v}=S_{f}(v_{3}',(v,w)v_{2})$ by \eqref{4.6'} and Definition \ref{df3.1}, then $f_{S_{f}}=f$. i.e., $\nu \mu=1$. On the other hand, for $S\in \Cor\fusion{M^{1}}{M^{2}(0)}{M^{3}(0)}$, again by \eqref{n4.7} and \eqref{4.6'}, we have:
$$S_{f_{S}}(v_{3}',(v,w)v_{2})=f_{S}(v_{3}'\otimes v\otimes v_{2})w^{-\deg v}=S(v_{3}',(v,w)v_{2}).$$
Moreover, $S_{f_{S}}$ and $S$ satisfy the same recursive formulas by \eqref{n4.16}, \eqref{n4.17}, \eqref{n3.4}, and \eqref{n3.6}, then it follows from an easy induction that $S_{f_{S}}=S$. i.e., $\mu \nu =1$, and so $\mu$ is an isomorphism. Now we have our main result: 

\begin{thm}\label{thm4.19}
	Let $M^{1},M^{2}$, and $M^{3}$ be $V$-modules, with conformal weight $h_{1}$, $h_{2}$, and $h_{3}$, respectively. Assume $M^{2}(0)$ and $M^{3}(0)$ are irreducible $A(V)$-modules. Then we have the following isomorphism of vector spaces:
	\begin{equation}\label{n4.30}
	\begin{aligned}
	I\fusion{M^1}{\bar{M}(M^{2}(0))}{\bar{M}(M^{3}(0)^{\ast})'}\cong I\fusion{M^{1}}{\bar{M^{2}}}{\bar{M^{3}}}&\cong 	(M^{3}(0)^{\ast}\otimes_{A(V)}B_{h}(M^{1})\otimes _{A(V)} M^{2}(0))^{\ast},\\
	I&\mapsto f_{I},\quad f_{I}(v_{3}'\otimes v\otimes v_{2})=\<v_{3}', o(v)v_{2}\>,
	\end{aligned}
	\end{equation}
	for all $v_{3}'\in M^{3}(0)^{\ast}$, $v\in M^{1}$, and $v_{2}\in M^{2}(0)$, where $h=h_{1}+h_{2}-h_{3}$, and $M^{2}=\bar{M}/\Rad(\bar{M})$ and $M^{3}=(\tilde{M}/\Rad{\tilde{M}})'$ are quotient modules of the generalized Verma module $\bar{M}(M^{2}(0))$ and $\bar{M}(M^{3}(0))$, respectively.
\end{thm}
\begin{proof}
	This is a direct consequence of Corollary \ref{coro2.6}, Theorem \ref{thm3.14}, and Theorem \ref{thm4.17}, together of which give us the isomorphism: $I\fusion{M^1}{\bar{M}(M^{2}(0))}{\bar{M}(M^{3}(0)^{\ast})'}\cong I\fusion{M^{1}}{\bar{M^{2}}}{\bar{M^{3}}}\cong \Cor\fusion{M^{1}}{M^{2}(0)}{M^{3}(0)}\cong	(M^{3}(0)^{\ast}\otimes_{A(V)}B_{h}(M^{1})\otimes _{A(V)}M^{2}(0))^{\ast},$
	such that $I\mapsto f_{I}$ as in \eqref{n4.30}.
\end{proof}

Recall that $V$-modules $\bar{M^2}$ and $\bar{M^3}'$ are irreducible if condition \eqref{n3.25} is satisfied (see Proposition \ref{prop3.11}). By the isomorphism \eqref{n4.29}, condition \eqref{n3.25} translates to the following: 

For any $f\in (M^{3}(0)^{\ast}\otimes_{A(V)}B_{h}(M^{1})\otimes _{A(V)} M^{2}(0))^{\ast}$, one has: 
\begin{equation}\label{n4.31}
\sum_{i\geq 0}\binom{n}{i} f(v_{3}'\otimes b(i)v\otimes v_{2})=0,
\end{equation}	
for all $b\in V$, $n\in \Z$ such that $\wt b-n-1>0$, $v\in M^{1}$, $v_{3}'\in M^{3}(0)^{\ast}$, and $v_{2}\in M^{2}(0)$. 

\begin{coro}
	Let $M^{1},M^{2}$, and $M^{3}$ be $V$-modules, with conformal weight $h_{1}$, $h_{2}$, and $h_{3}$, respectively. Suppose $M^2$ and $M^3$ are irreducible, and condition \eqref{n4.31} is satisfied, then we have an isomorphism: 
	$I\fusion {M^1}{M^2}{M^3}\cong (M^{3}(0)^{\ast}\otimes_{A(V)}B_{h}(M^{1})\otimes _{A(V)} M^{2}(0))^{\ast}.$
\end{coro}

Suppose $M^2$ and $M^3$ are $V$-modules (not necessarily irreducible) that are generated by their corresponding bottom levels $M^2(0)$ and $M^3(0)$, which are irreducible $A(V)$-modules. Then by \eqref{n3.41} and \eqref{n4.30}, we have the following estimate of the fusion rule: 
\begin{equation}\label{n4.32}
\dim I\fusion {M^1}{M^2}{M^3}\leq \dim 	(M^{3}(0)^{\ast}\otimes_{A(V)}B_{h}(M^{1})\otimes _{A(V)} M^{2}(0))^{\ast}.
\end{equation}
Finally, when $V$ is rational, by Theorem \ref{thm4.19} and Corollary \ref{coro3.15}, we have: 
\begin{coro}\label{coro4.19}
	Let $V$ be a rational VOA, and let $M^{1}$, $M^{2}$, and $M^{3}$ be $V$ modules, with conformal weight $h_{1}$, $h_{2}$, and $h_{3}$, respectively. Suppose $M^{2}$ and $M^{3}$ are irreducible, then 
	\begin{equation}\label{n4.33}
	I\fusion{M^{1}}{M^{2}}{M^{3}}\cong (M^{3}(0)^{\ast}\otimes_{A(V)}B_{h}(M^{1})\otimes _{A(V)} M^{2}(0))^{\ast}.
	\end{equation}
\end{coro}

\subsection{Examples} In this subsection, we will use \eqref{n4.30} and the estimating formula \eqref{n4.32} and compute the fusion rules for certain modules over the Virasoro VOAs and the Heisenberg VOAs.
\begin{example}\label{ex4.22}
	A counter-example that shows $I\fusion{M^1}{M^2}{M^3}$ is not isomorphic to $(M^{3}(0)^{\ast}\otimes _{A(V)}A(M^{1})\otimes _{A(V)}M^{2}(0))^{\ast}$ was presented in Section 2 in \cite{L}. It was given as follows:
	
	Recall that the (universal) Virasoro VOA $M_{c}=M(c,0)/\<L(-1)v_{c,0}\>$ defined in \cite{FZ} has Zhu's algebra $A(M_{c})\cong \C [t]$, with $[\w]^{n}\mapsto t^{n}$. Let $M(c,h)$ be the Verma module of highest weight $h$ and central charge $c$ over the Virasoro algebra, then $M(c,h)$ is a module over $M_{c}$, and we have the following equalities held in $A(M(c,h))$:
	$$[b]\ast [\w]^{n}=[(L(-2)+L(-1))^{n}b],\qquad [\w]^{n}\ast [b]=[(L(-2)+2L(-1)+L(0))^{n}b],$$
	for all $b\in M(c,h)$ and $n\in \N$. Hence there is an identification of $\C[t]\cong A(M_{c})$-bimodules:
	\begin{equation}\label{4.30}
	\begin{aligned}
	\C[t_{1},t_{2}]&\cong A(M(c,h))\\
	f(t_{1},t_{2})&\mapsto f(L(-2)+2L(-1)+L(0), L(-2)+L(-1))v_{c,h},
	\end{aligned}
	\end{equation}
	where $C[t_{1},t_{2}]$ is a bimodule over $\C[t]$ on which the actions are given by:
	$$t^{n}.f(t_{1},t_{2})=t_{1}^{n}f(t_{1},t_{2}),\qquad f(t_{1},t_{2}).t^{n}=t_{2}^{n}f(t_{1},t_{2}).$$
	For $h_{1},h_{2}\in \C$ such that $M(c,h_{1})$ and $M(c,h_{2})$ are irreducible, it is proved (see (2.37) in \cite{L}) that $I\fusion {M(c,h_{1})}{M_{c}}{M(c,h_{2})}=0,$
	while $\dim (M(c,h_{2})(0)^{\ast}\otimes _{A(M_{c})}A(M(c,h_{1}))\otimes _{A(M_{c})}M_{c}(0))^{\ast}=1$.
	
	Although $M^2=M_c$ is neither a generalized Verma module nor irreducible, we can still use \eqref{n4.30} and \eqref{n4.32} to obtain the correct fusion rules. Indeed, since $M_c$ and $M(c,h_2)$ are both generalized by their bottom levels, by \eqref{n4.32}, we have: 
	\begin{equation}\label{n4.35}
	\dim I\fusion {M(c,h_{1})}{M_{c}}{M(c,h_{2})} \leq \dim (M(c,h_{2})(0)^{\ast}\otimes _{A(M_{c})}B_{h}(M(c,h_{1}))\otimes _{A(M_{c})}M_{c}(0))^{\ast}.
	\end{equation}
	Moreover, since $h=h_{1}+0-h_{2}$, it follows from Lemma \ref{lm4.4'} and Lemma \ref{lm4.5'} that 
	$$B_{h}(M(c,h_{1}))= A(M(c,h_{1}))/\mathrm {span}\{(L(-1)+L(0)-h_{2})[b]:b\in M(c,h_{1})\}.$$
	Then $[L(-1)b]=-[(\deg b+h_{1}-h_{2})b]$ in $B_{h}(M(c,h_{1}))$. It follows from \eqref{4.30} that 
	$$B_{h}(M(c,h_{1}))\cong \C[t_{0}],\quad \mathrm{with}\quad 
	[(L(-2)-L(0)+h_{2})^{n}v_{c,h_{1}}]\mapsto t_{0}^{n},$$
	and $\C[t_{0}]$ is a $\C[t](\cong A(M_{c}))$-bimodule on which the actions are given by:
	$$f(t_{0}).t^{n}=t_{0}^{n}f(t_{0}),\quad \mathrm{and}\quad t.f(t_{0})=(t_{0}+h_{2})^{n}f(t_{0}).$$
	Hence we have 
	$B_{h}(M(c,h_{1}))\otimes_{A(M_{c})} M_{c}(0)\cong \C[t_{0}]\otimes _{\C[t]} M_{c}(0)\cong M_{c}(0),$
	and so
	$$ (M(c,h_{2})(0)^{\ast}\otimes _{A(M_{c})}B_{h}(M(c,h_{1}))\otimes _{A(M_{c})}M_{c}(0))^{\ast}
	\cong \Hom_{A(M_{c})}( M_{c}(0),M(c,h_{2})(0))=0,$$
	since $o(\omega) v_{c,0}=0$, $o(\omega)v_{c,h_{2}}=h_{2} v_{c,h_{2}}$ and $h_{2}\neq 0$. Thus, $I\fusion {M(c,h_{1})}{M_{c}}{M(c,h_{2})}=0$ by \eqref{n4.35}.
\end{example}

We give another example that shows that the bimodule $B_{h}(M^{1})$ in \eqref{n4.30} cannot be replaced by the $A(V)$-bimodule $A_{0}(M^{1})$ defined in \cite{HY} either.

\begin{example}\label{ex4.23}
	Let $V=M_{\hat{\h}}(1,0)$ be the Heisenberg VOA of level $1$ associated to a one-dimensional vector space $\h=\C \al$ with $(\al| \al)=1$. By Theorem 3.1.1 in \cite{FZ}, one has
	$A(M_{\hat{\h}}(1,0))\cong \C[x],$ with $
	[\al(-i_{1}-1)...\al(-i_{n}-1)\vac]\mapsto (-1)^{i_{1}+...+i_{n}}x^{n}.$
	
	Let $\la\in \h$, we have a $V$-module $M_{\hat{\h}}(1,\la)=M_{\hat{\h}}(1,0)\otimes_{\C} \C e^{\la}$, with conformal weight $h=\frac{(\la|\la)}{2}$. Note that $M_{\hat{\h}}(1,\la)$ is the Verma module over the Heisenberg Lie algebra $\hat{\h}$. Since $M_{\hat{\h}}(1,\la)$ is irreducible, it is automatically a generalized Verma module associated with its bottom level $\C e^{\la}$. By Theorem 3.2.1 in \cite{FZ}, we have:$$
	A(M_{\hat{\h}}(1,\la))\cong \C e^{\la}\otimes_{\C} \C[x], \quad \mathrm{with}\quad 
	[\al(-i_{1}-1)...\al(-i_{n}-1)e^{\la}]\mapsto (-1)^{i_{1}+...+i_{n}}e^{\la}\otimes x^{n},
	$$
	where the bimodule actions are given by 
	$x.(e^{\la}\otimes x^{n})=e^{\la}\otimes x^{n+1}+(\la|\al) e^{\la}\otimes x^{n},$ and $(e^{\la}\otimes x^{n}).x=e^{\la}\otimes x^{n+1}$ for all $n\in \N$. 
	By definition in Section 4 of \cite{HY}, $$A_{0}(M_{\hat{\h}}(1,\la))=A(M_{\hat{\h}}(1,\la))/span\{[(L(-1)+L(0)-(\la|\la)/2)b]:b\in M_{\hat{\h}}(1,\la)\}.$$
	
	Choose $\la\in \h$ such that $(\la|\al)\neq0$. Recall that $\w=\frac{1}{2}\al(-1)^{2}\vac$, and so
	$$L(-1)e^{\la}=\Res_{z}Y_{W}(\w ,z)e^{\la}=\sum_{i\geq 0}\al(-1-i)\al(i)e^{\la}=(\la|\al)\al(-1)e^{\la}. $$ 
	Then we have $[(\la|\al)\al(-1)e^{\la}]=[L(-1)e^{\la}]=-[(L(0)-(\la|\la)/2)e^{\la}]=0$ in $A_{0}(M_{\hat{\h}}(1,\la))$, and $[\al(-1)e^{\la}]=0$ in $A_{0}(M_{\hat{\h}}(1,\la))$. For any spanning element $[\al(-i_{1}-1)...\al(-i_{n}-1)e^{\la}]$ of $A_{0}(M_{\hat{\h}}(1,\la))$, we then have
	$[\al(-i_{1}-1)...\al(-i_{n}-1)e^{\la}]=(-1)^{i_{1}+...+i_{n}}[\al(-1)^{n}e^{\la}]=0$
	for $n>0$. Thus, $A_{0}(M_{\hat{\h}}(1,\la))\cong \C [e^{\la}]$, with the module actions given by:
	\begin{equation}\label{4.31}
	x.[e^{\la}]=(\la|\al)[e^{\la}],\quad\mathrm{and}\quad  [e^{\la}].x=0.
	\end{equation}
	Now choose $\mu\in \h$ such that $(\mu|\al)\neq 0$, it is well-known that 
	$\dim I\fusion{M_{\hat{\h}}(1,\la)}{M_{\hat{\h}}(1,\mu)}{M_{\hat{\h}}(1,\la+\mu)}=1.$
	But
	$$A_{0}(M_{\hat{\h}}(1,\la))\otimes _{A(M_{\hat{\h}}(1,0))} M_{\hat{\h}}(1,\mu)(0)\cong \C [e^{\la}]\otimes _{\C[x]} \C e^{\mu}=0,$$
	since it follows from \eqref{4.31} that $[e^{\la}]\otimes e^{\mu}=\frac{1}{(\mu|\al)}[e^{\la}]\otimes o(\al(-1)\vac) e^{\mu}=\frac{1}{(\mu|\al)}[e^{\la}].x\otimes e^{\mu}=0$ in the tensor product above. Then we have:
	$$\dim (M_{\hat{\h}}(1,\la+\mu)(0)^{\ast}\otimes_{A(M_{\hat{\h}}(1,0))}A_{0}(M_{\hat{\h}}(1,\la))\otimes _{A(M_{\hat{\h}}(1,0))} M_{\hat{\h}}(1,\mu)(0))^{\ast}=0\neq 1.$$
	This shows that the isomorphism \eqref{n4.30} is not true if one replaces $B_{h}(M^1)$ with $A_0(M^1)$.
	
	Now we verify \eqref{n4.30} in this case. Indeed, since $\h=\C \al$, then $(\la|\al)\neq 0$ and $(\mu|\al)\neq 0$ imply that $\la=m\al$ and $\mu=n\al$, with $m\neq 0$ and $n\neq 0$. Hence 
	$$h=\frac{(\la|\la)}{2}+\frac{(\mu|\mu)}{2}-\frac{(\la+\mu|\la+\mu)}{2}=-(\la|\mu)=-mn\neq 0.$$ 
	By definition \ref{df4.1'}, we have the following equality holds in $B_{h}(M_{\hat{\h}}(1,\la))$: 
	$$[(\la|\al)\al(-1)e^{\la}]=[L(-1)e^{\la}]=-[(L(0)-\frac{(\la|\la)}{2}+h)e^{\la}]=-(\la|\mu)[e^{\la}]$$
	Then for any spanning element $[\al(-i_{1}-1)...\al(-i_{n}-1)e^{\la}]$ of $B_{h}(M_{\hat{\h}}(1,\la))$, we have: 
	$$[\al(-i_{1}-1)...\al(-i_{n}-1)e^{\la}]=(-1)^{i_{1}+...+i_{n}}[\al(-1)^{n}e^{\la}]=(-1)^{i_{1}+...+i_{n}}\left(\frac{-(\la|\mu)}{(\la|\al)}\right)^{n} [e^{\la}].$$
	Thus $B_{h}(M_{\hat{\h}}(1,\la))= \C [e^{\la}]$, with the module actions given by 
	\begin{equation}\label{4.32}
	[e^{\la}].x=\frac{-(\la|\mu)}{(\la|\al)}[e^{\la}](\neq 0),\quad \mathrm{and}\quad x.[e^{\la}]=\frac{-(\la|\mu)}{(\la|\al)} [e^{\la}]+(\la|\al)[e^{\la}].
	\end{equation}
	Then by \eqref{4.32}, we have $B_{h}(M_{\hat{\h}}(1,\la))\otimes _{A(M_{\hat{\h}}(1,0))} M_{\hat{\h}}(1,\mu)(0)\cong \C [e^{\la}]\otimes _{\C[x]} \C e^{\mu}$
	is a one-dimensional vector space, with 
	$x.[e^\la]\otimes e^{\mu} =[e^\la].x\otimes e^{\mu}+(\la|\al) [e^{\la}]\otimes e^{\mu}=(\la+\mu|\al) [e^{\la}]\otimes e^{\mu}.$
	On the other hand, $x. e^{\la+\mu}= (\la+\mu|\al) e^{\la+\mu}$. Thus we have:
	$$\dim \Hom_{A(M_{\hat{\h}}(1,0))}(B_{h}(M_{\hat{\h}}(1,\la))\otimes _{A(M_{\hat{\h}}(1,0))} M_{\hat{\h}}(1,\mu)(0), M_{\hat{\h}}(1,\la+\mu)(0))=1.$$
	This shows \eqref{n4.30} is true for $M^{1}=M_{\hat{\h}}(1,\la)$, $M^{2}=M_{\hat{\h}}(1,\mu)$, and $M^{3}=M_{\hat{\h}}(1,\la+\mu)$. 
	
	Furthermore, the argument above also shows that $B_{h}(M_{\hat{\h}}(1,\la))\otimes _{A(M_{\hat{\h}}(1,0))} M_{\hat{\h}}(1,\mu)(0)$ is a one-dimensional vector space spanned by an eigenvector of $\h$ of eigenfunction $(\la+\mu|\cdot)$. Hence we have: 
	$$\Hom_{A(M_{\hat{\h}}(1,0))}(B_{h}(M_{\hat{\h}}(1,\la))\otimes _{A(M_{\hat{\h}}(1,0))} M_{\hat{\h}}(1,\mu)(0), M_{\hat{\h}}(1,\gamma)(0))=0,$$
	if $\gamma\neq \la+\mu$. On the other hand, for $\gamma \neq \la+\mu$, it is well-known that $ I\fusion{M_{\hat{\h}}(1,\la)}{M_{\hat{\h}}(1,\mu)}{M_{\hat{\h}}(1,\gamma)}=0.$
	Thus, the rank one Heisenberg VOA verifies \eqref{n4.30}.
\end{example}

Although the bimodule $B_{h}(M^{1})$ by its construction is a quotient module of $A(M^{1})$, the vectors spaces $M^{3}(0)^{\ast}\otimes_{A(V)} B_{h}(M^{1})\otimes_{A(V)}M^{2}(0),$ and $M^{3}(0)^{\ast}\otimes_{A(V)} A(M^{1})\otimes_{A(V)}M^{2}(0)$
might be isomorphic to each other, it is easy to see that the case of the rank one Heisenberg VOA in Example \ref{ex4.23} above is such an example.

\begin{remark}Note that in $A(M^{1})$ we have: 
	$[\w]\ast[u]-[u]\ast [\w]=\Res_{z} Y_{M^{1}}(\w,z)u(1+z)^{\wt \w-1}=[L(-1)u+L(0)u], $
	for all $u\in M^{1}$. Hence 
	$[(L(-1)+L(0)+h_{2}-h_{3})u]=[\w]\ast [u]-[u]\ast [\w]+(h_{2}-h_{3})[u],$
	and by Lemma \ref{lm4.5'}, we have $B_{h}(M^{1})=A(M^{1})/J$, where 
	$$J=span\{ [\w]\ast [u]-[u]\ast [\w]+(h_{2}-h_{3})[u]:u\in M^{1}\}.$$
	We have $M^{3}(0)^{\ast}\otimes J\otimes M^{2}(0)=0$ in $M^{3}(0)^{\ast}\otimes_{A(V)} A(M^{1})\otimes_{A(V)}M^{2}(0)$. Indeed, for any $v_{3}'\in M^{3}(0)^{\ast}$ and $v_{2}\in M^{2}(0)$, 
	\begin{align*}
	&v_{3}'\otimes ([\w]\ast [u]-[u]\ast [\w]+(h_{2}-h_{3})[u])\otimes v_{2}\\
	&=v_{3}'(o(\w)-h_{3})\otimes [u]\otimes v_{2}-v_{3}'\otimes [u]\otimes (o(\w)-h_{2})v_{2}\\
	&=v_{3}'(L(0)-h_{3})\otimes [u]\otimes v_{2}-v_{3}'\otimes [u]\otimes (L(0)-h_{2})v_{2}\\
	&=0.
	\end{align*}
	However, in general we do not have $M^{3}(0)^{\ast}\otimes_{A(V)} (A(M^{1})/J)\otimes _{A(V)} M^{2}(0)$ isomorphic to $M^{3}(0)^{\ast}\otimes_{A(V)}A(M^{1})\otimes _{A(V)}M^{2}(0)/(M^{3}(0)^{\ast}\otimes J\otimes M^{2}(0)),$ see Example \ref{ex4.22}.
\end{remark}

\section{Acknowledgements} I want to thank Professor Haisheng Li for the conversations and for pointing out an error in an earlier version of this paper. I am also grateful to the anonymous referee for valuable comments and suggestions.

\end{document}